\documentclass[11pt]{article}
\usepackage{amsfonts,amsmath}

\def\transpose{{\hbox{\tiny\it T}}}

\newcommand{\home}{\mbox{$\!\mathtt{\sim}$}}

\newcommand{\RL}{{\mathbb R}}

\newcommand{\IN}{{\mathbb Z}}

\newcommand{\IND}{{\mathbb I}}

\newcommand{\Fmax}{\mbox{$F_{\rm max}$}}

\newcommand{\LA}{\Lambda}

\newcommand{\la}{\lambda}

\newcommand{\smalloneovern}{{\textstyle\frac{1}{n}}}

\def\be{\begin{eqnarray}}
\def\ee{\end{eqnarray}}
\def\ben{\begin{eqnarray*}}
\def\een{\end{eqnarray*}}

\def\breakMed{\\[.25cm]}
\def\breakBig{\\[.5cm]}
\def\st{\hbox{s.t.}}

\def\cH{\check{H}}
\def\upp{+}

\def\LWotwo{L_{\infty,2}^{W_0}}
\def\MWotwo{\clM_{1,2}^{W_0}}

\def\haxi{\hat\xi}

\def\tlabel#1{\label{t:#1}}
\def\slabel#1{\label{s:#1}}

\def\elabel#1{\label{e:#1}}

\def\sq{$\Box$}

\def\qed{\ifmmode\sq\else{\unskip\nobreak\hfil
\penalty50\hskip1em\null\nobreak\hfil\sq
\parfillskip=0pt\finalhyphendemerits=0\endgraf}\fi\par\medbreak}

\newsavebox{\junk}
\savebox{\junk}[1.6mm]{\hbox{$|\!|\!|$}}
\def\lll{{\usebox{\junk}}}

\def\det{{\mathop{\rm det}}}
\def\limsup{\mathop{\rm lim\ sup}}
\def\liminf{\mathop{\rm lim\ inf}}

\def\state{{\sf X}}

\newcommand{\field}[1]{\mathbb{#1}}

\def\Re{\field{R}} 
\def\TT{\field{T}}

\def\nat{\field{Z}_+}

\def\Co{\field{C}}

\def\ind{\field{I}}

\def\One{\hbox{\large\bf 1}}

\def\cExpect{\check{\Expect}}

\def\cProb{\check{\Prob}}

\def\cpi{\check{\pi}}

\def\cP{{\check{P}}}

\newcommand{\cf}{\check f}
\newcommand{\cg}{\check g}

\newcommand{\cF}{\check F}

\newcommand{\cmu}{{\check{\mu}}}

\newcommand{\cV}{{\check{V}}}

\newcommand{\haF}{{\widehat{F}}} 
\newcommand{\haM}{{\widehat{M}}}
\newcommand{\haQ}{{\widehat{Q}}}

\newcommand{\haK}{{\widehat{K}}}

\newcommand{\haR}{{\widehat{R}}}
\newcommand{\haG}{{\widehat{G}}}

\newcommand{\one}{\hbox{\rm\large\textbf{1}}}

\def\bfm{{\bf m}}

\def\bfmW{{\mbox{\protect\boldmath$W$}}}

\def\bfPhi{\mbox{\protect\boldmath$\Phi$}}
\def\bfPsi{\mbox{\protect\boldmath$\Psi$}}

\def\haP{{\widehat P}}

\def\halambda{{\hat\lambda}}
\def\hagamma{{\hat\gamma}}

\def\halambda{{\hat\lambda}}

\def\til={{\widetilde =}}

\def\clA{{\cal A}}
\def\clB{{\cal B}}

\def\clF{{\cal F}}
\def\clG{{\cal G}}
\def\clH{{\cal H}}

\def\clM{{\cal M}}

\def\clQ{{\cal Q}}

\def\clS{{\cal S}}

\def\clW{{\cal W}}
\def\clX{{\cal X}}

\def\half{{\mathchoice{\textstyle \frac{1}{2}}%
{\frac{1}{2}}%
{\hbox{\tiny $\frac{1}{2}$}}%
{\hbox{\tiny $\frac{1}{2}$}} }}

\def\eqdef{\mathbin{:=}}

\def\Prob{{\sf P}}
\def\Probsub{{\sf P\! }}
\def\Expect{{\sf E}}

 \def\eq#1/{(\ref{#1})}

\def\epsy{\varepsilon}

\def\varble{\,\cdot\,}

\newtheorem{theorem}{Theorem}[section]
\newtheorem{corollary}[theorem]{Corollary}
\newtheorem{proposition}[theorem]{Proposition}
\newtheorem{lemma}[theorem]{Lemma}

\def\Lemma#1{Lemma~\ref{t:#1}}
\def\Proposition#1{Proposition~\ref{t:#1}}
\def\Theorem#1{Theorem~\ref{t:#1}}

\def\Section#1{Section~\ref{s:#1}}

\def\eq#1/{(\ref{e:#1})}

\newcommand{\beqn}[1]{\notes{#1}%
\begin{eqnarray} \elabel{#1}}

\newcommand{\eeqn}{\end{eqnarray} }

\newcommand{\beq}[1]{\notes{#1}%
\begin{equation}\elabel{#1}}

\newcommand{\eeq}{\end{equation}} 

\def\bdes{\begin{description}}
\def\edes{\end{description}}

\newcommand{\oo}{\overline}
 
\def\baromega{\oo\omega}

\def\bara{{\oo {a}}}

\def\baralpha{{\oo{\alpha}}}

\def\baralpha{{\oo{\alpha}}}

\def\proof{\noindent{\sc Proof. }} 

\def\notes#1{}

\def\mmid{\, \| \,}

\def\L{L_\infty}
\def\LV{L_\infty^V}
\def\Lv{L_\infty^v}

\def\LW{L_\infty^{W}}
\def\LWo{L_\infty^{W_0}}

\def\MW{\clM_1^W}
\def\MWo{\clM_1^{W_0}}

\def\Mv{\clM_1^v}

\def\llangle{\langle\!\langle}
\def\rrangle{\rangle\!\rangle}

\def\haK{{\widehat K}}

\def\cv{{\check{v}}}

\def\cLambda{{\check{\Lambda}}}
\def\cXi{{\check{\Xi}}}

\def\tilM{\widetilde{M}}

\def\bfmX{{\mbox{\protect\boldmath$X$}}}

\def\transpose{{\hbox{\tiny\it T}}}

\newlength{\noteWidth}
\long\def\notes#1{\ifinner
             {\tiny #1}
             \else
             \marginpar{\parbox[t]{\noteWidth}{\raggedright\tiny #1}}
             \fi}

\def\archival#1{} 
\def\done#1{}  

\newcounter{tasks}

\newenvironment{MyEnumerate}%
{\begin{quote}%
\begin{list}{\hspace{-.75cm}\hbox{\textrm{\rm (\roman{tasks})}\ }}{\usecounter{tasks}%
        \setlength{\labelsep}{0pt}
        \setlength{\leftmargin}{0pt}
        \setlength{\rightmargin}{0pt}
        \setlength{\labelwidth}{0pt}
        \setlength{\listparindent}{0pt}}}%
{\end{list}\end{quote}}

\newcounter{tasksA}

{\begin{quote}%
\begin{list}{\hspace{-.75cm}\hbox{\textrm{\rm (\alph{tasksA})}\ }}{\usecounter{tasksA}%
        \setlength{\labelsep}{0pt}
        \setlength{\leftmargin}{0pt}
        \setlength{\rightmargin}{0pt}
        \setlength{\labelwidth}{0pt}
        \setlength{\listparindent}{0pt}}}%
{\end{list}\end{quote}}

\begin{document}
 
\title{\large\bf
Large Deviations Asymptotics and the Spectral Theory\\
of Multiplicatively Regular Markov Processes}

\author
{
	{\bf\normalsize I. Kontoyiannis}\\
		\small
		Division of Applied Mathematics
		and Dept of Computer Science\\
		\small      
		Brown University, Box F,
		182 George St., Providence, RI 02912, USA\\
		\small      
                Email: {\tt yiannis@dam.brown.edu }
                Web: {\tt www.dam.brown.edu/people/yiannis/}
		  \thanks{I.K. was supported in part
		  by a Sloan Foundation Research Fellowship 
		and by NSF grant \#0073378-CCR.}
\and
        {\bf\small      S.P. Meyn}\\
		\small     
		Department of Electrical and Computer Engineering\\
		\small     
                University of Illinois at Urbana-Champaign,
                Urbana, IL 61801, USA\\
		\small     
		 Email: {\tt meyn@uiuc.edu }
		 Web: {\tt black.csl.uiuc.edu/$\home$meyn/}
		\thanks{S.M. was spported in part by
                NSF grant ECS 99-72957.}
}

\date{}
\maketitle

\begin{abstract}
In this paper we continue the investigation
of the spectral theory and exponential asymptotics
of primarily discrete-time Markov processes,
following Kontoyiannis and Meyn \cite{kontoyiannis-meyn:I}.
We introduce a new family of nonlinear Lyapunov 
drift criteria, which characterize distinct 
subclasses of geometrically ergodic Markov processes
in terms of simple inequalities for the nonlinear
generator. We concentrate primarily on the 
class of {\em multiplicatively
regular} Markov processes, which are characterized
via simple conditions
similar to (but weaker than) those of Donsker-Varadhan.
For any such process 
$\bfPhi=\{\Phi(t)\}$ with transition kernel $P$
on a general state 
space $\state$, 
the following are obtained.
\begin{description}
\item \textit{Spectral Theory}:
For a large class of (possibly unbounded)
functionals $F:\state\to\Co$,
the kernel $\haP(x,dy) = e^{ F(x)}P(x,dy)$
has a discrete spectrum in 
an appropriately defined Banach space.
It follows that there
exists a ``maximal'' solution $(\la,\cf)$
to the {\em multiplicative
Poisson equation}, defined
as the eigenvalue problem $\haP\cf=\la\cf$.
The functional $\Lambda(F) = \log(\lambda)$ 
is convex, smooth,
and its convex dual $\Lambda^{*}$
is convex, with compact sublevel sets.

\item \textit{Multiplicative Mean Ergodic Theorem}:
Consider the partial sums $\{S_t\}$ of the process
with respect to any one of the functionals
$F(\Phi(t))$ considered above.
The normalized mean $\Expect_x[\exp(S_t)]$ 
(and not the logarithm of the mean) 
converges to $\cf(x)$
exponentially fast, where $\cf$
is the above solution of the multiplicative
Poisson equation.
  
\item \textit{Multiplicative regularity}:
The Lyapunov drift criterion under which 
our results are derived is equivalent
to the existence of regeneration times
with finite exponential moments for the
partial sums $\{S_t\}$, with respect to any functional
$F$ in the above class.

\item \textit{Large Deviations}:
The sequence of empirical measures of $\{\Phi(t)\}$ satisfies 
a large deviations principle
in the ``$\tau^{W_0}$-topology,''  
a topology finer that the usual $\tau$-topology,
generated by the above class of functionals
$F$ on $\state$ which is strictly larger
than $L_\infty(\state)$.
The rate function of this LDP is $\Lambda^{*}$, 
and it is shown to coincide with the
Donsker-Varadhan rate function
in terms of relative entropy.

\item \textit{Exact Large Deviations Asymptotics}:
The above partial sums $\{S_t\}$ are shown to
satisfy an exact large deviations
expansion, analogous to that obtained by
Bahadur and Ranga Rao for independent random 
variables.  
\end{description}

\bigskip


\bigskip


\end{abstract}

\newpage 
\tableofcontents

\newpage

\section{Introduction and Main Results}
\slabel{intro}

Let $\bfPhi=\{ \Phi(t): t\in\TT\}$ be a Markov processes
taking values in a Polish state space $\state$, 
equipped with its associated Borel $\sigma$-field $\clB$.
The time index $\TT$ may be discrete, 
$\TT=\IN_+$, or continuous 
$\TT=\RL_+$, but we specialize 
to the discrete-parameter case
after \Section{metheory}.

The distribution of $\bfPhi$ is determined 
by its initial state $\Phi(0)=x\in\state$, and 
the transition semigroup
$\{P^t\,:\,t\in\TT\}$, where in discrete time 
all kernels $P^t$ are powers of the
1-step transition kernel $P$.
Throughout the paper we assume that 
$\bfPhi$ is  \textit{$\psi$-irreducible}
and \textit{aperiodic}.   This means that 
there is a $\sigma$-finite measure $\psi$
on $(\state,\clB)$ such that,
for any $A\in\clB$ satisfying $\psi(A)>0$
and any initial condition $x$,
\[
P^t(x,A) > 0,\qquad \hbox{for all $t$ sufficiently large.}
\]
Moreover, we assume that $\psi$ is \textit{maximal} 
in the sense that any other such $\psi'$ 
is absolutely continuous with respect to $\psi$ 
(written $\psi'\prec\psi$).

For a $\psi$-irreducible Markov process it is known  that 
ergodicity is equivalent to the existence of a solution 
to the \textit{Lyapunov drift criterion} (V3) below 
\cite{meyn-tweedie:book,dowmeytwe95a}.
Let $V:\, \state \to (0, \infty]$ be an extended-real valued 
function, with $V(x_0) < \infty$ for at least 
one $x_0\in\state$, and write
$\clA$ for the (extended) generator of the semigroup 
$\{P^t: t\in\TT\}$.  This is equal to $\clA=(P-I)$ in discrete time
(where $I=I(x,dy)$ denotes the identity kernel $\delta_x(dy)$),
and in continuous-time we think of $\clA$ as a generalization 
of the classical differential generator $\clA=\frac{d}{dt}P^t|_{t=0}$.

Recall that a function $s\colon\state\to \Re_+$ 
and a probability measure $\nu$ on $(\state,\clB)$
are called \textit{small} if for some measure
$m$ on $\IN$ with finite mean we have
\[
\sum_{t\geq 0}P^t(x,A) \, m(t) \ge s(x)\nu(A),\qquad x\in\state,\ A\in\clB.
\]
A set $C$ is called {\em small} if  
$s=\epsilon\IND_C$ is a small function for some $\epsilon>0$.
Also recall that an arbitrary kernel $\haP=\haP(x,dy)$ 
acts linearly on functions $f:\state\to\Co$ and 
measures $\nu$ on $(\state,\clB)$, via
\be
\haP f\, (\varble)=\int_{\state}\haP(\cdot,dy)f(y)
\;\;\;\mbox{and}\;\;\;
\nu \haP\, (\varble)=\int_{\state}\nu(dx)\haP(x,\varble),
\;\;\;\mbox{respectively.}
\label{eq:act}
\ee

We say that {\em the Lyapunov drift condition {\em (V3)} holds
with respect to the Lyapunov function $V$} 
\cite{meyn-tweedie:book}, if:
$$
\left. 
\mbox{\parbox{.85\hsize}{\raggedright
For a function  $W\colon\state\to [1,\infty)$, 
                a small set $C\subset\state$,
                and constants $\delta>0$, $b<\infty$,
\[
\clA V\leq -\delta W+b\ind_C\, ,\qquad 
\hbox{on $ S_V:= \{x\,:\, V(x) < \infty\}.$}
\]
}}
\right\}
\eqno{\bf (V3)}
$$
Condition (V3) implies that the set 
$S_V$ is absorbing (and hence full), 
so that $V(x) < \infty$ a.e.\ $[\psi]$; see
\cite[Proposition 4.2.3]{meyn-tweedie:book}.

As in \cite{meyn-tweedie:book,kontoyiannis-meyn:I}, 
a central role in our development will
be played by weighted $L_\infty$ spaces:
For any function $W\colon\state\to(0,\infty]$,
define the Banach space of complex-valued functions, 
\be
   \LW  \eqdef \Big \{ g\colon\state\to\Co\, \ \st\ \sup_x 
	\frac{|g(x)|}{W(x)}<\infty\Big\}\, ,
\label{eq:weightedL}
\ee
with associated norm $\|g\|_W \eqdef \sup_x |g(x)|/W(x)$.
We write $\clB^+$ for the set of functions
$s: \state \to [0,\infty]$ satisfying 
$\psi(s) \eqdef\int s(x)\, \psi(dx) > 0$,
and, with a slight abuse of notation, we
write  $A\in\clB^+$ if $A\in\clB$ and $\psi(A)>0$ (i.e.,
the indicator function  
$\IND_A$ is in $\clB^+$).
Also, we let $\MW$ denote the Banach space of signed 
and possibly complex-valued measures 
$\mu$ on $(\state,\clB)$ satisfying  
$\|\mu\|_W \eqdef \sup_{F\in\LW} |\mu|(F)<\infty$. 

The following consequences of~(V3) may be found in  
\cite[Theorem 14.0.1]{meyn-tweedie:book}.

\begin{theorem}
\tlabel{Fergodic}
{\em (Ergodicity)}
Suppose that $\bfPhi$ is a $\psi$-irreducible and aperiodic 
discrete-time chain, 
and that condition~{\em (V3)} is satisfied.
Then the following properties hold:
\begin{enumerate}
\item
{\em ($W$-ergodicity)\ }
The process is positive recurrent with a
unique invariant probability
measure $\pi\in\MW$ and for all $x\in S_V$,
\ben
\sup_{F\in L_\infty^W} \Big| P^t(x,F) - \pi(F) \Big |
&\to& 
        0, 
        \hspace{1.6in} t\to\infty,\\
\frac{1}{T} \sum_{t=0}^{T-1} F(\Phi(t))
&\to&  
        \pi(F)\eqdef\int F(y)\,\pi(dy),
        \hspace{0.2in}
        T\to\infty,\ a.s.\ [\Prob_x]
                        \quad F\in\LW\,,
\een
where $\Prob_x$ denotes
the conditional distribution of $\bfPhi$
given $\Phi(0)=x$. 
\item 
{\em ($W$-regularity)\ }
For any $A\in\clB^+$
there exists $c=c(A)<\infty$ such that 
\[
\Expect_x 
\Bigl[ \sum_{t=0}^{\tau_A - 1} W(\Phi(t)) \Bigr] 
	   	\le     \delta^{-1}V(x) + c,
                        \qquad x\in\state.
\]
where
$\Expect_x$ is the expectation with respect to
$\Prob_x$, and the hitting times $\tau_A$ are
defined as,
\begin{equation}
\tau_A \eqdef \inf \{ t \ge 1\,:\, \Phi (t) \in A\},\qquad A\in\clB.
\elabel{tauA}
\end{equation}

\item
{\em (Fundamental Kernel)\ }
There exists a linear operator $Z\colon \LW\to L_\infty^{V+1}$,
the {\em fundamental kernel}, such that
\[
\clA Z F = - F + \pi(F),\qquad F\in   \LW\, .
\]
That is, the function $\haF\eqdef ZF$ solves the 
{\em Poisson equation}, $\,\clA\haF = - F + \pi(F)\,$.
\end{enumerate}
\end{theorem}

\subsection{Multiplicative Ergodic Theory}
\slabel{metheory}

The ergodic theory outlined in \Theorem{Fergodic}
is based upon consideration of the semigroup 
of linear operators $\{P^t\}$ acting  on the Banach space $\LW$.
In particular, the ergodic behavior of the corresponding Markov
process can be determined via the generator $\clA$ 
of this semigroup. In this paper we show that the
foundations of the {\em multiplicative}
ergodic theory and of the large deviations 
behavior of $\bfPhi$ can be developed in 
analogy to the linear theory, by shifting
attention from the semigroup of linear operators
$\{P^t\}$ to the family of nonlinear, 
convex operators $\{ \clW^t\}$ defined, for
appropriate $G$, by
\[
\clW^t G\, (x)  :=  \log\Bigl(\Expect_x[e^{G(\Phi(t))}] \Bigr),
                           \qquad x\in\state\, , \ t\in\TT\,.
\]

Formally,  we would like to define
the `generator' $\clH$ associated with $\{ \clW^t\}$
by letting $\clH=(\clW-I)$ in discrete time and 
$\clH=\frac{d}{dt}\clW^t|_{t=0}$ in continuous time.
Observing that $\clW^t G=\log(P^te^G)$, in
discrete time we have 
$$\clH G=(\clW-I)G=\log (Pe^G)-G=\log(e^{-G}Pe^G),$$
and in continuous time we can similarly calculate,
\ben
\clH G
\;=\;
\lim_{t\to 0}\frac{1}{t}\bigl[\clW^t-I\bigr]G
\;=\;
\lim_{t\to 0}\frac{1}{t}\log\bigl(e^{-G}P^te^G\bigr)
\;=\;
e^{-G}\clA e^G\,,
\een
whenever all the above limits exist.
Rather than assume differentiability, we use these expressions
as motivation for the following rigorous definition
of the {\em nonlinear generator},  
\begin{equation}
  \clH(G) = \begin{cases} 
                   \
     \log(e^{-G}Pe^G)\quad & \hbox{discrete time ($\TT=\nat$)};
        \breakMed
			\
                     \phantom{\log(}   e^{-G} \clA e^G 
			   & \hbox{continuous time ($\TT=\Re_+$),}
            \end{cases}
\elabel{clH}
\end{equation}

\vspace{-0.1in}

\noindent
when $e^{G}$ is in the domain of the extended 
generator.
In continuous time, this is Fleming's nonlinear 
generator;
see \cite{fle78a} for a starting point,
and \cite{fen99a,feng-kurtz:preprint}
for recent surveys.

In this paper our main focus will be on the following 
`multiplicative' analog of (V3), where the role of the
generator is now played by the nonlinear generator $\clH$.
We say that the {\em Lyapunov drift criterion} (DV3) 
{\em holds with respect
to the Lyapunov function $V:\state\to(0,\infty]$}, if:
$$
\left. 
\mbox{\parbox{.85\hsize}{\raggedright
For a function  $W\colon\state\to [1,\infty)$, 
                a small set $C\subset\state$,
                and constants $\delta>0$, $b<\infty$,
\[
 \clH(V)\leq -\delta W+b\ind_C\,  ,\qquad \hbox{on $ S_V$}.
\]
}}
\right\}
\eqno{\hbox{\bf (DV3)}}
$$

\noindent
[This condition was introduced 
in \cite{kontoyiannis-meyn:I}, under the
name (mV3).]
Under either condition~(V3) or (DV3),
we let $\{C_W(r)\}$ denote the 
sublevel sets of $W$:
\begin{equation}
C_W(r) =\{ y : W(y) \le r\},\qquad r\in\RL.
\elabel{Cn}
\end{equation}
The main assumption in many of our results below will 
be that $\bfPhi$ satisfies (DV3), and also that 
the transition kernels satisfy a mild continuity
condition: We require that they
possess a density with respect
to some reference measure, uniformly over all initial
conditions $x$ in the sublevel set $C_W(r)$ of $W$.
These assumptions are formalized in 
condition (DV3\upp) below. 

%


$$
\left. 
\mbox{\parbox{.85\hsize}{\raggedright\begin{description}
\item[(i)]
The Markov process $\bfPhi$ is $\psi$-irreducible,
        aperiodic, and it satisfies condition~(DV3) with some
	Lyapunov function $V: \state \to [1,\infty)$;

\item[(ii)] 
There exists $T_0>0$ such that, for each $r< \|W\|_\infty$,
        there is a measure $\beta_r$ 
	with $\beta_r(e^V)<\infty$ and 
        $\Prob_{x}\{ \Phi(T_{0})\in A ,\ 
	\tau_{C^c_{W}(r)} >T_{0}\} \le \beta_r(A)$ for  
      all $  x\in C_{W}(r) $, $A\in\clB$.
      \end{description}
}}\right\}
\eqno{\hbox{\bf (DV3\upp)}}
$$



Condition~(DV3\upp) 
captures the essential ingredients of the large
deviations conditions imposed by 
Donsker and Varadhan in their pioneering work
\cite{donsker-varadhan:I-II,donsker-varadhan:III,donsker-varadhan:IV},
and is in fact somewhat weaker than those conditions.
In \Section{comparison} an extensive 
discussion of this assumption is given,
its relation to several well-known conditions
in the literature is described in detail.
In particular, part (ii)
of condition (DV3+) [to which we will often
refer as the ``density assumption'' in (DV3+)]
is generally the weaker of the two assumptions. 

In most of our results we assume 
that the function $W$ in (DV3) is unbounded,
$\|W\|_{\infty}:=\sup_x |W(x)|=\infty$.
When this is the case, we let
$W_0:\state\to[1,\infty)$ be a fixed
function in $L_\infty^W$, whose growth 
at infinity is strictly slower than $W$
in the sense that
\begin{equation}
\lim_{r\to\infty} \sup_{x\in\state}\,
        \Bigl [
        \frac{W_0(x)}{W(x)}\,\ind_{\{W(x)>r\}}
        \Bigr ]
= 0\,.
\elabel{Wo}
\end{equation}

Below we collect, from various parts of the
paper, the ``multiplicative'' ergodic 
results we derive from (DV3+), in analogy to
the ``linear'' ergodic-theoretic results stated 
in \Theorem{Fergodic}.

\begin{theorem}
\tlabel{Multergodic}
{\em (Multiplicative Ergodicity)}
Suppose that the discrete-time chain
$\bfPhi$ satisfies condition {\em (DV3\upp)} with $W$ unbounded,
and let $W_0\in\LW$ be as in {\em \eq Wo/}.
Then the following properties hold:
\begin{enumerate}
\item
{\em ($W$-multiplicative ergodicity)\ }
The process is positive recurrent with a unique invariant 
probability measure $\pi$ satisfying, for some
$\eta>0$,
\[
\pi(e^{\eta V})<\infty\
 \hbox{and}\quad 
\pi(e^{\eta W})<\infty.
\]
For any real-valued $F\in\LWo$,
there exist $\cF\in\LV$,
$\Lambda(F)\in\Co$, and
constants 
$b_0>0$, $B_0<\infty$,
such that
\begin{equation}
\left| \Expect_x 
        \left[ \exp 
                \Bigl(
                        \sum_{t=0}^{T-1}[F(\Phi(t))] - \LA(F)]
                \Bigr)
        \right]
                - e^{\cF(x)} 
\right|
\leq   e^{\eta V(x)+B_0 -b_0 T}
\,,
\elabel{f-muli-norm}
\end{equation}
for all $T\geq 1,\; x\in\state\,.$

\item 
{\em ($W$-multiplicative regularity)\ }
For any $A\in\clB^+$
there exist constants $\eta =\eta(A)>0$ and $c=c(A)<\infty$,
such that 
\[
\log\Bigl(
\Expect_x 
\Bigl[ \exp
\Bigl( 
	\eta\sum_{t=0}^{\tau_A-1} W(\Phi(t)) \Bigr)\Bigr] 
\Bigr)
                \le V(x) + c,
                        \qquad x\in\state.
\]

\item
{\em (Multiplicative Fundamental `Kernel')\ }
There exists a nonlinear operator
$\clG\colon \LWo \to \LV ,$
the {\em multiplicative fundamental kernel},
such that the function $\cF$ in (1.) can be expressed
as $\cF=\clG(F)$ for real-valued $F\in\LWo$, and $\cF$ solves   
the {\em multiplicative Poisson equation},
\be
\clH(\cF) = -F+\Lambda(F) \, .
\label{eq:mpe}
\ee
\end{enumerate}
\end{theorem}

\proof
Assumption (DV3) combined with \Theorem{v5} implies
that $\bfPhi$ is geometrically ergodic 
(equivalently, $V_0$-uniformly ergodic) 
for some Lyapunov function
$V_{0}\colon\state\to[1,\infty)$, hence the process 
is also positive recurrent.  
Moreover,   $v_{\eta_{0}}\eqdef e^{\eta_{0} V}\in L_{\infty}^{V_0}$
for some $0<\eta_{0}<1$.
By the geometric ergodic theorem of \cite{meyn-tweedie:book}
it follows that $\pi(v_{\eta_{0}})<\infty$.
 
Under (DV3), 
the stochastic process $\bfm=\{m(t)\}$ defined below
is a super-martingale with respect to 
$\clF_t = \sigma \{\Phi(s): 0\le s\le t\},\;t\geq 0$,
\begin{equation}
m(t) \eqdef  \exp\Bigl(V(\Phi(t)) + 
\sum_{s=0}^{t-1} 
[\delta W(\Phi(s)) - b\ind_C(\Phi(s))]\Bigr),\quad t \ge 0\,.
\elabel{expMart}
\end{equation}
From the super-martingale property and Jensen's 
inequality we obtain the bound,
\[
\Expect_x
 \Bigl[ \exp 
                \Bigl(\eta_{0} V(\Phi(t)) - \eta_{0}b +
	\sum_{s=0}^{t-1} \eta_{0} \delta W(\Phi(s)) \Bigr) \Bigr] 
	<  v_{\eta_{0}}(x)\, ,\qquad x\in\state.
\]
which gives the desired bound in (1.),
where $\eta\eqdef \delta\eta_{0}$.
The multiplicative ergodic limit  \eq f-muli-norm/ follows from 
\Theorem{mainSpectral}~(iii).   The existence of an inverse $\clG$ to $\clH$ 
is given in \Proposition{discreteTimeGa},  
which establishes the bound $\cF\in\LV$ stated in (1.), as well as result~(3.).
 
\Theorem{multReg}
shows that (DV3) actually characterizes $W$-multiplicative 
regularity, and provides the bound in  (2.). 
\qed

As in \cite{kontoyiannis-meyn:I}, central to our
development is the observation that the multiplicative
Poisson equation (\ref{eq:mpe}) can be written as an
eigenvalue problem. In discrete-time with $\LA=\LA(F)$,  
(\ref{eq:mpe}) becomes $(e^FP)e^{\cF}=e^\Lambda e^{\cF}$,
or, writing $f=e^F, \cf=e^{\cF}$ and $\la=e^\Lambda$,
we obtain the eigenvalue equation,
$$
P_f \cf=\la \cf,\qquad \mbox{for the kernel  $P_f(x,dy):=f(x)P(x,dy)$.}
$$

The assumptions of \Theorem{Multergodic} are most easily  
illustrated in continuous time.  Consider the following 
diffusion model on $\Re$, sometimes referred to as the 
\textit{Smoluchowski equation}.
For a given potential  $u:\Re\to\Re_+$, 
this is defined by the stochastic differential equation
\begin{equation}
         d X(t) =-u_x(X(t)) \, dt 
                 + \sigma d W(t)\, ,
                 \elabel{SmoluSDE}
\end{equation}
where $u_x\eqdef  \frac{d}{dx} u$, and 
$\bfmW=\{W(t): t\geq 0\}$ is a standard Brownian motion.
On $C^2$, the extended generator $\clA$ 
of $\bfmX=\{X(t):t\geq 0\}$ coincides with 
the differential generator given by,
\begin{equation}
\clA = \half \sigma^2 \,\frac{d^2}{dx^2} - 
                  u_x\frac{d}{dx}\;.
\end{equation}
When $\sigma>0$ this is an \textit{elliptic} diffusion, so that the
semigroup $\{P^t\}$ has a family of smooth, positive densities
$P^t(x,dy)=p(x,y;t)dy$, $x,y\in\Re$ \cite{kun90}.  
Hence the Markov process $\bfmX$ is
$\psi$-irreducible, with $\psi$ equal to Lebesgue measure on $\Re$.

A special case is  the one-dimensional Ornstein-Uhlenbeck process,
\begin{equation}
d X(t) = -\delta X(t)\, dt + \sigma d W(t)\, ,
\elabel{OU}
\end{equation}
where the corresponding potential function 
is $u(x)=\half \delta x^2$,
$x\in\Re$.   

\begin{proposition}
\tlabel{Smoluchowski}
The Smoluchowski equation satisfies {\em (DV3+)}
with $V=1 + u\sigma^{-2}$ and $W=1+u_{x}^{2}$, 
provided the potential function $u\colon\Re\to\Re_{+}$ 
is $C^2$ and satisfies:
\begin{itemize}
\item[\rm(a)]
\quad
$\displaystyle
\lim_{|x|\to\infty }u(x) = \infty$;
\item[\rm(b)]
\quad
$\displaystyle
\lim_{|x|\to\infty }
\frac{(u_x(x))^2} {|u_{xx}(x)|} = \infty,\qquad \liminf_{|x|\to\infty } (u_x(x))^2   >0$.
\end{itemize}
\end{proposition}

\proof
Let $V= 1 + u\sigma^{-2}$.   We then have,
 \[\begin{array}{rcl}  \clH(V)\eqdef 
 e^{-V}\clA e^{V} 
 &=& e^{-V}\Bigl\{ -u_{x}\Bigl(e^{V} \sigma^{-2} u_{x}\Bigr)
+\half \sigma^{2}\Bigl(e^{V}[u_{xx}\sigma^{-2}+ 
\sigma^{-4}u_{x}^{2}]\Bigr)   \Bigr\}
\breakBig  &  = &   
-\half \sigma^{-2} u_{x}^{2} + \half u_{xx} \,.
\end{array}\]
It is thus clear that the desired drift conditions hold.
The proof is complete since $P^{t}(x,dy)$ possesses 
a continuous density $p(x,y;t)$ for each $t>0$:  
We may take $T_{0}=1$, and for each $r$
we   take $\beta_{r}$ equal
to a constant times Lebesgue measure on $C_{W}(r)$.
\qed


\Proposition{Smoluchowski} does not admit an exact 
generalization to discrete-time models.  
However, the discrete-time   one-dimensional Ornstein-Uhlenbeck process,
\begin{equation}
 X(t+1)-X(t) = -\delta X(t)  +   W(t+1)\, ,\qquad t\ge 0, \ X(0)\in\Re,
\elabel{OUdt}
\end{equation}
does satisfy the conclusions of the proposition, again with $V= 1 + \epsilon_{0}x^{2}$
for some $\epsilon_{0}>0$, when  $\delta>0$ and $\bfmW $ is an i.i.d.\
Gaussian process with positive variance.


\medskip

\noindent
{\bf Notation. } Often in the transition from ergodic 
results to their multiplicative counterparts we have to 
take exponentials of the corresponding quantities. 
In order to make this correspondence transparent 
we have tried throughout the paper to follow, 
as consistently as possible, the convention that 
the exponential version of a quantity is written as
the corresponding lower case letter. For example,
above we already had
$f=e^F, \cf=e^{\cF}$ and $\la=e^\Lambda$.

\subsection{Large Deviations}

From now on we restrict attention
to the discrete-time case.

Part~1 of~\Theorem{Multergodic} extends the 
{\em multiplicative mean ergodic theorem} of 
\cite{kontoyiannis-meyn:I} to the larger class
of (possibly unbounded)
functionals $F\in\LWo$. In this section
we assume that~(DV3\upp) holds 
with an unbounded function $W$,
and we let a function $W_0\in\LW$ 
be chosen as in \eq Wo/.

For $n\geq 1$, let $L_n$ denote the empirical
measures induced by $\bfPhi$ on $(\state,\clB)$,
\begin{equation}
L_n \eqdef 
\frac{1}{n} \sum_{t=0}^{n-1} \delta_{\Phi(t)}
\qquad  n\geq 1\, ,
\elabel{OccEmp}
\end{equation}
and write $\langle \cdot , \cdot\rangle$ for the usual
inner product; for $\mu$ a measure and $G$ a function,
$\langle \mu , G\rangle = \mu(G) \eqdef \int G(y)\,\mu(dy)$,
whenever the integral exists.
Then, from \Theorem{mainSpectral}
it follows that for any real-valued
$F\in \LWo $ and any $a\in\RL$
we have the following version of 
the multiplicative mean ergodic theorem,
\be
\exp \Bigl( - n\LA(a F ) \Bigr)
\Expect_x
        \Bigl[ \exp
                \Bigl(
                        a n\langle L_n, F \rangle
                \Bigr)
        \Bigr]
        \;\to\;
        \cf_a(x)
\,, \qquad
n\to\infty,\, x\in\state\,,
\label{eq:mmet}
\ee
where $\cf_a\eqdef e^{\clG(aF)}$ is the eigenfunction 
constructed in part~3 of~\Theorem{Multergodic},
corresponding to the function $a F$.

In~\Section{ldp}, strong large deviations results 
for the sequence of empirical measures $\{L_n\}$
are derived from the multiplicative mean ergodic 
theorem in (\ref{eq:mmet}), 
using standard techniques 
\cite{dawson-gartner:87,chaganty-sethuraman:93,dembo-zeitouni:book}.
First we show that, for any initial condition $x\in\state$,
the sequence $\{L_n\}$
satisfies a large deviations principle
(LDP) in the space $\clM_1$ of all 
probability measures on $(\state,\clB)$ 
equipped with the {\em $\tau^{W_0}$-topology},
that is, the topology generated by the system 
of neighborhoods 
\be
N_F(c,\epsilon)
\eqdef 
\bigl\{\nu\in\clM_1:|\nu(F)-c|<\epsilon\}\,,
\quad\mbox{for real-valued}\;F\in \LWo ,\,c\in\RL,\,\epsilon>0  \,.
\label{eq:nbhds}
\ee
Moreover, the rate function $I(\nu)$ that governs 
this LDP is the same as the Donsker-Varadhan 
rate function, and can be characterized in
terms of relative entropy,
$$
I(\nu):=\inf H(\nu\odot \cP\|\nu\odot P)\,,
$$
where the infimum is over all transition kernels $\cP$
for which $\nu$ is an invariant measure, $\nu\odot \cP$ denotes
the bivariate measure $[\nu\odot \cP](dx,dy)\eqdef \nu(dx)\cP(x,dy)$
on $(\state\times\state, \clB\times\clB)$,
and $H(\varble{\|}\varble)$ denotes the relative entropy,
\begin{equation}
\begin{array}{rcl}
H(\mu\|\nu)=
        \left\{ \begin{array}{ll}
                        \int d\mu\log\frac{d\mu}{d\nu},&\;\;\;\mbox{when}
                                \frac{d\mu}{d\nu}\;\mbox{exists}
\breakMed
                        \infty,                        &\;\;\;\mbox{otherwise.}
                \end{array}
        \right.
\end{array}
\elabel{entropyGeneral}
\end{equation}
[Throughout the paper we follow the usual convention that
the infimum of the empty set is $+\infty$.]
As we discuss in \Section{DV} and
\Section{ldp},  
the density assumption in (DV3+)~(ii) is
weaker than the continuity assumptions
of Donsker and Varadhan, but it cannot be
removed entirely.

Further, the precise convergence in (\ref{eq:mmet})
leads to exact large deviations
expansions analogous to those obtained by
Bahadur and Ranga Rao \cite{bahadur-rao:60} for
independent random variables, and to the 
local expansions established in 
\cite{kontoyiannis-meyn:I} for
geometrically ergodic chains.
For real-valued, non-lattice functionals $F\in\LWo$,
in~\Theorem{Bahadur-RaoNL}
we obtain the following:
For $c>\pi(F)$ and $x\in\state$,
\be
\Probsub_x\Bigl\{\sum_{t=0}^{n-1}F(\Phi(t)) \geq nc\Bigr\}
       \;\sim\;
        \frac{\cf_a(x)}{a\sqrt{2\pi n\sigma_a^2}}
        e^{-nJ(c)},
        \quad
        n\to\infty,
\label{eq:introBR}
\ee
where $a\in\RL$ is chosen such that $\frac{d}{da}\LA(aF)=c$,
$\cf_a(x)$ is the eigenfunction appearing in the 
multiplicative mean ergodic theorem (\ref{eq:mmet}),
$\sigma_a^2=\frac{d^2}{da^2}\LA(aF)$,
and the exponent $J(c)$ is given in terms of $I(\nu)$ as
\be
J(c)\eqdef \inf \bigl\{
 I(\nu) : \hbox{$\nu$ is a probability measure on $(\state,\clB)$
  satisfying $\nu(F)\geq c$}
 			\bigr\}\,.
			 \label{eq:rf}
\ee
A corresponding expansion is given for
lattice functionals.

These large deviations results extend the classical
Donsker-Varadhan LDP
\cite{donsker-varadhan:I-II,donsker-varadhan:III}
in several directions: First, our conditions are 
weaker. Second, when (DV3\upp) holds with an 
unbounded function $W$, 
the $\tau^{W_0}$-topology 
is finer and hence stronger than either the topology of
weak convergence, or the  $\tau$-topology, 
with respect to which the LDP for 
the empirical measures $\{L_n\}$ is usually 
established 
\cite{GOR:79,bolthausen:87,deuschel-stroock:book}. 
Third, apart from the LDP we also 
obtain precise large deviations expansions 
as in (\ref{eq:introBR}) for the partial sums 
with respect to (possibly unbounded) functionals 
$F\in \LWo .$

Following the Donsker-Varadhan papers,
a large amount of work has been 
done in establishing large deviations
properties of Markov chains under a
variety of different assumptions; see 
\cite{dembo-zeitouni:book,deuschel-stroock:book}
for detailed treatments. Under 
conditions similar to those in this paper, 
Ney and Nummelin have proved ``pinned'' 
large deviations principles in 
\cite{ney-nummelin:87a,ney-nummelin:87b}.
In a different vein, under much 
weaker assumptions (essentially under 
irreducibility alone) de Acosta 
\cite{deacosta:90} and Jain \cite{jain:90} 
have proved general large deviations 
lower bounds, but these are, in general, 
not tight.

 
One of the first 
places where the Feller continuity assumption 
of Donsker and Varadhan was relaxed is
Bolthausen's work \cite{bolthausen:87}.
There, a very  stringent condition on the chain is imposed, 
often referred to  in the literature as Stroock's uniform 
condition (U). 
In \Section{doeblin}  we   argue that (U) is much more restrictive
than the conditions we impose in this paper.
In particular, condition (U)   implies
Doeblin recurrence as well as the density 
assumption in (DV3+)~(ii).

More recently, Eichelsbacher and Schmock 
\cite{E-S:01} proved an LDP for the empirical 
measures of Markov chains, again under the uniform condition (U).
This LDP is proved in a strict subset of $\clM_1$, 
and with respect to a topology finer than the usual 
$\tau$-topology and similar in spirit to the 
$\tau^{W_0}$ topology introduced here.
In addition to (U),  the results of \cite{E-S:01} require strong 
integrability conditions that 
are {\em a priori} hard to verify:
In the above notation, in \cite{E-S:01}
it is assumed that for at least one 
unbounded function $W_0:\state\to\RL$,
we have
$\Expect_x[\exp\{a |W_0(\Phi(n))|\}]<\infty,$
uniformly over $n\geq 1$, for all real $a>0$.
This assumption is closely related to our
condition (DV3), and, as
we show in \Section{met}, (DV3) in particular 
provides a means for identifying a natural
class of functions $W_0$ satisfying this bound.

\section{Structural Assumptions}
\slabel{comparison}

There is a wide range of interrelated tools that have been 
used to establish large deviations 
properties for Markov processes and to develop parts
of the corresponding multiplicative ergodic theory.
Most of these tools rely on a functional-analytic
setting within which spectral properties of the process 
are examined. A brief survey of these approaches is
given in \cite{kontoyiannis-meyn:I}, where the main
results relied on the {\em geometric ergodicity} 
of the process. In this section we show how the 
assumptions used in prior work may be expressed 
in terms of the drift criteria introduced here
and describe the operator-theoretic setting upon
which all our subsequent results will be based.

\subsection{Drift Conditions}
\slabel{conditions}

Recall that the {\em (extended) generator} $\clA$ of $\bfPhi$ 
is defined as follows:  For a function $g:\state\to\Co$,
we write $\clA g=h$ if for each initial condition
$\Phi(0)=x\in\state$ the process
$\ell(t) \eqdef\sum_{s=0}^{t-1} h(\Phi(s)) -  g(\Phi(t)),$
$t\geq 1,$
is a local martingale with respect to the
natural filtration $\{\clF_t
=\sigma(\Phi(s),\, 0\le s\le t)\,:\,t\geq 1\}$.
In discrete time, the extended generator is simply
$\clA=P-I$, and its domain contains all measurable 
functions on $\state$.

The following drift conditions are considered in \cite{meyn-tweedie:book}
in discrete time,
\[
\begin{array}{rrcl}
\hbox{\textbf{(V2)}}   \qquad &  \clA V&\leq& -\delta + b\ind_C
\breakMed
\hbox{\textbf{(V3)}}   \qquad &  \clA V&\leq& -\delta W + b\ind_C
\breakMed
\hbox{\textbf{(V4)}}   \qquad &  \clA V&\leq& -\delta V + b\ind_C\,,
\end{array}
\]
where in each case $C$ is small, $V\colon\state\to(0,\infty]$ is
finite a.e. $[\psi]$, and $b<\infty$, $\delta>0$ are constants.
We further assume that $W$ is bounded below by unity in (V3), 
and that $V$ is bounded from below by unity in (V4).
It is easy to see that (V2)--(V4) are stated in order
of increasing strength: (V4) $\Rightarrow$ (V3) $\Rightarrow$ (V2).

Analogous multiplicative versions of these drift criteria 
are defined as follows,
\[
\begin{array}{rrcl}
\hbox{\textbf{(DV2)}}   \qquad &  \clH V&\leq& -\delta + b\ind_C
\breakMed
\hbox{\textbf{(DV3)}}   \qquad &  \clH V&\leq& -\delta W + b\ind_C
\breakMed
\hbox{\textbf{(DV4)}}   \qquad &  \clH V&\leq& -\delta V + b\ind_C\,,
\end{array}
\]
where $\clH$ is the nonlinear generator 
defined in \eq clH/.
The following implications follow easily 
from the definitions:

\begin{proposition}
\tlabel{ev}
For each $k=2,3,4$, the drift condition~{\em (DV$k$)}
implies~{\em (V$k$)}.
\end{proposition}

\proof
We provide a proof only for $k=3$ since all are similar.
Under (DV3),
$P e^V \le e^{V-W+b\ind_C}.$
Jensen's inequality gives $e^{PV} \le P e^V$, 
and taking logarithms gives (V3).
%
\qed

We find that \Proposition{ev} gives a poor bound in general.
\Theorem{v5} shows that (DV2) actually implies (V4).  Its proof
is given in the Appendix, after the proof of \Theorem{multReg}.

\begin{theorem}
\tlabel{v5}
{\em ((DV2) $\Rightarrow$ (V4)) }
Suppose $\bfPhi$ is $\psi$-irreducible and aperiodic.
If~{\em (DV2)} holds for some $V\colon\state\to(0,\infty]$,
then~{\em (V4)} holds for some $V_0$ which is {\em equivalent} to
$v_\eta\eqdef  e^{\eta V}$ for some $\eta>0$, in 
the sense that,
\[
V_0\in L_\infty^{v_\eta} \quad \hbox{and }\quad 
v_\eta\in L_\infty^{V_0}\,.
\]
\end{theorem}

\subsection{Spectral Theory Without Reversibility}
\slabel{STWR}

The spectral theory described in this paper and in 
\cite{kontoyiannis-meyn:I} is based on various 
operator semigroups
$\{\haP^n: n\in\nat\}$, where each $\haP^n$ is the $n$th
composition of a possibly non-positive kernel $\haP$.  
Examples are the transition kernel $P$; 
the multiplication kernel $I_G(x,dy)=G(x)\delta_x(dy)$.
for a given function $G$;
the {\em scaled kernel} defined by
\begin{equation} 
P_f(x,dy)  \eqdef  f(x)P(x,dy)\,, 
\elabel{If}
\end{equation}
for any function $F\colon\state\to \Co$ with
$f=e^{F}$;
and also the \textit{twisted kernel},
defined for a given function $h\colon\state\to (0,\infty)$ by
\begin{equation}
\cP_h(x, A) \eqdef [ I_{Ph}^{-1} P I_h]\, (x,A)
 =  \frac{\int_A P (x, dy) h (y) }
           {P h (x)} 
                     \qquad x\in\state,\ A\in\clB.
\elabel{cPh} 
\end{equation} 
This is a {\em probabilistic} kernel 
(i.e., a positive kernel with
$\cP_h(x,\state)=1$ for all $x$)
provided $Ph\, (x)<\infty$, $x\in\state$.
It is a generalization of the twisted kernel considered 
in \cite{kontoyiannis-meyn:I}, where
the function $h$ was taken as $h=\cf$ 
for a specially constructed $\cf$.
It may also be regarded as  a version of Doob's $h$-transform
\cite{pinsky:book}.

The most common approach to spectral decompositions 
for probabilistic semigroups $\{ P^n\}$
is to impose a reversibility condition 
\cite{fukstr86a,boveckgaykle00a,reytho00a}.
The motivation for this assumption
comes from the $L_2$ setting in which these
problems are typically posed, and the well-known fact that 
the semigroup $\{P^n\}$ is then self-adjoint.
We avoid a Hilbert space setting here and instead consider the 
weighted $L_\infty$ function spaces 
defined in (\ref{eq:weightedL}); cf.
\cite{kar85a,kar85b,horspi92a,meyn-tweedie:94b,kontoyiannis-meyn:I}.

The weighting function is determined by the particular 
drift condition satisfied by
the process.  In particular,
under (DV3) it follows from 
the convexity of $\clH$ (see 
\Proposition{Hconvex}) that for any  $0<\eta\le 1$ we have the bound,
\begin{equation}
 \clH(\eta V)\leq -\delta\eta W+b\eta\ind_C\,  ,\qquad \hbox{on $ S_V$\,,}
\elabel{JensenV3}
\end{equation}
which may be equivalently expressed 
as $P v_\eta\le e^{\eta[-\delta W + b\ind_C]} v_\eta$,
where $v_\eta \eqdef  e^{\eta V}$.  This bound implies that 
$P_f\colon L_\infty^{v_\eta}  \to  L_\infty^{v_\eta} $ is a bounded 
linear operator for any function $f$ satisfying 
$\|F^+\|_W\le \eta\delta$ (where $F^{+}\eqdef \max(F,0)$), and 
any $0\le \eta\le 1$.

Under any one of the above Lyapunov drift criteria, 
we will usually consider the function $v$ defined
in terms of the corresponding Lyapunov function $V$ on
$\state$ via $v=e^V$. For any such function 
$v\colon\state\to[1,\infty)$ and any
linear operator $\haP\colon L_\infty^{v}  \to  L_\infty^{v}$,
we denote the induced operator norm by,
\begin{equation}
\lll \haP \lll_{v}
\eqdef   \sup \Bigl\{ \frac{\|\haP h\|_{v}}{\|h\|_{v}} 
:  h\in L^{v}_\infty,\ \|h\|_{v}\neq 0\Bigr\}.
\elabel{Vnorm}
\end{equation}
The \textit{spectrum} $\clS(\haP)\subset \Co$ of
$\haP$ is the set of $z\in\Co$ 
such that the inverse $[Iz - \haP]^{-1}$ 
does not exist as a bounded linear operator on $ L_\infty^{v} $.
We let $\xi=\xi(\{\haP^n\})$ denote the {\em spectral radius}
of the semigroup $\{\haP^n\}$,
\begin{equation}
\xi(\{\haP^n\}) \eqdef  \lim_{n\to\infty} \lll \haP^n \lll_v^{1/n}.
\elabel{spectralRad}
\end{equation}
In general, the quantities $\lll\haP\lll_v$ and $\xi$
depend upon the particular weighting function $v$. 
If $\haP$ is a positive operator,
then $\xi$ is greater than or equal to the  \textit{generalized principal eigenvalue}, or \textit{g.p.e.} (see e.g.\ \cite{nummelin:book}),  
and they are actually equal 
under suitable regularity assumptions 
(see \cite{balaji-meyn,kontoyiannis-meyn:I}, and \Proposition{gpe} below).

As in \cite{kontoyiannis-meyn:I}, we say that $\haP$ 
admits a {\em spectral gap} if there exists $\epsilon_0>0$
such that the set
$\clS(\haP)\cap\{z:|z|\geq\xi -\epsilon_0\}$
is finite and contains only poles of finite multiplicity;
recall that 
$z_0\in \clS(\haP)$ is a
\textit{pole of (finite) multiplicity $n$} if:
\begin{MyEnumerate}  
\item
$z_0$ is isolated in $\clS(\haP)$, i.e.,
for some $\epsilon_{1}>0$ we
have $\{ z\in \clS(\haP) : |z-z_{0}| \le \epsilon_{1} \} =\{z_{0}\}$;
\item
The associated projection operator 
\begin{equation}\elabel{projection}
\haQ\eqdef \frac{1}{2\pi i}\int_{\partial\{
        z:|z-z_0|\leq\epsilon_1
        \}}[Iz-\haP]^{-1}dz\,,
\end{equation}
can be expressed as a finite linear combination
of some $\{s_i\}\subset L_\infty^v$,
$\{\nu_i\}\subset \clM_1^v$,
\ben
\haQ = \sum_{i,j=0}^{n-1} m_{i,j} [ s_i\otimes\nu_j]\,,
\een
where $[s\otimes\nu](x,dy) \eqdef s(x)\nu(dy)$.
\end{MyEnumerate}
See \cite[Sec.~4]{kontoyiannis-meyn:I} for more details.
Moreover, we say that $\haP$ is {\em $v$-uniform} 
if it admits a spectral gap and also there exists a unique 
pole $\lambda_\circ\in \clS(\haP)$
of multiplicity one,  satisfying $|\lambda_\circ| = \xi(\{\haP^t\})$.


Recall
that a Markov process $\bfPhi$ is called {\em geometrically ergodic}
\cite{kontoyiannis-meyn:I} 
or equivalently {\em $V$-uniformly ergodic} \cite{meyn-tweedie:book}
if it is positive recurrent,
and the semigroup converges in
the induced operator norm,
\[
\lll P^n - \One\otimes\pi \lll_V \to 0, \qquad n\to\infty\,,
\]
where $\One$ denotes the constant function $\One(x)\equiv 1.$
It is known that this is characterized by  condition~(V4).
Under this assumption, in \cite{kontoyiannis-meyn:I}
we proved that $\bfPhi$ satisfies a ``local'' large
deviations principle. In this paper
under the stronger condition~(DV3\upp) we show that
these local results can be extended to a full large deviations
principle.

The following result, taken from  
\cite[Proposition 4.6]{kontoyiannis-meyn:I},
says that geometric ergodicity 
is equivalent to the existence of a spectral gap:

\begin{theorem}
\tlabel{gap}
{\em (Spectral Gap \&\ (V4)) }
Let $\bfPhi$ by a $\psi$-irreducible and 
aperiodic Markov chain.
\begin{enumerate}
  \item[$(a)$]
        If $\bfPhi$ is geometrically ergodic
        with Lyapunov function 
        $V$,
        then its transition kernel $P$
        admits a spectral gap in $L_\infty^V$
        and it is $V$-uniform.
  \item[$(b)$]
        Conversely, if $P$ is $V_0$-uniform,
        then $\bfPhi$ is geometrically ergodic
        with respect to some Lyapunov function 
        $V\in L_\infty^{V_0}$. 
\end{enumerate}\end{theorem}

Next we want to investigate the corresponding 
relationship between
condition (DV3) and when the kernel $P$ has 
a discrete spectrum in $L_\infty^v$.
First we  establish an 
analogous `near equivalence' between assumption 
(DV3) and the notion of $v$-separability, and
in \Theorem{SepImpliesDiscrete} we show that 
$v$-separability implies the discrete spectrum
property.

For any 
$v\colon\state\to [1,\infty]$, finite a.e.\ $[\psi]$,
we say that the linear operator
$\haP\colon  L_\infty^{v}\to  L_\infty^{v}$ is
\textit{$v$-separable} if it can be approximated
uniformly by kernels with finite-rank. 
That is, for each $\epsilon > 0$, there exists 
a finite-rank operator $\haK_\epsilon$ such that
$\lll \haP - \haK_\epsilon \lll_v   \leq   \epsilon.$
Since the kernel $\haK_\epsilon$ has a finite-dimensional 
range space, we are assured of the existence of 
an integer $n\ge 1$,
functions $\{s_i : 1 \leq i \leq n \}\subset  L_\infty^v $, 
and probability measures
$\{\nu_i : 1 \leq i \leq n \}\subset  \clM_1^{v}$,   
such that $\haK_\epsilon$ may be expressed,
\begin{equation}
\haK_\epsilon (x, dy) = \sum_{i=1}^n s_i \otimes \nu_i\, .
\elabel{Ke}
\end{equation} 
Note that the eigenvalues of $\haK_\epsilon$ may be interpreted
as a \textit{pseudo-spectrum}; see \cite{davies:00}.

The following equivalence, established in the Appendix,
illustrates the intimate relationship between the 
essential ingredients of the Donsker-Varadhan conditions,
and the associated spectral theory as developed in this paper. 
Note that in \Theorem{gapB} the density assumption from part~(ii)
of~(DV3\upp) has been replaced by the more natural and weaker 
statement that $I_{C_W(r)} P^{T_0}$ is $v$-separable for all 
$r$.\footnote{The notation $I_A\haP$ for a set $A\in\clB$ and a kernel
$\haP$ is used to denote the kernel $\IND_A(x)\haP(x,dy)$.}
The fact that this is indeed weaker than the assumption
in (DV3)~(ii) follows from \Lemma{uniSep} in the Appendix.
Applications of \Theorem{gapB} to diffusions on $\Re^{n}$ and refinements in
this special case are developed in \cite{huimeysch04a}.

\begin{theorem}
\tlabel{gapB}
{\em ($v$-Separability \& (DV3)) }
Let $\bfPhi$ be a $\psi$-irreducible and aperiodic Markov chain
and let $T_0>0$ arbitrary.
The following are equivalent:
\begin{enumerate}
  \item[$(a)$] 
Condition {\rm (DV3)} holds with $V\colon \state\to[1,\infty)$; 
$W$ unbounded;
and $I_{C_W(r)} P^{T_0}$ is $v$-separable for  
all $r$, where $v=e^V$.

 \item[$(b)$]
The kernel $P^{T_0}$ is $v_0$-separable for some 
	unbounded function $v_0\colon\state\to [1,\infty)$.
\end{enumerate} 
\end{theorem}

We say that a linear operator $\haP:\Lv\to\Lv$
has a {\em discrete spectrum in $\Lv$} if its
spectrum $\clS$ has the property that
$\clS\cap K$ is finite, and contains
only poles of finite multiplicity,
for any compact set $K\subset\Co\setminus\{0\}$.
It is shown in  \Theorem{SepImpliesDiscrete}
that the spectrum of $P$ is discrete under the conditions of (b) above.

Taking a different operator-theoretic approach,
Deuschel and Stroock \cite{deuschel-stroock:book} 
prove large deviations results for the empirical
measures of stationary Markov chains under the
condition of hypercontractivity (or hypermixing). 
In particular, their conditions imply that for
some $T_0$, the kernel $P^{T_0}(x,dy)$ is a bounded
linear operator from $L_2(\pi)$ to $L_4(\pi)$, 
with norm equal to 1.

\subsection{Multiplicative Regularity}
\slabel{reg}

Recall the definition of the empirical
measures in~\eq OccEmp/, and the hitting times $\{\tau_A\}$ 
defined in \eq tauA/.
The next set of results characterize the drift 
criterion (DV3) in terms
of the following regularity assumptions:
\begin{quote}
\textbf{Regularity}
\nobreak
\begin{MyEnumerate}
\item 
A set $C\in\clB$ is called \textit{geometrically regular} 
if for any $A\in\clB^+$ there exists $\eta=\eta(A)>0$ such that
\[
\sup_{x\in C} \Expect_x [\exp(\eta \tau_A)] < \infty.
\]
The Markov process $\bfPhi$ is called \textit{geometrically regular} if there exists a geometrically regular set $C$, and $\eta>0$ such that
\[
\Expect_x [\exp(\eta\tau_C)] < \infty,\qquad x\in\state.
\]

\item 
A set $C\in\clB$ is called 
\textit{$H$-multiplicatively regular}  
($H$-m.-regular) if for any $A\in\clB^+$,
there exists  $\eta=\eta(A)>0$ satisfying,
\[
\sup_{x\in C} \Expect_x \Bigl[ 
\exp(\eta\tau_A\langle L_{\tau_A}, H\rangle)\Bigr] <\infty.
\]

The Markov process $\bfPhi$ is \textit{$H$-m.-regular} if there exists 
an $H$-m.-regular  set $C\in\clB$, and $\eta>0$ such that
\[ 
\Expect_x \Bigl[ \exp( \eta\tau_C\langle L_{\tau_C}, H\rangle)\Bigr] <\infty
                        \,,\qquad x\in\state\, .
\]
\end{MyEnumerate}
\end{quote}

In \cite[Theorem 15.0.1]{meyn-tweedie:book} a precise equivalence is
given between geometric regularity and the existence of a solution
to the drift inequality (V4).  The following analogous result 
shows that (DV3) characterizes multiplicative 
regularity.  A proof of \Theorem{multReg} 
is included in the Appendix.

\begin{theorem}
\tlabel {multReg}
{\em (Multiplicative Regularity $\Leftrightarrow$ (DV3)) }
For any $H\colon\state\to [1,\infty)$, the following
are equivalent:
\begin{MyEnumerate}
\item 
$\bfPhi$ is $H$-m.-regular;

\item 
The drift inequality~{\em (DV3)} holds 
for some
$V:\state\to(0,\infty)$ and with
$H\in \LW$.
\end{MyEnumerate} 
If either of these equivalent conditions hold, then for any $A\in\clB^+$,
there exists  $\epsilon>0$, $1\ge\eta>0$, and $B<\infty$ satisfying,
\[
	\Expect_x \Bigl[ 
	\exp\Big(\epsilon\tau_A\langle L_{\tau_A}, H\rangle 
	+ \eta V(\Phi(\tau_A))\Big)\Bigr] 
	\le \exp(\eta V(x) + B),\qquad x\in\state,
\] 
where $V$ is the  solution to {\em (DV3)} in {\em (ii)}.
\end{theorem}

In a similar vein,
in \cite{wu:01} the following condition is 
imposed for a diffusion on $\state = \Re^n$:
\begin{equation}
\parbox{.75\hsize}{\raggedright%
\textit{For any $n\ge 1$ there exists $K_n\subset\state$ compact, 
such that for any compact set $K\subset \state$,}
\[
\sup_{x\in K} \Expect_x[e^{n\tau_{K_n}}] < \infty.
\]
}
\elabel{Wu}
\end{equation}
In  \cite{wu:01,reytho01a}  it is shown that this condition is closely
related to the existence of a solution to (DV3), where the function
$W$ is further assumed to have compact sublevel sets. 
Under these assumptions, and under continuity assumptions
similar to those imposed in \cite{varadhan:book},  
it is possible to show that the operator $P^n$ 
is \textit{compact} for all $n>0$ 
\cite[Theorem 2.1]{reytho01a}, or \cite[Lemma 3.4]{deljun99}.
 
\archival{For a strong Feller chain with 
$\lll P\lll_{v}<\infty$ one can show that $I_{K}P^{2}$
is $v$-separable for any compact set 
$K\subset\hbox{support}\,\pi$.}

We show in \Proposition{wu}  that the bound assumed 
in \cite{wu:01} always holds under (DV3+).
We say that $G\colon\state\to\Re_{+}$ is \textit{coercive} if the sublevel
set $\{x : G(x)\le n\}$ is precompact for each $n\ge 1$.   Coercive functions
exist only when  $\state$ is $\sigma$-compact.
 
\begin{proposition}
\tlabel{wu} 
Let $\bfPhi$ be a $\psi$-irreducible and aperiodic Markov chain on $\state$.
Assume moreover that $\state=\Re^{n}$;  
that condition~{\rm (DV3\upp)} holds with $V\colon \state\to[1,\infty)$ 
continuous; $W$ unbounded;
and  the kernels  $\{ I_{C_W(r)} P^{T_0}: r\ge 1\}$ are 
$v$-separable for some $T_{0}\ge 1$.  Then,  there exists a sequence 
of compact sets $\{K_{n }: n\ge 1\}$   satisfying  \eq Wu/. 
 \end{proposition}

\proof  
\Lemma{SpectralRadius} combined with  \Proposition{coercive} 
implies that we may construct functions $(V_1, W_1)$
from $\state$ to $[1,\infty)$,   and a constant $b_1$ satisfying 
the following:  $\sup \{V(x) :  x\in C_{W_{1}} (r)\}  < \infty$
for each $r$; $W_1,V_1\in \LV$;  $W_1$  is coercive; and
$ \clH(V_1) \le V_1 - W_1 + b_1.$
\Lemma{VerySpecial} combined with continuity of $V$ 
then implies that  \eq Wu/ also holds, with
$K_r = \hbox{closure of}\, C_{W_1}(n_{r})$ 
for some sequence of positive integers 
$\{n_{r}\}$. 
\qed

%

 \Proposition{wu} has a partial converse:

\begin{proposition}
\tlabel{wub}
Suppose the chain $\bfPhi$ is $\psi$-irreducible and aperiodic.
Suppose moreover that $\state=\Re^{n}$;
that the support of $\psi$ has non-empty interior;
that $P$  has the   Feller property;
and that   there exists a sequence of compact sets $\{K_{n }: n\ge 1\}$
satisfying \eq Wu/.   Then Condition~{\rm (DV3)} holds with $V,W\colon \state\to[1,\infty)$  continuous  and   coercive.
\end{proposition}

\proof
\Proposition{fellerWu} asserts that there exists a  solution  
to the inequality
$\clH(V)\leq -\half W+b\ind_C$ with $(V,W)$  
continuous 
and coercive, $C$ compact, and $b<\infty$. 
Under the assumptions of the proposition, compact sets are small  
(combine  Proposition 6.2.8 with Theorem~5.5.7 of \cite{meyn-tweedie:book}).  
We may conclude that $C$ is small, and hence that (DV3) holds.
\qed

%
%

\subsection{Perron-Frobenius Theory}
\slabel{PF}

As in \cite{kontoyiannis-meyn:I} we find strong connections between the 
theory developed
in this  paper, and the Perron-Frobenius theory of positive semigroups, as 
developed in  \cite{nummelin:book}.

Suppose that $\{\haP^{n}: n\in\nat\}$ 
is a  semigroup of positive operators.
We assume that $\{\haP^{n}\}$ has 
finite spectral radius $\haxi$ in $\Lv$.
Then, the resolvent kernel defined by 
$\haR_{\lambda}\eqdef [I\lambda -\haP]^{-1}$
is a bounded linear operator on $\Lv$ for each $\lambda>\haxi$. 
We assume moreover that the semigroup is $\psi$-irreducible,
that is, whenever $A\in\clB$ satisfies $\psi(A)>0$, then 
$\sum_{k=0}^{\infty}\haP^{k}(x,A)>0$, for all $x\in\state$.
If $\bfPhi$ is a $\psi$-irreducible Markov chain, then for 
\textit{any} measurable function $F\colon\state\to\Re$, 
the kernel $\haP=P_{f}$
generates a $\psi$-irreducible semigroup.  
In general, under $\psi$-irreducibility of the semigroup,
one may find many solutions to the minorization condition,
\begin{equation}
\haR_\lambda (x, A) 
=
\sum_{k=0}^{\infty} \lambda^{-k-1} \haP^k 
\ge s(x) \nu (A),\qquad
x\in \state,\; A \in \clB,
\elabel{small}
\end{equation}
with $\lambda >0$, $s\in\clB^{+}$, and $\nu\in\clM^{+}$,
that is,
$s\colon\state\to\Re_{+}$ is measurable
with $\psi(s)>0$, and $\nu$ is a positive measure on 
$(\state,\clB)$ satisfying $\nu(\state)>0$.
The pair $(s,\nu)$ is then called small, 
just as in the probabilistic setting.  

Theorem~3.2 of \cite{nummelin:book}  states that there exists
a constant    $\halambda  \in (0,\infty]$, the \textit{generalized 
principal eigenvalue}, or {\em g.p.e.},
such that, for any  small function $s\in\clB^+$,
\begin{equation}
\sum_{k=0}^{\infty} \lambda^{-k-1} \haP^k s(x) 
\quad
 \left \{
\begin{array}{llr}
=\infty &   \hbox{for all $x\in\state $,}
    & \lambda < \halambda 
\\  \\
< \infty & \hbox{for a.e.\ $x\in\state$\ $[\psi]$,}
    & \lambda > \halambda  .
\end{array}
\right.
\elabel{gpeInt}
\end{equation}
The semigroup is said to be \textit{$\halambda$-transient} if
for one, and then all small pairs $(s,\nu)$, satisfying $s\in\clB^{+}$,
$\nu\in\clM^{+}$, we have
$\sum_{k=0}^{\infty} \halambda^{-k-1} \nu \haP^k s < \infty;$
otherwise it is called \textit{$\halambda  $-recurrent}.  
 
\Proposition{gpe} shows that the generalized 
principal eigenvalue coincides with the spectral
radius when considering positive semigroups that admit a spectral gap. 
Related results may be found in Theorem 4.4 and 
Proposition 4.5 of \cite{kontoyiannis-meyn:I}.

\begin{proposition}

\tlabel{gpe}
Suppose that $\{\haP^{n}: n\in\nat\}$ is a 
$\psi$-irreducible, positive semigroup.
Suppose moreover that the semigroup admits a spectral
  gap in $\Lv$, with finite spectral radius $\haxi$.  Then:
  \begin{MyEnumerate}
\item  $\haxi=\halambda$.
\item  The semigroup is $\halambda$-recurrent.
\item  $\haP$ is $v$-uniform.
\item
For any $\lambda >\haxi$, and any $(s,\nu)$ that solve \eq small/
with $s\in\clB^{+}$, $\nu\in\clM^{+}$,  the function
$h:= [I\hagamma-(\haR_{\lambda}- s\otimes\nu)]^{-1}s, 
\in\Lv$  is an eigenfunction.
  \end{MyEnumerate}
\end{proposition}

\proof
Suppose that either (i) or (ii) is false.  In either case, 
for all small pairs $(s,\nu)$,
\[
\lim_{\lambda\downarrow \haxi}  \nu \haR_{\lambda}s
=
\sum_{k=0}^{\infty} \haxi^{-k-1} \nu \haP^k s < \infty.
\]
It then follows that the projection operator $\haQ$ defined in \eq projection/
satisfies $\nu \haQ s=0$ for all small $s\in \Lv$, $\nu\in 
\clM_{1}^{v}$.  This is only
possible if $\haQ=0$, which is impossible under our assumption that the semigroup
admits a spectral gap.  

To complete the proof,  observe that the semigroup generated by the
kernel $\haR_{\lambda}$ also admits a spectral gap, with
spectral radius $\hagamma = (\lambda-\haxi)^{-1}$.   It follows that there is a closed ball
$D\subset\Co$ containing $\hagamma$ such that  the two
kernels below are bounded linear operators on $\Lv$ 
for each $\gamma\in D\setminus \{\hagamma\}$,
\[
X_{\gamma}=[I\gamma - \haR_{\lambda}]^{-1}\, ,
\qquad 
Y_{\gamma}=[I\gamma - (\haR_{\lambda}- s\otimes\nu)]^{-1}\, .
\]
From (i) and (ii) we know that $\haR_{\lambda}$ is $\hagamma$-recurrent, which implies that
$\nu Y_{\hagamma}s =1$, and that $\haP h = \haxi h$ (see   \cite[Theorem~5.1]{nummelin:book}).    Moreover, again from (i), (ii),  since  $\nu Y_{\hagamma}s <\infty$ 
it follows that the spectral radius of 
$(\haR_{\lambda}- s\otimes\nu)$ is strictly less
than $\hagamma$, which implies  (iii).   
Finally, since    $\lll Y_{\hagamma}\lll_{v}<\infty$ we may 
conclude that $h\in\Lv$, and this establishes (iv).
\qed 

On specializing to the kernels $\{P_{f}: F\in\LWo\}$ 
we obtain the following corollary.
Define for any measurable function $F\colon\state\to (-\infty,\infty]$:
\begin{equation} 
\begin{array}{rl}
 \hbox{(i)}  \   &   
 \Lambda(F) = \log(\lambda(F)) = 
 \hbox{ the logarithm of the g.p.e.
for $P_{f}$.}
\breakMed
 \hbox{(ii)}  \   & 
\Xi(F) = \log(\xi(F)) = \hbox{the logarithm of the spectral radius of $P_{f}$.}
\end{array}
\elabel{LambdaXiDefns}
\end{equation}

\begin{lemma}\tlabel{WinV}
Consider a $\psi$-irreducible Markov chain, and a measurable function
$G\colon\state\to \Re_{+}$.   If   $\Xi(G)<\infty$ then  $G\in\LV$.    
\end{lemma}

\proof
We have $\lll P_{g}^{n}\lll_{v}<\infty$ for some $n\ge 1$ when $\Xi(G)<\infty$.
Consequently, since $G$ and $V$ are assumed positive, we have
$g(x)\le P_{g}^{n}v\, (x)\le  \lll P_{g}^{n}\lll_{v} v(x),$
for all $x\in\state$.
\qed

\begin{proposition}
\tlabel{LambdaEqualsXi}
Under {\em (DV3+)} the functional $\Xi$ is finite-valued and 
convex on $\LWo$, and may be identified as the logarithm 
of the generalized principal eigenvalue:
\[
\Xi(F)=\Lambda(F),\qquad F\in\LWo.
\]
\end{proposition}

\proof
\Theorem{gapB} implies that $P_f$ is $v$-separable, and \Proposition{gpe} 
then gives the desired equivalence.  Convexity is established in 
\Lemma{SpectralRadiusConvex}.
\qed 

The spectral radius of the twisted kernel 
given in \eq cPh/ also has a simple representation, 
when the function $h$ is chosen as a solution 
to the multiplicative Poisson equation: 

\begin{proposition}
\tlabel{twisted}
Assume that the Markov chain $\bfPhi$
satisfies condition~{\em (DV3\upp)} with $W$ unbounded.   
For real-valued $F\in\LWo$,  the twisted kernel 
$\cP_{\cf}$ satisfies {\em (DV3\upp)} with 
Lyapunov function $\cV\eqdef V- \cF +c$ for $c\ge 0$ sufficiently large.   
Consequently, the semigroup generated by the
twisted kernel has a discrete spectrum in $L_{\infty}^{\cv}$, and
its log-spectral radius has the representation,
\[
\cXi(G) = \Xi(F+G) \,,\qquad G\in\LWo.
\]
\end{proposition}

\proof
The kernels $P_{f}$ and $\cP_{\cf}$ are related by a scaling and a 
similarity transformation,
\[
\cP_{\cf}= \lambda(f)^{-1} I_{\cf}^{-1}P_{f} I_{\cf} .
\]
It follows that (DV3\upp)~(i) is satisfied with the Lyapunov function $\cV$, and we have
$\cV\ge 1$ for sufficiently large $c$ since $\cf\in\Lv$.  
The representation of $\cXi$ also follows from the above relationship
between $\cP_{\cf}$ and $P_{f}$.

The density condition  (DV3\upp)~(ii) follows similarly.  Letting $b_r=\|  \lambda(f)^{-1}f\ind_{C_{W}(r)}\|_{\infty}$,
we have, under the transition law $\cP_{\cf}$,
\[
   \cProb_{x}\{ \Phi(T_{0})\in A ,\ \tau_{C_{W}(r)} >T_{0}\} \le 
   \cf^{-1} (x)b_r^{T_{0}}\check{\beta}_r(A),\quad A\in\clB,\ x\in C_{W}(r),
\]
where $\check{\beta}_{r}(dx) = \beta_{r}(dx)\cf(x)$.  To establish  (DV3\upp)~(ii)
it remains to show that
$\cf^{-1}$ is bounded on $C_{W}(r)$.

  Since the set $C_{W}(r)$
is small for the semigroup $\{ P_{f}^t : t\ge 0\}$, there exists $\epsilon>0$, $T_{1}<\infty$, and
a probability distribution $\nu$ such that
\[
P_{f}^{T_1}(x,dy)\ge \epsilon\nu(dy),\qquad x\in C_{W}(r),\  y\in\state.
\]
It follows that
\[
 \lambda(f)^{-T_1}\cf
 =
P_{f}^{T_1}\cf\ge \epsilon\nu(\cf),\qquad x\in C_{W}(r).
\]
Consequently, $\cf^{-1}$ is bounded on $C_{W}(r)$. 
\qed

\subsection{Doeblin and Uniform Conditions}
\slabel{doeblin}

The uniform upper bound in condition~(DV3\upp)~(ii) 
is easily verified in many models.  
Consider first the special case of a discrete time chain
$\bfPhi$ with a countable state space $\state$, and
with $W$ such that $C_W(r)$ is finite for all $r<\|W\|_\infty$.  
In this case we may
take $T_0=1$ in~(DV3+)~(ii), and set
\[
\beta_r(A) = \sum_{x\in C_W(r)} P(x,A),
                        \qquad A\in\clB\, .
\]
This is the starting point for the bounds obtained in \cite{balaji-meyn}.


A common assumption for general state space models is the
following:


\begin{quote}
\textbf{Condition~(U)}  \quad
There exist $1\leq T_1\leq T_2$
and a constant $b_0\geq 1$, such that
\be
P^{T_1} (x,A)\leq b_0 \frac{1}{T_2} \sum_{t=1}^{T_2}  P^t(y,A)\,,
\qquad x,y\in\state,\;A\in\clB.
\elabel{Uhat}
\ee
\end{quote}


See  \cite{deuschel-stroock:book,dembo-zeitouni:book},
as well as \cite{varadhan:book,iscneynum85a,jensen:91}.
It is obvious that \eq Uhat/ implies
the validity of the upper bound in our 
assumption~(DV3\upp)~(ii).
Somewhat surprisingly, Condition~(U) also
implies a corresponding lower bound, and 
moreover we may take the bounding
measure equal to the invariant measure $\pi$:
\begin{proposition}
\tlabel{uhat}
Suppose that $\bfPhi$ is an aperiodic, $\psi$-irreducible chain.  
Then, condition~{\rm (U)}
holds if and only if there is a probability
measure $\pi$ on $(\state,\clB)$, a constant $N_0\ge 1$, 
and a sequence of non-negative numbers $\{\delta_n : n\ge N_0\}$,
satisfying,
\be
\begin{array}{rcl}
| P^n(x,A) -\pi(A) | & \leq &\delta_n \pi(A)\,,
\qquad\;A\in\clB,\;x\in\state,\ n\ge N_0;
\breakMed
\lim_{n\to\infty} \delta_n &=& 0.
\end{array}
\label{e:altU}
\ee
\end{proposition}


\proof
It is enough to show that condition~(U) 
implies the sequence of bounds given in \eq altU/.
 
Condition~(U) implies the following minorization,
\[
 \sum_{t=1}^{T_2} P^t(y,A)
\ge \epsilon \nu(A),\qquad A\in\clB, y\in\state,
\]
where $\epsilon = T_2 b_0^{-1}$, and $\nu(A) = 
P^{T_1} (x_0,A)$, $A\in\clB$, with $x_0\in\state$ arbitrary.
Since the chain is assumed aperiodic and $\psi$-irreducible, 
it follows that the chain is \textit{uniformly ergodic}, 
a property somewhat stronger than Doeblin's condition
\cite[Theorem~16.2.2]{meyn-tweedie:book}.  
Consequently, there exists an invariant probability measure
$\pi$, and constants  $B_0<\infty, b_0>0$ such that,
\begin{equation}
\lll P^n - \one\otimes\pi \lll_1 \le  e^{-b_0n + B_0},\qquad n\in\TT.
\elabel{uniErgo}
\end{equation}

Condition~(U) then gives the following upper bound:  On
multiplying \eq Uhat/ by $\pi(dy)$,
and integrating over $y\in\state$, we obtain,
\[
P^{T_1}(x,A)
\leq b_0   \pi(A),\qquad x\in\state,\;A\in\clB.
\]
Let $\Gamma$ denote the bivariate measure given by,
$\Gamma(dx,dy) = \pi(dx) P^{T_1}(x,dy)$, for $x,y\in\state$.
The previous bound implies that $\Gamma$ has a density 
$p(x,y;T_1)$ with respect to $\pi\times\pi$, 
where $p(\varble,\varble;T_1)$ is jointly measurable, and may be chosen so that it satisfies the strict upper bound,
$p(x,y;T_1) \le  b_0,$ for $x,y\in\state.$
The probability measure $\Gamma$ has common 
one-dimensional marginals
(equal to $\pi$).  Consequently,  we must
have $\int p(x,y;T_1)\pi(dx) = 1$ a.e.\ $y\in\state$ $[\pi]$.  

For $n\ge 2T_1$ we \textit{define} the density $p(x,y;n)$ via,
\[
p(x,y;n) \eqdef  \int P^{n- T_1}(x,dz)p(z,y;T_1),\qquad x,y\in\state.
\]
We have the upper bound
$\sup_{x,y}p(x,y;n)\le b_0$ for all $n\ge T_1$ since $P^k$ is an 
$L_\infty$-contraction for any $k\ge 0$.  Combining this bound with \eq uniErgo/
gives the strict bound,
\[
\begin{array}{rcl}
|p(x,y;n)-1 |  &=& \Bigl| \int P^{n-T_1} (x,dz) (p(z,y;T_1) - 1)  \, \Bigr| 
\breakMed
& = & \Bigl| \int P^{n-T_1} (x,dz) p(z,y;T_1) - \int \pi(dz) p(z,y;T_1)  \, \Bigr| 
\breakMed
&\le & b_0 \lll P^{n-T_1} - \pi \lll_1 \le b_0  e^{B_0 - b_0(n-T_1)},\qquad n\ge T_1, x,y\in\state.
\end{array}
\]
This easily implies the result.
\qed


Note that, for the special case of reflected Brownian
motion on a compact domain, a similar result is established
in \cite{bass-hsu:91}.

We have already noted in the above proof
that the lower bound in (\ref{e:altU}) 
implies the Doeblin condition, which
is known to be equivalent to (V4) with $V$ bounded
for a $\psi$-irreducible chain
\cite[Theorem~16.2.2]{meyn-tweedie:book}.    
Consequently, condition~(U) frequently 
holds for models on compact state spaces
but it rarely holds for models on $\Re^n$.
We summarize this and related correspondences with drift
criteria here. 

\newpage

\begin{proposition} 
\tlabel{uniAssumptions}
Suppose that $\bfPhi$ is an aperiodic, $\psi$-irreducible chain.  
\begin{MyEnumerate}
\item 
If $\bfPhi$ satisfies Doeblin's condition, 
then {\rm (DV4)} holds with respect to the
Lyapunov function $V\equiv 1$.

\item 
If $\bfPhi$ satisfies condition~{\rm (U)}
and $V_0\colon\state\to [1,\infty)$ is given 
with  $\lll P\lll_{v_0}<\infty$,
then  {\rm (DV4)} 
holds for a function $V\colon\state\to[1,\infty)$ that is equivalent to 
$V_0$. And, trivially, 
part~{\em (ii)} of condition {\em (DV3+)} also holds.
\end{MyEnumerate}
\end{proposition}
 
\proof
Result (i) is a consequence of  
\cite[Theorems~16.2.3 and 16.2.3]{meyn-tweedie:book} which state 
that the state space $\state$ is small under these assumptions, 
and hence (DV4) holds with $V\equiv 1$.

To prove (ii) we define,
\[
V(x) \eqdef 1+\log\Bigl(\Expect_x\Bigl[ \exp \Bigl(\epsilon \sum_{i=0}^{{T_1}-1} r^i V_0(\Phi(i)) \Bigr)\Bigr]\Bigr),
\qquad x\in\state,
\]
where $r>1$ is arbitrary, and $\epsilon>0$ is to be determined. The functions $V$ and $V_0$  are equivalent when $\epsilon \le {T_1}^{-1}r^{-{T_1}+1}$ since then by H\"older's inequality,
\[
V(x) \le 1+\frac{1}{{T_1}} \sum_{i=0}^{{T_1}-1}  \log\bigl(\Expect_x[ \exp ({T_1} \epsilon  r^i V_0(\Phi(i)) )]\bigr),
\qquad x\in\state,
\]
and the right hand side is in $L_\infty^{V_0}$  since    $\lll P^i\lll_{v_0}<\infty$
for $i\ge 0$ under the assumptions of (ii).  Moreover, we have $V\ge \epsilon V_0$ by considering only
the first term in the definition of $V$.  Hence $V\in L_\infty^{V_0}$ and 
$V_0\in L_\infty^V$, which shows that $V$ and $V_0$ are equivalent.
We assume henceforth that this bound holds on $\epsilon$.

H\"older's inequality also gives the bound,
\[
\begin{array}{rcl}
P e^V &=& \Expect_x\Bigl[ \exp \Bigl(\epsilon \sum_{i=0}^{n-1} r^i V_0(\Phi(i+1)) \Bigr)\Bigr]
\breakMed
&\le &
 \Expect_x\Bigl[ \exp \Bigl(pr^{-1}\epsilon \sum_{i=1}^{n-1} r^i V_0(\Phi(i)) \Bigr)\Bigr]^{1/p}
 \Expect_x\Bigl[ \exp \Bigl(q r^{T_1-1} \epsilon  V_0(\Phi(T_1)) \Bigr)\Bigr]^{1/q},
\end{array}
\]
where we set $p=r>1$ and $q=r(r-1)^{-1}>1$. 
Under Condition~(U) we   have $\|P^{T_1} e^{V_0}\|_\infty<\infty$.
Consequently,  
provided $\epsilon>0$ is chosen so
that 
$q r^{T_1-1} \epsilon <1$ we then have, for some constant $b_1$,
\[
\clH(V) \eqdef \log(Pe^V) - V
\le  -(1-r^{-1}) V + b_1.
\]
This implies the result since the state space is small.
\qed

\subsection{Donsker-Varadhan Theory}
\slabel{DV}

In Donsker and Varadhan's classic papers 
\cite{donsker-varadhan:I-II,donsker-varadhan:III,donsker-varadhan:IV}
there are two distinct sets of assumptions that are imposed 
for ensuring the existence of a large deviations principle,
roughly corresponding to parts~(i) and~(ii) of our
condition~(DV3\upp).


\newpage

\noindent
{\bf Lyapunov criteria. }
The Lyapunov function criterion of
\cite{donsker-varadhan:IV,varadhan:book} is essentially equivalent 
to (DV3), with the additional constraint that the function 
$W$ has compact sublevel sets; see conditions 
(1)--(5) on \cite[p.~34]{varadhan:book}. 
In the general case (when $\state$ is not
compact) this implies that (DV3) holds with
an unbounded $W$.

It is worth noting that the nonlinear generator 
is implicitly already present in the Donsker-Varadhan work, 
visible both in the form of the rate function, 
and in the assumptions imposed in 
\cite{donsker-varadhan:III,donsker-varadhan:IV,varadhan:book}.

\medskip

\noindent
{\bf Continuity and density assumptions. }
In \cite{varadhan:book} two additional conditions
are imposed on $\bfPhi$. It is assumed that the
chain satisfies a strong version of the Feller
property, and that for each $x$, $P(x,dy)$ has a
continuous density $p_x(y)$ with respect 
to some reference measure $\alpha(dy)$ 
which is independent of $x$. 

These rather
strong assumptions are easily seen to imply 
condition~(DV3\upp)~(ii) when $W$ is 
coercive, so that the sets $C_W(r)$ are 
pre-compact.

\section{Multiplicative Ergodic Theory}
\slabel{met}
  
\subsection{Multiplicative Mean Ergodic Theorems}
\slabel{met-summary}

The main results of this section are
summarized in the following two theorems.
In particular, the multiplicative mean
ergodic theorem given in \eq mainMMET/
will play a central role in the proofs
of the large deviations limit theorems
in~\Section{ldp}.
For all these results we will assume that 
$\bfPhi$ satisfies (DV3) with an unbounded
function $W$.  As above, we let $\clB^+$ denote the set 
of functions $h\colon \state \to [0,\infty]$ 
with $\psi(h)>0$; for $A\in\clB$ we write $A\in\clB^+$ 
if $\psi(A)>0$; and let 
$\clM^{+}$ denote the set of positive measures
on $\clB$ satisfying $\mu(\state)>0$.

As in \eq Wo/ in the Introduction, we choose an arbitrary  
measurable function $W_0:\state\to[1,\infty)$ 
in $\LW$,
whose growth at infinity is strictly slower than $W$.
This may be expressed in terms of the weighted $L_\infty$ norm via,
\begin{equation}
\lim_{r\to\infty} \|W_0\ind_{C_W(r)^c} \|_W = 0\, ,
\elabel{littleOh}
\end{equation}
where $\{ C_W(r)\}$ are the sublevel sets of $W$ defined in \eq Cn/.
The function $W_0$ is fixed throughout this section. 

Given $F\in \LWo $ and an arbitrary $\alpha\in\Co$, we recall from
\cite{kontoyiannis-meyn:I} the notation
$\haP_\alpha\eqdef e^{\alpha F}P$,
 and
\[
\clS_\alpha\eqdef \clS(\haP_\alpha) \eqdef  
\mbox{spectrum of } \haP_\alpha\;\mbox{in } 
L_\infty^v\,,
\]
where $v\eqdef e^V$ and $V$ is the Lyapunov
function in (DV3+).  

Next, we collect the main results of this section
in the following theorem. Recall the definition 
of the empirical measures $\{L_n\}$ from \eq OccEmp/.


\begin{theorem}
\tlabel{mainSpectral}
{\em (Multiplicative Mean Ergodic Theorem)}
Assume that the Markov chain $\bfPhi$
satisfies condition~{\em (DV3\upp)} with an unbounded 
$W$. For any $m>0,M>0$
there exist $\bara>m,\baromega>0$ such that 
for any real-valued $F\in \LWo $ with $\|F\|_{W_0}\leq M$,
and any $\alpha$ in the compact set
\[
\Omega=\Omega(\bara,\baromega):=\{\alpha=a+i\omega\in\Co\,:\,
|a|\leq \bara,\;\mbox{and}\;|\omega|\leq\baromega\}\,,
\]
we have:
\begin{MyEnumerate}
        \item
        There is a maximal, isolated eigenvalue $\lambda(\alpha F)
\in \clS_\alpha$ satisfying $| \lambda(\alpha F)|=\xi(\alpha F)$.  
Furthermore, 
$\LA(\alpha F)\eqdef \log(\lambda(\alpha F))$
is analytic as a function of $\alpha\in\Omega$,
and for real $\alpha$ it coincides with the 
log-generalized principal eigenvalue of \Section{PF}.
        \item 
Corresponding to each eigenvalue $\la(\alpha F)$,
there is an eigenfunction $\cf_\alpha \in L_\infty^v$
and an eigenmeasure $\cmu_\alpha \in \clM_1^v$, 
where $v\eqdef e^V$, normalized so that 
$\cmu_\alpha(\cf_\alpha)=\cmu_\alpha(\state)=1$.
The function $\cf_\alpha$ solves the 
{\em multiplicative Poisson equation},
$$\haP_\alpha \cf_\alpha = \lambda(\alpha F) \cf_\alpha\,,$$
and the measure $\cmu_\alpha$ is a corresponding 
eigenmeasure:
$\cmu_\alpha \haP_\alpha = \lambda(\alpha F)\cmu_\alpha.$
        \item 
There exist constants
$b_0>0$, $ B_0 < \infty$, independent of $\alpha$,
such that for all $x\in\state$, $\alpha\in\Omega$, $n\ge 1$,
\be
\Bigl| \Expect_x 
    \Bigl[ \exp 
        \Bigl ( 
           n[
                \alpha \langle L_n, F \rangle
                - \LA(\alpha F)
            ]
        \Bigr)
    \Bigr]
                - \cf_\alpha(x)
\Bigr|
\leq
                        |\alpha| v(x) e^{B_0 - b_0 n} \, .
                \elabel{mainMMET}
\ee
\end{MyEnumerate}
\end{theorem}
   
\proof
\Lemma{uniSep} in the Appendix shows that 
$(P_{f_0})^{2T_0+2}$ is $v_\eta$-separable
for any $F_0\in\LWo$,
and \Theorem{SepImpliesDiscrete} then implies that 
the spectrum of $P_{f_0}$ is discrete.
It follows that solutions to the  eigenvalue problem 
for $P_{f_0}$ exist with $\cf_0\in L_\infty^{v_\eta}$, 
$\cmu_0\in\clM_1^{v_\eta}$.  The eigenvalue satisfies   
$|\lambda(F_{0}) |=\xi(F_{0}) <\infty$.
Smoothness of $\Lambda$ is established in
\Proposition{Lcts}.

\Theorem{localMMET} establishes the limit (iii) 
for $\alpha\in\Co$ in a neighborhood of the origin.  

Consider then the twisted kernel $\cP=\cP_{\cf_a}$,
where $a$ is real.  \Proposition{twisted} states that
this satisfies (DV3+) with 
Lyapunov function $\cV\eqdef V/\cf_a$. 
 An application of \Theorem{localMMET}
to this kernel then implies   a uniform bound 
of the form (iii) for $\alpha$ in a 
neighborhood of $a$. For any given $\bara>0$ we may appeal to compactness
of the line-segment $\{ a\in\Re : |a|\le \bara\}$ to construct $\baromega>0$
such that \eq mainMMET/  holds for $\alpha\in\Omega$.   
\qed

%

\done{omitted part of global MMET: 
If $F$ is real-valued, then the twisted kernel 
$\cP_{\cf}$ satisfies (DV3), and 
$\ind_{C_W(r)} \cP_{\cf}^{T_0}$ is 
$v_{\eta_0}$-separable each $0\le\eta_0\le1$, 
$r\ge 1$.  Should this go somewhere else?
Write a proposition that also mentions the fact that
$\Lambda(F+G)=\cLambda(F)$?}

We note that this result has many immediate extensions.  
In particular,
if condition (DV3\upp) is satisfied, then this condition also holds 
with $(V,W)$ replaced by $ (1-\eta+\eta V,  W)$ for any $0<\eta<1$.  
Consequently, $\cf \in L_{\infty}^{v_{\eta}}$ for \textit{any} 
$0<\eta\le 1$ when $F\in\LWo$.

\medskip

Part~(iii) of the theorem is at the heart 
of the proof of all the large deviations properties
we establish in \Section{ldp}. For example,
from \eq  mainMMET/ we easily obtain that,
for any $F\in\LWo$, the log-moment
generating functions of the partial sums
$$S_n=\sum_{i=0}^{n-1}F(\Phi(i))=n\langle L_n,F\rangle$$
converge uniformly and exponentially fast:
\be
\frac{1}{n}
\log \Expect_x
  \bigl [ \exp( \alpha n\langle L_n, F \rangle  )
  \bigr ]
	\to
	\LA(\alpha F),\;\;\;\;n\to\infty.
\label{eq:mmetL}
\ee
We therefore think of $\LA(\alpha F)$ as the limiting
log-moment generating function of the partial sums
$\{S_n\}$ corresponding to the function $F$, 
and much of our effort in the following
two section will be devoted to examining the regularity 
properties of $\LA$ and its convex dual $\LA^*$.
 
Following \cite{kontoyiannis-meyn:I},
next we give a weaker multiplicative mean
ergodic theorem for $\alpha$ in 
a   neighborhood of the imaginary axis.
Recall the following terminology:
The {\em asymptotic variance} $\sigma^2(F)$
of a function $F:\state\to\Re$ is defined to be 
variance obtained in the corresponding 
Central Limit Theorem for the partial
sums of $F(\Phi(n))$, assuming it exists.
For a $V$-uniformly ergodic
(or, equivalently, a geometrically ergodic)
chain, the asymptotic variance is finite for 
any function $F $ satisfying $F^2\in L_\infty^V$,   and
\cite[Theorem17.0.1]{meyn-tweedie:book} gives the representation,
\begin{equation}
\sigma^2(F) = \lim_{n\to\infty}
n \Expect_\pi[ (\langle L_n, F\rangle - \pi(F) )^2]\, .
\elabel{avariance}
\end{equation} 

A function $F:\state\to\Re$ is called
{\em lattice} if there are
$h>0$ and $0\le d < h$, such that
$[F(x)-d]/h$
is an integer for all
$x\in\state$.
The minimal $h$ for which this
holds is called the {\em span} of $F$.
If the function $F$ can be written as a sum,
$F=F_0 + F_\ell,$
where $F_\ell$ is lattice with span $h$
and $F_0$ has zero asymptotic variance 
then $F$ is called {\em almost-lattice} (and $h$
is its span).
Otherwise, $F$ is called {\em strongly non-lattice.}
The lattice condition is discussed in more detail
in \cite{kontoyiannis-meyn:I}.  
The proof of the following result follows from \Theorem{mainSpectral}
and the arguments used in the proof of \cite[Theorem~4.2]{kontoyiannis-meyn:I}.   
\begin{theorem}
{\em (Bounds Around the $i\omega$-Axis)}
\tlabel{mainSpectral2}
Assume that the Markov chain
$\bfPhi$ satisfies condition {\em (DV3\upp)}
with an unbounded $W$, 
and that $F\in \LWo $
is real-valued.
\begin{itemize}
\item[{\em (NL)}]
If $F$ is strongly non-lattice,
then for any $m>0$ and
$0<\omega_0 <\omega_1< \infty$,
there exist $\bara>m$, $b_0>0$, $B_0 < \infty$
(possibly different than in \Theorem{mainSpectral}),
such that
\be
\Bigl | \Expect_x 
  \Bigl [ \exp 
    \Bigl (
        n[\alpha \langle L_n, F \rangle - \LA(aF)]
    \Bigr )
  \Bigr ] 
\Bigr|
\leq
  v(x)
  e^{B_0-b_0 n}\, ,  \qquad  x \in \state,\; n \geq 1,
\label{eq:ImMMET}
\ee
for all $\alpha=a+i\omega$ with $|a| \le \bara$
and $\omega_0\leq |\omega|\leq\omega_1$,
where $v\eqdef e^V$.
 
\item[{\em (L)}]
If $F$ is almost-lattice with span $h>0$,
then 
for any $m>0$ and $\epsilon>0$, there exist
$\bara>m$, $b_0>0$, and $B_0 < \infty$
(possibly different than above and
in \Theorem{mainSpectral}),
such that (\ref{eq:ImMMET}) holds
for all $\alpha=a+i\omega$ with $|a| \le \bara$
and $\epsilon\leq|\omega|\leq 2\pi/h - \epsilon.$
\end{itemize}
\end{theorem}

\subsection{Spectral Theory of {\Large $v$}-Separable Operators}
\slabel{separable}

The following continuity result allows perturbation analysis 
to establish a spectral gap under (DV3).  
Recall that we set  $v_\eta:=e^{\eta V}$;  for any real-valued 
$F\in \LW$ we define  $f\eqdef e^F$; and we let  $P_f$ denote  
the kernel $P_f(x,dy)\eqdef f(x)P(x,dy)$. 

\begin{lemma}
\tlabel{FrechetLemma}
Suppose that $\bfPhi$ is $\psi$-irreducible and aperiodic,
and that condition~{\em (DV3)} is satisfied.
Then, for  $0<\eta\le1$, $n\ge 1$,
there exists   $b_{\eta,n}<\infty$,
such that for any   $F,G\in\LWo$,  
\[
\lll P_f  - P_g\lll_{v_\eta} \le b_{\eta,n} \|F-G\|_{W_0},
\]
whenever $\|F\|_{W_{0}} \le n$,
and $\|G\|_{W_{0}} \le n$. Moreover, for any 
$h\in L_\infty^{v_\eta}$ the map $F\mapsto P_f h$
is Frechet differentiable as a function 
from $\LWo$ to $L_\infty^{v_\eta}$.
\end{lemma}

\proof 
We have from the definition of the induced operator norm,
\[
\begin{array}{rcl}
\lll P_f  - P_g \lll_{v_\eta}
&=& \sup_{x\in\state} 
\Bigl(
 |f(x)-g(x)| \frac{P v_\eta \, (x)}{v_\eta(x)}
\Bigr) 
\breakBig
&\le & 
\sup_{x\in\state} |f(x)-g(x)| \exp\bigl(-\eta\delta W(x) + \eta b\bigr) \, .
\end{array}
\]
Also, we have the elementary bounds, for all $x\in\state$,
\[
\begin{array}{rcl}
 |f(x)-g(x)| = |e^{F(x)}-e^{G(x)} |
&\le & 
|F(x)-G(x)| e^{|F(x)|+|G(x)|}
\breakMed
&\le &
 \|F-G\|_{W_0} W_0(x) \exp\bigl((\|F\|_{W_0} +\|G\|_{W_0} )W_0(x) \bigr) 
\breakMed
&\le &
 \|F-G\|_{W_0} \exp\bigl((1+\|F\|_{W_0} +\|G\|_{W_0} )W_0(x)\bigr)  \, .
\end{array}
\]
Combining these bounds gives,
\begin{equation}
\lll P_f  - P_g\lll_{v_\eta}
\le
\|F-G\|_{W_0} \sup_{x\in\state}\Bigl(
\exp\bigl( (1+\|F\|_{W_0} +\|G\|_{W_0} )W_0(x) 
-\eta\delta W(x) + \eta b \bigr)\Bigr).
\elabel{PgBdd}
\end{equation}
The supremum is bounded under the assumptions of the proposition, 
which establishes the desired bound.

We now show that, for any  
given $h\in L_\infty^{v_\eta}$, $F\in\LWo$, 
the map $G\mapsto  I_{G-F}P_{f} h$
represents the Frechet derivative of $P_{f}h$.  We begin with   
the mean value theorem,
\[
 P_{f} h- P_{g} h-   I_{G-F}P_{f} h
 =
 (G-F)[ P_{f_{\theta}}h -   P_{f} h]
\]
where $F_{\theta} = \theta F + (1-\theta)G$ for some   $\theta\colon\state\to (0,1)$.
The  bounds leading up to \eq PgBdd/ 
then lead to the   following bound, for all $x\in\state$,
\[ \begin{array}{rcl}
 \lefteqn{\bigl|[P_{f} h- P_{g} h-   I_{G-F}P_{f} h]\, (x) \bigr| }   &   &   
 \breakMed  &  \le &   
\Bigl( \|G-F\|_{W_{0}}  W_{0} (x)\Bigr)
\Bigl( \|F-G\|_{W_0}  
\exp\bigl( (1+\|F\|_{W_0} +\|F_{\theta}\|_{W_0} )W_0(x) 
-\eta\delta W(x) + \eta b \bigr)\Bigr).  \end{array}\]
It follows that  there exists $ b_1<\infty$ such that 
\[
\lll [P_{f} h- P_{g} h-   I_{G-F}P_{f} h]\lll_{v_{\eta}}
\le
 b_{1} \|F-G\|_{W_0}^{2} \,  \qquad G\in\LWo,\ \| F-G\|_{W_{0}}\le 1\, ,
\]
which establishes Frechet differentiability.
\qed


Next we present a  \textit{local} result, 
in the sense that it holds for all 
$F$ with sufficiently small  $\LW$-norm, 
where the precise bound on $\|F\|_W$ is not explicit.  
Although a value can be computed as in
\cite{kontoyiannis-meyn:I}, it is not of  
a very attractive form.  Note that \Theorem{localMMET} does
not require the density condition used in (DV3+).

The definition of the empirical measures $\{L_n\}$ 
is given in \eq OccEmp/.
\begin{theorem}
\tlabel{localMMET}
{\em (Local Multiplicative Mean Ergodic Theorem) }
Suppose that $\bfPhi$ is $\psi$-irreducible and aperiodic,
and that condition~{\em (DV3)} is satisfied.
Then there exists $\epsilon_0>0$,  
$0<\eta_0\leq 1$, such that for any complex-valued 
$F\in \LW$ satisfying $\|F\|_{W} \le \epsilon_0$, 
and any $0<\eta\le \eta_0$:
\begin{MyEnumerate}
\item 
There exist solutions $\lambda$, $\cf$ and $\cmu$
to the eigenvalue problems 
\begin{equation}
P_f \cf = \lambda \cf,\quad 
\cmu P_f   = \lambda \cmu\,.
\elabel{eigen}
\end{equation}
These solutions satisfy
$\cf\in  L_\infty^{v_\eta}$,
$\cmu\in \clM_1^{v_\eta}$,
$\cmu(\state)=\cmu(\cf)=1$,
and the eigenvalue $\la=\la(F)\in\Co$
satisfies $|\la| = \xi(\{P_f^t\})$. Moreover, the solutions are 
uniformly continuous on this domain: 
For some $b_\eta<\infty$,  
\[
|\Lambda(F)-\Lambda(G)|\le b_\eta \|F-G\|_{W},\qquad
|\cf-\cg|_{v_\eta} \le b_\eta \|F-G\|_{W}\, ,
\] 
whenever $F,G\in\LW$ satisfy $ \|F\|_W\le \epsilon_0$, 
$\|G\|_W\le \epsilon_0$.

\item  
There exist positive constants $B_0$ and $b_0$ such that,
for all $g\in  L_\infty^{v_\eta} $,  $x\in\state$,  $n\ge 1$,
we have 
\be
\Bigl| \Expect_x \bigl[ \exp (n \langle L_n, F\rangle 
                - n\LA(F)) g (\Phi(n))\bigl] 
                - \cf(x)\cmu(g) \Bigr| 
& \leq &
                \| g \|_{v_\eta} e^{\eta V(x)+B_0-b_0  n}  
                \nonumber\breakMed
\Bigl| \Expect_x \bigl[ \exp (n \langle L_n, F\rangle 
                - n\LA(F))\bigr]
                - \cf(x)\Bigr|
& \leq & 
                        \| F\|_W e^{\eta V(x)+B_0 -b_0 n}  
                \elabel{mainMMETb}
\ee
with $\cf,\cmu,\la(F)$ given as in (i).
\item
If $V$ is bounded on the set $C$ used in 
{\em (DV3)} then we may take $\eta_0=1$.
\end{MyEnumerate}
\end{theorem}

\proof  
Assumption (DV3) combined with \Theorem{v5} 
implies that  $P$ is $v_\eta$-uniform for all $\eta>0$
sufficiently small (when $V$ is bounded on $C$ then (DV3) 
implies $v$-uniformity, so we may take $\eta=1$).

It follows that the inverse $[I-P+ \One\otimes\pi]^{-1}$ exists as
a bounded linear operator on $ L_\infty^{v_\eta}$
\cite[Theorem~16.0.1]{meyn-tweedie:book}.
An application of \Lemma{FrechetLemma}
implies that the kernels $P_f$ converge to $P$ in norm
\[
\lll P- P_f\lll_{v_\eta} \to 0,\qquad \hbox{as}\ \|F\|_W\to 0\, ,
\qquad 0<\eta\le 1\, .
\]
Consequently,  there exists $\epsilon_1>0$ 
such that $[Iz- P_f+ \One\otimes\pi]^{-1}$ is bounded for all $z\in\Co$
satisfying $|z-1|<\epsilon_1$,  and all 
$F\in \LW$ satisfying $\|F\|_{W} \le \epsilon_1$.
  
We  have the explicit representation, writing
$\Delta\eqdef [(z-1)I + I_{1-f} P]$, $H\eqdef [I- P+ \One\otimes\pi] $,
\[
\begin{array}{rcl}
[Iz- P_f+ \One\otimes\pi]^{-1}
&=&
[H + \Delta]^{-1}
\\
&=&
[I+ H^{-1}\Delta]^{-1} H^{-1}\,.
\end{array}
\]
The first term on the right hand side exists 
as a power series in $ H^{-1}\Delta$,
provided 
\begin{equation}
 \lll\Delta \lll_{v_\eta}<  (\lll H^{-1}\lll_{v_\eta})^{-1}\, .
\elabel{PotBdd0}
\end{equation}
Moreover, in this case we obtain the bound,
\begin{equation}
\lll [Iz- P_f+ \One\otimes\pi]^{-1}\lll_{v_\eta} 
\le
	\frac{\lll H^{-1}\lll_{v_\eta}}
	{1- \lll\Delta \lll_{v_\eta}\lll H^{-1}\lll_{v_\eta}} <\infty.
\elabel{PotBdd1}
\end{equation}
 
For any $F\in\LW$   we have the upper bound,
$| F| \le [\lll F\lll_W \delta^{-1}] \delta W$, 
where $\delta>0$ is given in (DV3).
Recalling the definition of the log-generalized principal
eigenvalue functional $\LA$ from \Section{PF}, and
assuming that $\theta\eqdef \lll F\lll_W \delta^{-1}<1$,    
we may apply the convexity of $\Lambda$ (see 
\Lemma{SpectralRadiusConvex}) to  obtain the upper bound,
\begin{equation}
|\Lambda(F)|  
       \le  \Lambda( \theta \delta W) \le  \theta \Lambda(\delta W)
\le \theta b =
\lll F\lll_W \delta^{-1} b
\elabel{PotBdd2}
\end{equation}
where $b$ is given in (DV3).  

From \eq PotBdd2/ we conclude that there is a constant $\epsilon_0>0$ such that 
$\epsilon_0<\half \epsilon_1$, and \eq PotBdd0/ together with the bound
$| \lambda(F) - 1 |< \half \epsilon_1$ hold whenever $ \lll F\lll_W <\epsilon_0$.
For such $F$, it follows that \eq PotBdd1/ holds,
and hence  $P_f$ is $v_\eta$-uniform.  Setting   
$\cH \eqdef [I\lambda(F)- P_f+ \One\otimes\pi]$ we may express 
the eigenfunction and eigenmeasure explicitly as:
\begin{eqnarray*}
\cf&\eqdef&
c_1
 \cH^{-1} \One\, , \qquad c_1 \eqdef 
\Bigl(
\frac{\pi \cH^{-1} \One}{\pi \cH^{-2} \One}\Bigr)
\\
\cmu&\eqdef&
c_2
\pi \cH^{-1}  
\, , \qquad c_2 \eqdef 
\Bigl(
\frac{1}{\pi \cH^{-1}\One}\Bigr)\, .
\end{eqnarray*}
The  remaining results follow as in \cite[Theorem~4.1]{kontoyiannis-meyn:I}.
\qed
 
In order to extend~\Theorem{localMMET} to a non-local result
we   invoke the density condition in~(DV3\upp)~(ii). 
In fact, any such extension seems to 
require some sort of a density assumption.

Recall that, in the notation of \Section{STWR} and
\Section{PF},
we say that the spectrum $\clS$ in $\Lv$ of a
linear operator $\haP:\Lv\to\Lv$ is
{\em discrete}, if 
for any compact set $K\subset\Co\setminus\{0\}$,
$\clS\cap K$ is finite and contains 
only poles of finite multiplicity.
We saw earlier  that condition (DV3+) implies that 
$P^{2T_0+2}$ is $v$-separable.
Next we show in turn that any $v$-separable linear operator 
$\haP$ has a discrete spectrum in $\Lv$.

\archival{In the introduction of 
\cite{kendall-montana:02} the authors mention in
passing
that a continuous transition  density may be expressed as a sum,
$p(x,y)  = \sum s_i(x) p_i(y)$,
with $\{s_i\}$ continuous small functions, and $\{p_i\}$ probability
densities on $\state$.
The main contribution of the paper is to construct counterexamples   to
show
that an analogous representation with measurable small functions does not
hold
in general if the transition density $p$ is only measurable. 
The point of the paper is to demonstrate that a one-step minorization with
a small function
may not exist, even if there is a measurable transition density.   The
representation above
is used in the intro to emphasize how different the situation is for
continuous densities.}

\begin{theorem}
\tlabel{SepImpliesDiscrete}
{\em ($v$-Separability $\Rightarrow$ Discrete Spectrum) }
If the linear operator $\haP\colon\Lv\to\Lv$ is bounded
and $\haP^{T_0}\colon\Lv\to\Lv$ is $v$-separable
for some $T_0\ge 1$, then 
$\haP$ has a discrete spectrum in $\Lv$. 
\end{theorem}

\proof
Assume first that $T_0=1$.
For a given $\epsilon>0$, set $\haP = K+\Delta$ with 
$\lll \Delta\lll_v<\epsilon$, and with $K$ a 
finite-rank operator.  
Write $K=\sum_{i=1}^n s_i \otimes \nu_i$, and for each $z\in\Co$ define
the complex numbers $\{m_{ij}(z)\}$ via
\[
          m_{ij}(z)  =  \langle  \nu_i , [Iz -\Delta]^{-1} s_j \rangle,
											\qquad  1\le i,j \le n\, .
\]
Let $M(z)$ denote the corresponding $n\times n$ matrix, 
and set $\gamma(z)=\det(I-M(z))$. 
The function $\gamma$ is analytic on $\{|z|>\lll \Delta\lll_v\}$ 
because on this domain we have 
\[
 [Iz -\Delta]^{-1} =\sum z^{-n-1}\Delta^n,\qquad  \lll  [Iz -\Delta]^{-1}\lll_v\le (|z|-\lll \Delta\lll_v)^{-1}<\infty .
\]
Moreover, this function satisfies 
$\gamma(z)\to 1$ as $|z|\to\infty$, from which we may 
conclude that the equation $\gamma(z)=0$ has at most a 
finite number of solutions in any compact subset 
of $\{|z|>\lll \Delta\lll_v\}$.  

As argued in the proof of
\Theorem{localMMET}, if $\gamma(z)\neq 0$, then we
have,
\[
\begin{array}{rcl}
 [Iz -\haP]^{-1}     &  =& [(Iz -\Delta) - K]^{-1}   
 \\
   &  =& [Iz -\Delta]^{-1}[I - K[Iz -\Delta]^{-1}]^{-1}    \, .
\end{array}
\]
Conversely, this inverse does not exist when $\gamma(z)=0$.
Recalling that $\epsilon \ge \lll \Delta\lll_v$, we conclude that 
$\clS(\haP)\cap  \{z : |z|>\epsilon \} = \{z : \gamma(z)=0\}.$
The right hand side denotes a finite set, and $\epsilon>0$
is arbitrary.  Consequently,  it follows  
that the spectrum of $\haP$ is discrete.
 
If $T_0>1$ then from the foregoing we may conclude that the
spectrum of  $\haP^{T_0}$ is discrete.  The conclusion
then follows from the identity
\[
\hspace{1.0in}
\bigl[Iz-\haP \bigr]^{-1} 
	= \sum_{k=0}^{T_0-1} z^{-k+T_0-1}
	\Bigl(  \haP^k \bigl[Iz^{T_0} - \haP^{T_0}\bigr]^{-1} \Bigr),
					\qquad z\in\Co.
\hspace{1.0in}
\Box
\]

For each $n\geq 1$, we define the nonlinear 
operators $\Lambda_n$ and $\clG_n$ the space
of real-valued functions $F\in \LWo$, via,
\[
\begin{array}{rcl}
\Lambda_n(F)  &\eqdef&   
\frac{1}{n}
\log \Expect_x
  \bigl [ \exp( n\langle L_n, F \rangle  )
  \bigr ]
  \breakMed
\clG_n(F)  &\eqdef&   
\log \Expect_x
  \Bigl [ \exp
    \Bigl (
        n[\langle L_n, F \rangle - \LA(F)]
    \Bigr )
  \Bigr ] \, ,
  \qquad F\in \LWo,\ x\in\state
  \, .
\end{array}
\]
The following result 
implies that both sequences of operators 
$\{\clG_n\}$ and $\{\Lambda_n\}$ are convergent.
Smoothness properties of the limiting 
nonlinear operators  
are established in Propositions~\ref{t:Lcts} and \ref{t:discreteTimeG}.

\begin{proposition}
\tlabel{discreteTimeGa}
Suppose that~{\em (DV3\upp)} holds with an unbounded
function $W$.  Then there exists a nonlinear operator $\clG\colon\LWo\to\LV$
such that $\cf=e^{\clG(F)}$ is a solution to the multiplicative Poisson equation for
each $F\in\LWo$.  Moreover,  for each $F_{0}\in \LWo$ and $\delta_{0}>0$
we have,
\begin{eqnarray*}
\sup_{\|F-F_{0}\|_{W_{0}}\le \delta_{0}}
\| \clG_n(F) - \clG(F)   \|_{V} &  \to &  0,
\breakMed
\sup_{\|F-F_{0}\|_{W_{0}}\le \delta_{0}}
\| \Lambda_{n}(F) -  \Lambda(F) \|_{V} &  \to & 0\,,\qquad n\to\infty.
\end{eqnarray*} 
\end{proposition}

\proof
Note that the second bound follows from the first.
So, let $\delta_{0}>0$ and $F_{0}\in\LWo$ be given, and consider an arbitrary
$F\in \LWo$ satisfying $\|F-F_{0}\|_{W_{0}}\le \delta_{0}$.
We define $\cF_{n}\eqdef\clG_{n}(F)$
for  $n\ge 0$,  and $\cF = \clG(F)\eqdef \log(\cf)$,  with $\cf$ given in 
\Theorem{mainSpectral}.    
We  show below   that for any $\eta>0$, there exists $b(\eta)<\infty$ such that
for all such $F$,
\begin{equation}\begin{array}{rcl}
|\cF(x)|  & \le& \eta V(x) + b(\eta) \, ,\qquad x\in\state\, ;
\breakMed  |\cF_{n}  (x)|& \le& \eta V(x) + b(\eta) \, ,\qquad x\in\state\, ,\  n\ge 1. 
\end{array}\elabel{cFnBdd}
\end{equation}
Taking this for granted for the moment,  observe that we then have,
for any $r\ge 1$, $n\ge 1$,
\[
\sup_{\|F-F_{0}\|_{W_{0}}\le \delta_{0}}
\||\cF_{n}-\cF|\ind_{C_{V}(r)^{c}}\|_{V}\le  2[ \eta + b(\eta)r^{-1}]\,.
\] 
Moreover,
\Theorem{mainSpectral}   implies that for any $r\ge 1$,
\[
\sup_{\|F-F_{0}\|_{W_{0}}\le \delta_{0}}
\||\cF_{n}-\cF|\ind_{C_{V}(r)}\|_{V}\to 0\,, \quad \hbox{\it exponentially fast as $n\to\infty$,}
\] 
provided we have the uniform bound \eq cFnBdd/.  Putting these two conclusions together, and letting $r\to\infty$
then gives,
\[
\limsup_{n\to\infty}  \sup_{\|F-F_{0}\|_{W_{0}}\le \delta_{0}}
\|\cF_{n}-\cF \|_{V}\le 2\eta.
\]
This then proves the desired uniform convergence, since $\eta>0$ is arbitrary.

\medskip

We now prove the uniform bound \eq cFnBdd/.  We begin with consideration 
of the functions $\{\cF:\|F-F_{0}\|_{W_{0}} \le \delta_{0}\}$, since the corresponding bounds on $\{\cF_{n}\}$ then follow
relatively easily.

We know that $\cf\in L_\infty^{v_{\eta}}$ 
from \Theorem{mainSpectral}.
(If (DV3\upp) holds, then it also holds with 
$V$ replaced by $(1-\eta)+\eta V$ for any $0<\eta<1$.)
This implies that 
$\cF(x)  \le \eta V(x) + \log \|\cf\|_{v_{\eta}}$,
for $x\in\state.$
Hence it remains to obtain a lower bound.  

Let $\tau= \min\{k\ge 1: |\cF(\Phi(k)) |\le r\}$, with $r\ge 1$
chosen so that $\{x:|\cF(x) |\le r\} \in\clB^+$.
The stochastic process below is a positive local martingale,
\[
m(t)=
\exp
                \Bigl(
                        t \langle L_t , (F-\Lambda(F)) \rangle
                \Bigr)\cf(\Phi(t)) \, ,\qquad t\in\nat.
\]
The local martingale property combined with Fatou's Lemma then gives
the bound,     
\[
\cf(x) 
\ge
\Expect_x
        \Bigl[ \exp
                \Bigl(
                        \tau \langle L_\tau , (F-\Lambda(F)) \rangle
                \Bigr)\cf(\Phi(\tau)) 
        \Bigr] \, ,\quad x\in\state,
\]
and then by Jensen's inequality and the definition of $\tau$,
\begin{equation}\begin{array}{rcl}
\cF  (x)
&  \ge &   
\Expect_x
        \Bigl[    \cF  (\Phi(\tau)) 
                   +     \tau \langle L_\tau , (F-\Lambda(F)) \rangle
                      \Bigr] 
                      \breakMed  &  \ge &    
-r 
-
\Expect_x
        \Bigl[      \tau \langle L_\tau , |F+\Lambda(F)| \rangle
                      \Bigr] \, ,\qquad x\in\state.
 \end{array} \elabel{GnA}
\end{equation}
The right hand side is bounded below by $-k_0 (V+1)$ for some finite $k_0$
by (V3) and  \cite[Theorem~14.0.1]{meyn-tweedie:book}.  However, this bound 
can be improved.  Since $F\in \LWo$, and since $W\in\LV$ with $(W_{0},W)$
satisfying \eq Wo/,
we can find, for any $\eta_{0}>0$, a constant
$b_{0}(\eta_{0})$ and a small set $S_{\eta_{0}}$  satisfying
\begin{equation}
|F+\Lambda(F)|  \le  b_{0}(\eta_{0}) \ind_{S_{\eta_{0}}}  +  \eta_{0} V.\elabel{GnB}
\end{equation}
Small sets are \textit{special} 
(see \cite{nummelin:book}),  which implies that
\begin{equation}
\sup_{x\in\state}
\Expect_x
        \bigl[      \tau \langle L_\tau ,  S_{\eta_{0}}   \rangle
                      \bigr]  <\infty.
\elabel{GnC}
\end{equation}
Moreover,  it follows from \cite[Theorem~14.0.1]{meyn-tweedie:book} that
for some $b_{0}<\infty$,
\begin{equation}
\Expect_x
        \bigl[      \tau \langle L_\tau ,   V   \rangle
                      \bigr]  \le b_{0} V(x),\qquad x\in\state.\elabel{GnD}
\end{equation}
Combining the bounds (\ref{e:GnA}--\ref{e:GnD})
    establishes \eq cFnBdd/ for $\cF$.
    
From \eq mainMMET/ in  \Theorem{mainSpectral} we have, for any $\eta >0$,
constants $B_{\eta}<\infty, b_{\eta}>0$ such that, whenever $\|F-F_{0}\|_{W_{0}} \le 1$,
\[
\cF_{n}(x)\le \cF(x) + \log\bigl( 1 +   \exp(\eta V(x)-\cF(x)+B_\eta - b_\eta n)\bigr) \, ,\qquad n\ge 1.
\]
From the forgoing we see that the right hand side is bounded by $2\eta V + b(2\eta)$
for some $b(2\eta)<\infty$ and all $n$.

To complete the proof, we show that a corresponding lower bound holds:
By definition of $\cf_{n}$ and an application of Jensen's inequality we have  for all $n\ge 0$,
\[
 \cf_{n}(x)\cf^{-1}(x)= \cExpect_{x}[\cf^{-1}(\Phi(n))]\ge (\cExpect_{x}[\cf(\Phi(n))])^{-1}  
\]
where the expectation is with respect to the process with transition kernel $\cP_{\cf}$.
On taking logarithms, and appealing to the mean ergodic limit for
the twisted process, for constants $B_{\eta}<\infty, b_{\eta}>0$,
\[
\cF_{n}(x)-\cF(x)\ge - \log(\cExpect_{x}[\cf(\Phi(n))] ) \ge -\log\bigl(\cpi(\cf) +  \exp(\eta V(x) +  B_{\eta}-nb_{\eta})\bigr)\, ,\qquad n\ge 1.
\] 
This together with the bounds obtained on $\cF$ shows that \eq cFnBdd/  does hold.
 \qed

\section{Entropy, Duality and Convexity}
\slabel{convex}

In this section we consider structural properties of 
the operators $\clG$, $\clH$ and the functional
$\Lambda$. As above, we assume 
throughout that $\bfPhi$ satisfies (DV3\upp) with
an unbounded function $W$, 
and we choose and fix an arbitrary function $W_0\in\LW$ 
as in \eq littleOh/. 
Also, throughout this section we restrict 
attention
to real-valued functions in $\LWo$ and
real-valued measures in $\MWo$ 
since one of our goals is to establish 
convexity and present Taylor series
expansions of $\clG$, $\clH$,  and  
$\Lambda$ acting on $\LWo$.  
Recall from \Proposition{gpe} 
that the log-generalized principle eigenvalue $\Lambda$
coincides with the log-spectral radius $\Xi$ on this domain.

The convex dual of the functional  
$\Lambda\colon\LWo\to \Re$ is 
defined for $\mu\in\MWo$ via,
\begin{equation}
    \Lambda^* (\mu)  :=   \sup \{\langle \mu, F \rangle - \Lambda(F):  
	F\in \LWo  \} \, . 
\elabel{I} 
\end{equation}   
A probability measure $\mu\in\MWo$ and a function 
$F\in \LWo$ form a \textit{dual pair} if the 
above supremum is attained, so that
$\Lambda(F) +    \Lambda^* (\mu)  = \langle \mu, F \rangle.$

The main result of this section is a proof that   
$\Lambda^*$ can be expressed in terms of relative entropy
(recall \eq entropyGeneral/) provided that we extend the 
definition to include bivariate measures on
$(\state\times\state,\clB\times\clB).$
Throughout this section we
let $M$ denote a generic function on $\state\times\state$, and
$\Gamma$ a generic measure on $(\state\times\state,\clB\times\clB)$.
The definitions of $\LW$ and $\MW$ are extended as follows:
\be
L^{W}_{\infty,2} 
&\eqdef &
 \displaystyle
 \Bigl\{ 
        M : \|M\|_W:=\sup_{(x,y)\in\state\times\state} 
        \Bigl(\frac{|M(x,y)|}{ [W(x)+W(y)]}\Bigr)  <\infty
                \Bigr\}
	\label{eq:biL}\\
\clM^{W}_{1,2} 
&\eqdef & 
  \displaystyle
  \Bigl\{ 
        \Gamma : \|\Gamma\|_W :=
		\int_{\state\times\state}[ W(x)+W(y)]\, |\Gamma(dx,dy)|<\infty
                \Bigr\}\,.
	\label{eq:biM}
\ee

The following proposition shows that 
consideration of the bivariate chain $\bfPsi$,
\begin{equation}
\Psi(k) = \textstyle \begin{pmatrix} 
                      \Phi(k+1) \\ \Phi(k) 
                       \end{pmatrix},\qquad
k\ge 0,\ \Psi(0)\in\state\times\state\, ,
\elabel{Psi}
\end{equation}
allows us to extend the domain of $\Lambda$ 
to include bivariate functions, and then $\Lambda^*$ 
is defined on bivariate measures via
\be
\Lambda^*(\Gamma)  
\eqdef   \sup_{M\in \LWotwo} (\langle \Gamma , M \rangle - \Lambda(M) ),
                \qquad \Gamma\in \MWotwo.
\elabel{I2} 
\ee
For any univariate measure $\mu$ and transition kernel
$\cP$, we write $\mu\odot\cP$ for the bivariate measure
$\mu\odot\cP(dx,dy):=\mu(dx)\cP(x,dy).$
In particular, \Proposition{PsiAndPhi} shows that if $\bfPhi$ satisfies
(DV3+) with an unbounded $W$, then so does $\bfPsi$.

\begin{proposition}\tlabel{PsiAndPhi}
The following implications hold for any Markov chain $\bfPhi$,
with corresponding bivariate chain $\bfPsi$:
\begin{MyEnumerate}
\item  
If $\bfPhi$ is $\psi$-irreducible, then 
$\bfPsi$ is $\psi_{2}$-irreducible, 
with $\psi_{2}\eqdef \psi\odot P$;
 
\item
If $C$ is a small set for $\bfPhi$, then 
$\state\times C$ is small for $\bfPsi$;
\item  
If $C\in\clB$, $\mu$, and
$T_{0}\ge 1$  satisfy $P^{T_{0}}(y,A)\le \mu(A)$
for $y\in C$, $A\in\clB$,  then on setting $C_{2}=\state\times C$ 
and $\mu_{2}= \mu\odot P$ we have,
\[
P_{2}^{T_{0}+1}((x,y), A_{2})\le \mu_{2}(A_{2}),\qquad
 (x,y)\in C_2,\ A_{2}\in\clB\times\clB,
 \] 
where $P_2$ denotes
the transition kernel for   $\bfPsi$;

\item
If $\nu\in\clM^+$ is small for $\bfPhi$ then 
$\nu_{2}\eqdef \nu\odot P$ is small for $\bfPsi$;
  
\item
Suppose that $\bfPhi$ satisfies the drift condition 
{\em (DV3)}. Then $\bfPsi$ also satisfies the 
following version of {\em (DV3)},
\[
\clH_{2}(V_{2}) \le - \delta W_{2} + b\ind_{C_{2}}\,,
\qquad \hbox{on $ S_{V_2}$}\, ,
\]
where $\clH_{2}$ is the nonlinear generator for $\bfPsi$,
$C_2=\state\times C,$
and 
\[
V_{2} (x,y) = V(y) + \half \delta W(x),\
W_2(x,y) = \half (W(x)+W(y)),\ x,y\in\state.
\]
\end{MyEnumerate}
\end{proposition}

\proof
To prove (i) consider any set 
$A_{2}\in\clB\times\clB$ with $\psi_{2}(A_{2})>0$.
Define 
\[
g(x) = \int_{y\in\state}P(x,dy)\ind_{A_{2}}(x,y)\,, \quad x\in\state.
\]
Then we have $\psi(g)>0$, and hence by $\psi$-irreducibility of $\bfPhi$,
$\sum_{k=0}^{\infty} P^{k}g\, (x) > 0,$ for all $x\in\state$.
It follows immediately that
$\sum_{k=0}^{\infty} P^{k}_{2}\ind_{A_{2}}\, (x,y) > 0,$
for all $x,y\in\state,$
from which we deduce that $\bfPsi$ is $\psi_{2}$-irreducible.
This proves (i), and (ii)-(iv) are similar.

To see (v), observe that under (DV3),
\ben
\hspace{0.6in}
\log P_{2}e^{V_{2}}\, (x,y)  
&  = &   
	\log \int P(y,dz)e^{\half \delta W(y) + V(z)}\\
&  \le &
  \half\delta W(y) + [V(y) -  \delta W(y)+b\ind_C(y)]\\
&  = & V_{2}(x,y) -  \delta W_{2}(x,y) + b \ind_{C_{2}}(x,y)  
			,\qquad
x\in\state,\ y\in S_V.
\hspace{0.6in}\Box
\een

We show in \Theorem{bivariateLambda} that
the convex dual may be expressed as relative 
entropy when $\Gamma$ is a probability measure 
in $\MWotwo$,
\begin{equation}
\Lambda^*(\Gamma)
=H(\Gamma\|\cpi\odot P)
=\int_{\state\times\state}
	 \log\Big( \frac{d \Gamma}{d[\cpi\odot P]}(x,y) \Big)\,
							\Gamma(dx,dy),
\elabel{GammaEntropy}
\end{equation}
where $\cpi$ is the first marginal of $\Gamma$
and $\cpi\odot P$ denotes the bivariate
measure $[\cpi\odot P](dx,dy) = \cpi(dx) P(x,dy)$.
When $\Lambda^*(\Gamma)<\infty$, we show in 
\Lemma{marginals} that  the two marginals agree.  
Consequently, 
$\Gamma$ may be expressed as,
$\Gamma(dx,dy) = [\cpi \odot \cP](dx,dy) =\cpi(dx)\cP(x,dy),$
where $\cP$ is a transition kernel
and $\cpi$ is an invariant measure for $\cP$.

\begin{theorem}
{\em (Identification of $\Lambda^*$ as Relative Entropy) }
\tlabel{bivariateLambda}
Suppose that~{\em (DV3\upp)} holds with an unbounded
function $W$. Then:
\begin{MyEnumerate}
\item
For any probability measure $\Gamma\in\MWotwo$, 
if $\Lambda^*(\Gamma)<\infty$ then
the one-dimensional marginals $\{\Gamma_1,\Gamma_2\}$ agree.  
Consequently, letting $\cpi=\Gamma_1$ denote 
the first marginal of $\Gamma$
we can write, for some transition kernel $\cP$,
\[
\Gamma(dx,dy) = \cpi(dx)\cP(x,dy)\,,
\]
where $\cpi$ is an invariant measure for the transition 
kernel $\cP$.

\item
If $\Lambda^*(\Gamma)<\infty$ for
some probability measure $\Gamma\in\MWotwo$, 
then
\begin{equation}
\Lambda^*(\Gamma) = H(\Gamma \mmid \cpi\odot P)
                \eqdef 
\int_{\state\times\state} 
\log\Bigl( \frac{d \Gamma}{d [\cpi\odot P]  } (x,y) \Bigr) 
                                \Gamma(dx,dy) \, , 
\elabel{entropy}
\end{equation}
where $[\cpi\odot P](dx,dy)\eqdef  \cpi(dx)P(x,dy)$
and $\cpi$ is the first marginal of $\Gamma$.

\item
For any $c>0,$ the set $\{ \Gamma \in \MWotwo : \Lambda^*(\Gamma)\le c\}$
is a bounded subset of $\MWotwo$.  
\end{MyEnumerate}
\end{theorem}

\proof 
Any probability measure $\Gamma$ on $(\state\times\state,\clB\times\clB)$ 
can be decomposed as
$\Gamma(dx,dy) = \cpi(dx)\cP(x,dy),$
where $\cpi$ is the first marginal for $\Gamma$.  We show in
\Lemma{marginals} that the marginals of $\Gamma$ must
agree when $\Lambda^*(\Gamma)<\infty$, and this establishes (i).

Finiteness of $\Lambda^*(\Gamma)$ also implies that 
$\Gamma$ is absolutely continuous with respect to $\pi\odot P$.  
This follows from  
\Proposition{Rate}~(iv) below, applied to the bivariate chain $\bfPsi$.
Consequently, the transition kernel can be expressed,
$\cP(x,dy) =  m(x,y) P(x,dy),$ for $x,y\in\state,$
for some measurable function $m\colon\state\times\state\to[0,\infty]$.

With $M=\log m$, \Proposition{EntropyBdd1}  gives the upper bound,
\[
\Lambda^*(\Gamma) \le \langle \Gamma, M\rangle=H(\Gamma\|\cpi\odot P)\,.
\]
We apply \Proposition{EntropyBdd3} to obtain a corresponding
lower bound: There is a sequence $\{M_k : k\ge 1\}\subset L_{\infty}$
such that $M_k\to M$ point-wise,  $|M_{k}|\le |M|$ for all $k\ge 1$, and 
$\Lambda(M_k) \to \Lambda(M),$ as $k\to\infty$.
Moreover, we have $\Lambda(M)=0$ since $\cP(x,dy)=m(x,y) P(x,dy)$ is
transition kernel for a positive recurrent Markov chain, 
and hence `$1$-recurrent \cite{nummelin:book}.
Consequently,
\[
\Lambda^*(\Gamma) \ge \langle \Gamma, M_k-\Lambda(M_k)\rangle\to 
 \langle \Gamma, M\rangle  ,\qquad k\to\infty\, .
\]
 We thus obtain the identity $\Lambda^*(\Gamma) = \langle \Gamma, M\rangle $,
which is precisely (ii).

Finally, part~(iii) follows from 
\Proposition{Rate}~(iii) combined with
\Proposition{bivariateL}.
\qed

\subsection{Convexity and Taylor Expansions}
\slabel{convexTaylor}

We now return to consideration of  the univariate chain $\bfPhi$,
and establish some regularity and smoothness properties
for the (univariate) functional $\LA$ and the 
nonlinear operators $\clH$ and $\clG$.

We recall the definition of the twisted kernel
$\cP_h$ from \eq cPh/, and 
for any $h\colon\state\to (0,\infty)$
we define the bilinear and quadratic forms,
\begin{equation}
\begin{array}{rcl}
 \llangle F,G \rrangle_h &\eqdef & \bigl[ \cP_h (FG) 
               - ( \cP_h F) ( \cP_h G) \bigr] 
\breakMed
 \clQ_h( F) &\eqdef  &    \llangle F,F \rrangle_h
\end{array}
\qquad F,G\in \LWo\, .
\elabel{langle}
\end{equation}
When $h\equiv 1$ we remove the subscript so that 
$\llangle F,G \rrangle \eqdef   P (FG) - (P F) (P G) $, 
and  $\clQ( F) \eqdef  P(F^2) - (PF)^2$.
It is well-known that 
$\sigma^2(F)\eqdef \pi(\clQ(ZF))$ is equal to   the asymptotic 
variance given in \eq avariance/   
\cite[Theorem~17.5.3]{meyn-tweedie:book}, where  one version of 
the fundamental kernel $Z\colon\Lv\to\Lv$
is given by $Z = [I-P + \One\otimes \pi]^{-1}$;
see \cite{meyn-tweedie:book,kontoyiannis-meyn:I} for details.

The fundamental kernels $\{Z_h\}$ for $\{\cP_h\}$ and the
quadratic forms $\{\clQ_h\}$ determine the second-order
Taylor series expansions for $\LA$, $\clG$ and $\clH$.  
We begin with an examination of $\LA$.

\begin{proposition}
\tlabel{Lcts} 
Suppose that~{\em (DV3\upp)} holds with an unbounded
function $W$.  Then the functional $\LA$ is finite-valued on $\LWo$, and
has the following properties:
\begin{MyEnumerate} 
\item
{\em $\LA$ is strongly continuous:}
For each $F_0\in \LWo$ there exists $B<\infty$, such that for all $F\in\LWo$
satisfying $\| F\|_{W_0} < 1$,
\[
|\Lambda(F_0+F)-\Lambda(F_0)| \le B \|F\|_{W_0};
\]

\item
{\em $\LA\colon\LWo\to\LV$ is smooth:}  For each $F,F_0\in\LWo$,
the function $\Lambda(F_0+a F)$ is analytic as  a function of $a$.
Moreover, we have the second-order Taylor expansion,
\[
\LA(F_0+a F) = \LA(F_0) + a  \pi_g( F)  
                +  \half a^2 \pi_g( \clQ_g( Z_gF) ) 
                + O(a^3) ,  \qquad a\in\Re, 
\]
where  $g=\cf_0 \eqdef e^{\clG(F_0)}$, and $\pi_g$ is 
the invariant probability measure of $\cP_g$.
\end{MyEnumerate}
\end{proposition}

\proof
Part~(i) follows from  \Proposition{LambdaEqualsXi} combined with
\Lemma{FrechetLemma}.

To establish (ii) we note that $\Lambda_{n}(F_0+a F)$ 
is an analytic function of 
$a$ for each initial $x$, and $F_{0},F\in\LWo$.   
\Proposition{discreteTimeGa}
states that this converges to  $\Lambda(F_{0}+aF)$, 
which is convex and hence also
continuous on $\Re$, and the convergence is uniform for $a$ 
in compact subsets of $\Re$.   
This implies that the limit is an analytic function of $a$.

The second-order Taylor series expansion follows as in the proof 
of property~P4 in the Appendix of \cite{kontoyiannis-meyn:I}.
\qed

%


We now consider $\clH$, viewed as a nonlinear operator from $\LWo$ to $\LV$.
\Proposition{Hconvex} establishes smoothness and
pointwise convexity of $\clH$, and
\Proposition{discreteTimeG}
gives analogous results for $\clG$.
See \cite[Chapter~3]{boyvan02a} for related results 
for finite-dimensional positive matrices, and various applications
to optimization.

\begin{proposition}
\tlabel{Hconvex} 
Suppose that~{\em (DV3\upp)} 
holds with an unbounded function $W$.  
\begin{MyEnumerate}
\item
{\em $\clH\colon\LWo\to \LV$ is pointwise convex:}
For any $F_1,F_2\in \LWo$, and for  
any $\theta\in (0,1)$ we have,
\[
\clH (\theta F_1 + (1-\theta) F_2) \leq  
\theta \clH (F_1) + (1 - \theta) \clH (F_2)\, ,
\] 
where inequalities between functions are interpreted pointwise.

\item
{\em $\clH$ is smooth:}
We have the second-order Taylor expansions, for any $F,F_0\in\LWo$,
\[
\clH(F_0+a F) =
      \clH(F_0) +  a \clA_g F 
                + \half a^2 \clQ_g(F) 
                + O(a^3) ,  \qquad a\in\Re, 
\]
where $g=\cf_0:=e^{\clG(F_0)}$ and $\clA_g$ is the generator of $\cP_g$.
\end{MyEnumerate}
\end{proposition}

\proof
We first show that $\clH\colon\LWo\to \LV $. 
To see this, take any $F\in \LWo$.   
Since $W\in  \LV $ and $W_0\in\LW$ satisfies 
\eq Wo/, we can find $b(F)<\infty$ such that 
$ |F|\le  V + b(F)$.  It follows from (DV3) that,
\[\begin{array}{rcl} \log( P e^{F})&  \ge  &     \log(Pe^{-V}) - b(F) \ge -(V-\delta W + b + b(F))
\breakMed \log(Pe^F) & \le &    \log(Pe^V) + b(F) \le V-\delta W + b + b(F),
\end{array}\]
which shows that $\clH(F)\in \LV $.  
Given these bounds, the smoothness result (ii) is a consequence of
elementary calculus.

To establish convexity, we let $H_i = \clH(F_i)$ and $f_i=e^{F_i}$,
so that $P f_i = e^{H_i} f_i$, $i=1,2$.  
An application of  H\"older's inequality
gives the bound,
\[
P (f_1^\theta f_2^{(1-\theta)}) 
 \leq  (P f_1)^\theta(P f_2)^{(1-\theta)} 
 \exp (\theta H_1 + (1-\theta) H_2) f_1^\theta f_2^{(1-\theta)}.
\]
With $F \eqdef   \theta F_1 + (1 - \theta) F_2 
= \log (f_1^\theta f_2^{(1-\theta)})$ we then have
\[ 
\hspace{1.5in}
\clH(F)  =
 \log (P f / f)   
 \le
\theta \clH(F_1) +(1-\theta) \clH(F_2)\, . 
\hspace{1.5in}
\Box
\] 

We can also obtain a Taylor-series approximation for $\clG$, 
but it is convenient
to consider a re-normalization to avoid additive constants.  Define,
\[
\clG_0(F) = \clG(F) - \pi(\clG(F) ),\qquad F\in\LWo.
\]

\archival{The $n$th order derivative of $\clG$ is bounded 
by $V^{n-1}W$ -- see Proposition{logMV3}}

\begin{proposition}
\tlabel{discreteTimeG}
Suppose that~{\em (DV3\upp)} holds with an unbounded
function $W$.    For each $F_0\in\LWo$, $0<\eta \le 1$,
there is  $ \epsilon_0>0$, $ b_0<\infty$, such that
\[
        \| e^{\clG_0(F_0+F)} - e^{\clG_0(F_0)} \|_{v_\eta} \le b_0\|F \|_{W}, 
\]
whenever $ \|F \|_{W_0} <\epsilon_0$.
We have the Taylor series expansion,
\[
\clG_0(F_0+a F) = \clG_0(F_0) + a Z_{\cf_0} F  
                +  \half a^2 Z_{\cf_0} \clQ_{\cf_0}( Z_{\cf_0}F)
                + O(a^3) ,  \qquad a\in\Re, 
\]
where $Z_{\cf_0}$ is the fundamental
kernel for $\cP_{\cf_0}$, normalized so that $\pi Z_{\cf_0} F = 0$,
$F\in \LWo$.
\end{proposition}
 
\proof 
The strong continuity follows from strong continuity of $P_g$
given in \Lemma{FrechetLemma}.

The Taylor-series expansion is established first with $F_0=0$.  Given   
$F\in\LWo$, $a\in\Re$, we let $f_a =\exp(a F)$, and let $\cf_a$
be  the    solution to the eigenfunction equation  given by
\[ 
\cf_a = [I\lambda_a-P_{f_a} + \One\otimes\pi]^{-1}\One
\,.
\]
Under assumption (DV3) alone we have seen in   
\Theorem{localMMET} that  this is an eigenfunction in $\Lv$ for small
$|a|$.  We also have $\cF_a =\log(\cf_a) = \clG_0(F_a) + k(a)$, with
$k(a) = \pi(\cF_a)$.   
In the analysis that follow, our
consideration will focus on $\cF_a$ rather than
$\clG_0(F_a)$ since constant terms will be eliminated through our normalization.

We note that the first derivative may be written explicitly as,
\[
\frac{d}{da} \cg_a
=
[I \lambda_{a}-P_{f_a} + \One\otimes\pi]^{-1}
( \frac{d}{da}\lambda_a I 
	-  I_FP_{f_a})
[I\lambda_{a}-P_{f_a} + \One\otimes\pi]^{-1}\One\,.
\]
Observe that    the derivative is in $L_\infty^{v}$ since
both $I_FP_{f_a}$ and 
$ [I-P_{f_a} + \One\otimes\pi]^{-1} $ are bounded linear 
operators on $L_\infty^{v}$.  Similar conclusions hold for all 
higher-order derivatives.

We define the twisted kernel as above,
\[
\cP_a (x, A)  
\eqdef
\cP_{\cf_a}(x,A)
=  \frac{\int_A P (x, dy) \cf_a (y)}{
P \cf_a (x)} 
\qquad x\in\state,\ A\in\clB.
\]
As in \cite{kontoyiannis-meyn:I} we may verify that the
function $\haF_a  
= \frac{d}{da} \cF_a$ 
is a solution to Poisson's equation,
\[
\cP_a \haF_a = \haF_a - F + \pi_a(F),
\qquad 
\pi_a(F)=  \frac{d}{da} \Lambda(a F)\, ,
\]
where $\pi_a$ is invariant for $\cP_a$.  Setting $a=0$
gives the first term in the Taylor series expansion for $\clG_0$.

To obtain an expression for the second term we differentiate
Poisson's equation:
\begin{equation}
 \frac{d}{da} 
\bigl(\haF_a - F +  \frac{d}{da} \Lambda(a F) \bigr)
=
 \frac{d}{da} \bigl(
\cP_a \haF_a \bigr)
=
\bigl(
 \frac{d}{da} 
\cP_a  \bigr) \haF_a
+
\cP_a \haF^{(2)}_a .
\elabel{SecondFish2}
\end{equation}
We wish to compute the second derivative,
$\haF^{(2)}_a = \frac{d^2}{da^2} \log(\cf_a),$
which requires a formula for the derivative of  $\cP_a$:
For any $G\in \LV$,
\begin{equation}
\begin{array}{rcl}
 \frac{d}{da}\Big(\cP_a G\Big)
&=&\displaystyle
\frac{P (\cf_a' G) P \cf_a  -  P(\cf_a G) 
                       ( P\cf_a')}{ (P \cf_a)^2 }
\breakBig
&=&\displaystyle
\cP_a (\haF_a G) 
 - \bigl( \cP_a G \bigr)\bigl( \cP_a \haF_a \bigr) 
	= \llangle \haF_a, G\rrangle_{\cf_a}\, .
\end{array}
\elabel{SecondFish1}
\end{equation}
Letting $H_a =   \llangle \haF_a, \haF_a \rrangle_{\cf_{a}}$,
the identities  \eq SecondFish2/ and  \eq SecondFish1/ then give,
\begin{equation}
\cP_a \haF^{(2)}_a  = \haF^{(2)}_a  
- H_a +\Lambda''(aF)\, .
\elabel{CLTfish}
\end{equation}
Letting $Z_a$ denote the fundamental kernel for
$\cP_a$ we conclude that
\[
 \haF^{(2)}_a  - \pi( \haF^{(2)}_a) 
		= Z_a H_a  
              	= Z_a\llangle \haF_a,\haF_a\rrangle_{\cf_a} \, .
\]
Evaluating all derivatives at the origin
provides the quadratic approximation for $\clG_{0}$, 
\[
\clG_{0}(a F) =  a Z F  + \half a^2 Z [ \llangle \haF, \haF\rrangle]
                          + O(a^3)
\]
where $Z$ is the fundamental kernel for $P$, normalized so that $\pi Z=0$,
and $\haF=ZF$.

To establish the Taylor-series expansion at arbitrary $a_0\in\Re$ we 
repeat the above arguments, applied to the Markov chain with transition
kernel $\cP_{a_0}$.  This   satisfies (DV3\upp) with   
$\cV=c+V-\cF_{a_0}$ for sufficiently large $c>0$, by
\Proposition{twisted}.
\qed

\subsection{Representations of the Univariate Convex Dual}
\slabel{uniDual}

The following result provides bounds on the (univariate)
convex dual functional $\Lambda^*$, and gives 
some alternative representations:

\begin{proposition}
\tlabel{Rate} 
Suppose that~{\em (DV3\upp)} holds with an unbounded
function $W$.   
Then, for any probability measure $\mu\in \MWo$:
\begin{MyEnumerate}
\item
	$\displaystyle
    	\Lambda^* (\mu)  =  
	\sup \{\langle \mu, F \rangle - \Lambda(F):  F\in L_{\infty} \ \hbox{and} \ \cF\in L_\infty  \}$.
\item
	$\displaystyle \Lambda^*(\mu) = \sup\Bigl\{ \langle \mu, -\clH(H)\rangle  
	: H\in L_\infty  \Bigr\}$.
\item 
	There exists  $\epsilon_0>0$, independent of $ \mu\in \MWo$, such that
	\[
	\Lambda^*(\mu) \ge \epsilon_0 \Bigl( \frac{\|\mu - \pi \|_{W_0}^2}%
                                { 1+ \|\mu - \pi \|_{W_0}} \Bigr)
                                 \,, \qquad \mu\in\MWo\, .
	\]
\item
If $\mu$ is not absolutely continuous with respect to $\pi$, 
then $\Lambda^*(\mu)=\infty$.
\end{MyEnumerate}
\end{proposition}

The proof is provided after the following bound.
\begin{lemma}\tlabel{BddCheck}
Suppose that~{\em (DV3\upp)} holds with an unbounded
function $W$.   
Then, $\cF\in L_{\infty}$ provided  the following conditions hold:
$F\in\LV$;   $\Lambda(F)=0$;
and     $F=F\ind_{C_{V}(r)}$ for some $r\ge 1$. \end{lemma}

\proof 
From the local martingale property  we   have,
\[\begin{array}{rcl}\cF(x)   &  = &   
\log \Expect_x\Bigl [
 \exp\Bigl(\sum_{i=0}^{\tau_{C_{V}(r)} -1}   F(\Phi(i))       \Bigr)\cf(\Phi( \tau_{C_{V}(r)} )    )\Bigr] 
\breakMed  &  = &   
F(x)
+
\log
\Expect_x\bigl [
 \exp\bigl( \cf(\Phi( \tau_{C_{V}(r)} )   \bigr )\bigr] .
 \end{array}\]
 This then gives the bound,   
$ \|\cF\|_{\infty} \le \|F\|_\infty + \|\cF\ind_{  C_{V}(r ) } \|_\infty   
<\infty $.\qed

\noindent{\sc Proof of \Proposition{Rate}.}
For any $F\in\LWo$, and any $r\ge 1$ we write,
$F_{r}= \ind_{C_{V}(r)} [F- \gamma_{r}],$
where $\gamma_{r}\in\Re$ is chosen so that $\Lambda(F_{r})=0$.  Its
existence follows from \Proposition{Lcts}.  

From \Proposition{LambdaTight} we can show that $\gamma_{r}\to\Lambda(F)$, and then
also that $\Lambda(F_{r})\to 0$ as $r\to\infty$.
Consequently, $\Lambda^* (\mu)  =  
	\sup \{\langle \mu, F_{r} \rangle 
	- \Lambda(F_{r}):  F\in \LWo,\ r\ge 1  \},$
and \Lemma{BddCheck} implies that 
$\cF_{r}\in L_{\infty}$ for each $r$, which completes the proof of (i).

  
Part (ii) is essentially a reinterpretation of  (i):    
From the   equation  $\clH(\cF) = - F + \Lambda(F)$ and part~(i) we   obtain
the upper bound,
\[
\begin{array}{rcl}
\Lambda^*(\mu) & =  &
 \sup \{\langle \mu, F \rangle - \Lambda(F):  \cF\in L_\infty \} 
\breakMed
&\le&
 \sup \{\langle \mu, -\clH(\cF) \rangle:  \cF\in L_\infty  \}  
\breakMed
&\le &
 \sup \{\langle \mu, -\clH(G) \rangle:  G\in L_\infty  \}  .
\end{array}
\]  
Conversely, for any function 
$G\in L_\infty$, the function $F \eqdef  -\clH(G)$
satisfies $\Lambda (F) = 0$, $F\in L_\infty$. 
This gives the desired lower bound,
$\Lambda^*(\mu)  
        \ge \langle \mu, F \rangle  
         =  \langle \mu,-\clH(G)\rangle,$
for $G\in L_\infty.$

\archival{Check tightness of inequality below -- If OK, then this will
give Frechet differentiability of Lambda, since we have a uniform bound 
on the second derivative. }

Result (iii) is obtained from the mean value theorem, 
justified by \Proposition{Lcts}:  For any $F\in\LWo$, $ \epsilon\ge 0$,
there is $0\le \tilde\epsilon\le \epsilon$ such that
$\Lambda (\epsilon F)  =  \epsilon \pi (F) + \half 
\epsilon^2 \Lambda'' (\tilde\epsilon F).$
Let $B_{0}=\sup \{  \Lambda''(\epsilon G) :
\|G\|_{W_{0}} \le 1,\;0\leq\epsilon\leq 1 \}$.
Note that $B_0<\infty$ by the Lemma following the proof.
Then, whenever
$\|F\|_{W_{0}}\le 1$, $\epsy\le 1$, we have
$\Lambda (\epsilon F)  \le  \epsilon \pi (F) + \half  B_{0} \epsilon^2.$
The definition of the convex dual then gives,
\[
\begin{array}{rcl}
\epsilon \mu (F)  =  \langle \mu, \epsilon F \rangle  & \leq & \Lambda^*(\mu) 
+ \Lambda (\epsilon F) \\ 
& \leq & \Lambda^*(\mu) + \epsilon \pi (F) + \half  B_{0} \epsilon^2  \, ,
\end{array}\]
and since this holds for any 
$\|F\|_{W_{0}}\le 1$, we have the absolute bound,
\[
\begin{array}{rcl}
| \mu (F) - \pi (F) | & \leq &  \frac{1}{ \epsilon} \Lambda^*(\mu) 
+ \half B_{0} \epsilon^2   \, ,\qquad    \|F\|_{W_{0}}\le 1\, .
\end{array}
\]
Letting $\epsilon   =  \sqrt{\Lambda^*(\mu)}$ we obtain
\[
\| \mu - \pi \|_{W_0}
=
\sup_{\|F\|_{W_0}\le 1} | \mu (F) - \pi (F) |
 \leq    
  \sqrt{\Lambda^*(\mu)}
  +   \half B_{0} \Lambda^*(\mu)\,, \qquad |\Lambda^*(\mu)|<1,
\]
which implies the desired     lower bound on $\Lambda^*$.

To prove (iv), write $\mu = p \mu_0 + (1-p)\mu_1$ 
where $\mu_0,\mu_1$ are probability measures on 
$(\state,\clB)$ such that
$\mu_1\prec\pi$ is absolutely continuous 
and  $\mu_0$ is singular with respect to $\pi$.
Let $S$ denote the support of $\mu_0$.  We have $\Lambda(F)=0$
whenever $F\in L_\infty$ is supported on $S$, and hence
\ben
\Lambda^* (\mu)  
&\ge&
 	\sup \{\langle \mu, F \rangle - \Lambda(F):  
	F\in \LWo,\ F=\ind_S F  \}\\
&=&
	p \sup \{\langle \mu_0, F \rangle  :  F\in \LWo,\ F=\ind_S F  \},
\een
which is infinite, as claimed.
\qed

\begin{lemma}
$B_{0}=\sup \{  \frac{d^2}{da^2} \Lambda(a G) :
\|G\|_{W_{0}} \le 1,\;0\leq a\leq 1 \}<\infty$. 
\end{lemma}

\proof  (sketch)  
Let $\cP_a = \cP_{\cg_a}$ and let $ \pi_a$ denote the 
invariant distribution for given $\|G\|_{W_{0}} \le 1$, and $ a\in[0,1]$. 
We let  
$Z_a$   the fundamental kernel for
$\cP_a$, normalized  so that $\pi_a(Z_a G) = \pi_a(G)$, and  we 
let $\haG_a = Z_a G$.  \Proposition{Lcts} then  gives the representation,
\[
\Lambda''(aG) = 
\pi_a( \clQ_a( Z_a G) ) = \pi_a\bigl(  P_a(\haG_a^2) - (P_a \haG_a)^2 \bigr).
\]
The proof is completed on showing that
\[
\sup \|\pi_a\|_v <\infty,\quad \sup  \| \haG_a\|_V<\infty,
\]
where the supremum is over all $a$ and $G$ in this class.
This follows from the arguments above -- see in particular \eq cFnBdd/ and
the surrounding arguments.
\qed


In the following proposition we give another
characterization of dual pairs $(\mu,G)$ for $\Lambda^*$.

\begin{proposition}
\tlabel{Iinv}
Suppose that~{\em (DV3\upp)} holds with an unbounded
function $W$.   We then have:
\begin{MyEnumerate}
\item 
For any $H\in\LWo$, $\pi \bigl( \clH(H) \bigr)  \geq 0,$
with equality if and only if $\clH(H)=0$,  in 
which case $H$ is  constant a.e.\ $[\pi]$.  

\item 
If $\mu\in\MWo $ is not invariant under $P$
then there is $H\in L_{\infty}$ satisfying 
$\mu \bigl( \clH(H) \bigr) <0$. 

\item 
Suppose that $\mu\in \MWo$, and that there exists $G\in \LWo$ 
satisfying,
\[
\Lambda^*(\mu) 
=  
\langle \mu,-\clH(G) \rangle 
=
\sup\Bigl\{ 
  \langle  \mu, -\clH(H)\rangle  
          : H\in \LWo\Bigr\} \,.
\]
Then $\mu$ is invariant under the twisted
kernel $\cP_g$.
\end{MyEnumerate} 
\end{proposition}
    
\proof
The first result is simply Jensen's inequality: 
\[
\begin{array}{rcl}
\displaystyle
\pi  \bigl( \clH(H) \bigr)
   & = & 
\displaystyle
\int \log 
\Bigl( 
\Expect_x \bigl[ 
   \exp \bigl( H (\Phi(1)) - H (\Phi(0)) \bigr) 
   \bigr]  
\Bigr) \pi (dx) 
\breakBig
 & \geq & 
\Expect_\pi \bigl[  H (\Phi(1)) - H (\Phi(0)) \bigr]   = 0.
\end{array}
\]
If equality holds, it then follows that $e^{H} $ is constant a.e.\ $[\pi]$.  

To prove (ii) let  
$F (x) = \epsilon [\ind_A-\gamma_{\epsilon} \ind_{B}] $, with $\epsilon>0$,
$A,B\in\clB^+$ small sets such that   $\sup_{A\cup B}V(x)$ is finite, 
and $\gamma_\epsilon >0$ is chosen
so that $\Lambda(F)=0$.  The function $H \eqdef \cF= \clG(F)$ is then bounded, by \Lemma{BddCheck}.   Moreover,
\[
\mu \bigl( \clH(H) \bigr) =   -\mu(F)  = -\epsilon [\mu (A) - \gamma_{\epsilon} \mu(B) ].
\]
Under~(DV3\upp) we may apply \Proposition{Lcts} 
to justify the Taylor series expansion,
\[
0 = \Lambda(F) = \epsilon [\pi(A) - \gamma_\epsilon \pi(B)]  + O(\epsilon^{2}),
\]
which gives  $\gamma_{\epsilon } =\pi(A)/\pi(B) + O(\epsilon)$.   
Choosing $A,B$ so that $\mu (A)/\mu(B) > \pi (A)/\pi(B)$ we see that this function $H$ 
satisfies the desired bound for $\epsilon>0$ sufficiently small.

We now prove (iii).  Applying \Proposition{Rate}~(ii), the 
convex dual $\cLambda^*$ for the kernel $\cP_g$ may be expressed
as
\[
\cLambda^*(\mu)  \eqdef   
- \Bigl(
     \inf_{H\in\LWo} 
\Bigl \langle \mu, \log \Bigl(\frac{\cP_g h}{h}\Bigr)  \Bigr\rangle 
        \Bigr) \,.
\]
For  any $H\in\LWo$ set $H' = H+G$ so that,
\[
\begin{array}{rcl}
\cLambda^* (\mu)
& = & 
- \Bigl( \inf_{H'\in\LWo} 
\Bigl\langle \mu, \log 
\Bigl( \frac{g}{ h'} \frac{Ph'}{P g}\Bigr) \Bigr  \rangle  \Bigr)
 \breakMed
& = &
 - \Bigl( \inf_{H'\in\LWo}\Bigl\langle \mu, \log \Bigl(
 \frac{Ph'}{h'}\Bigr)\Bigr \rangle \Bigr)   +  
\Bigl\langle \mu, \log \Bigl(\frac{P g}{g} \Bigr)\Bigr\rangle =
0.
\end{array}\]
Thus $\mu$ is invariant for $\cP_{g}$, by (i). 
\qed

\subsection{Characterization of the Bivariate Convex Dual}
\slabel{biDual}

We now turn to the case of bivariate functions
and measures.

Given any function of two variables 
$M\colon\state\times\state\to \Re$,
we let $m=e^M$ and extend the definition 
of the scaled kernel in \eq If/ via,  
\[
P_m(x,dy) \eqdef  m(x,y) P(x,dy)\, ,\qquad x,y\in\state\, .
\]
The following result  shows that the spectral radius of this kernel
coincides with that defined for the bivariate chain $\bfPsi$.  
The proof is routine.

\begin{proposition}
\tlabel{bivariateL}
Suppose that $P_m$ has finite spectral radius $\lambda_m$ in 
$v_\eta$-norm for all sufficiently small $\eta>0$.  
Let $P_2$ denote
the transition kernel for the bivariate chain $\bfPsi$.
\begin{MyEnumerate}
\item 
$I_m P_2$ has the same spectral radius in $v_{\eta2}$-norm
for sufficiently small $\eta>0$, with $v_{\eta2}(x,y) = 
\exp(\eta[V(y) + \half \delta W(x)])$.

\item 
If $P_m$ has an eigenfunction $\cf$, then 
$I_m P_2$ also possesses an eigenfunction given by,
\[
        \cf_2(x_1,x_2) = m(x_1,x_2) \cf(x_2).
\]
\end{MyEnumerate}
\end{proposition}

For a Markov process with transition kernel $P$ satisfying
(DV3+), we say that $M$ and $\tilM$ are \textit{similar} 
if there exists $H\in\LV$ such that
\[
\tilM(x,y)= M(x,y) +H(x)-H(y)\qquad a.e.\ 
(x,y)\in\state\times\state\ [\pi\odot P]\, .
\]
The function $M$ is called \textit{degenerate} if it is similar to 
$\tilM\equiv 0$.  The log-generalized principal 
eigenvalues  agree ($\Lambda(M)=\Lambda(\tilM)$) whenever 
$M,\tilM$ are similar.  This is the basis of  the following two lemmas.

\begin{lemma}
\tlabel{marginals}
Suppose that~{\em (DV3\upp)} holds with an unbounded
function $W$.  
If  $\Gamma\in\MWotwo$ is a probability measure with
$\Lambda^*(\Gamma)<\infty$,
then $\Gamma\prec \pi\odot P$, 
and the one-dimensional marginals of $\Gamma$ agree.
\end{lemma}

\proof
The conclusion that $\Gamma\prec \pi\odot P$ 
follows from \Proposition{Rate}~(iv).

For any $M\in L_\infty(\state\times\state)$, $H\in L_{\infty}$, we have
$\Lambda(M) = \Lambda(\tilM)$, where $\tilM(x,y) \eqdef  M(x,y) +H(x)-H(y)$.
Hence, for all such $M,H$,
$$\Lambda^*(\Gamma) 
\ge \langle \Gamma , \tilM\rangle - \Lambda(\tilM) 
= \langle \Gamma , M \rangle 
                 - \Lambda(M) 
+ \langle \Gamma_1 , H\rangle - \langle \Gamma_2 , H\rangle 
$$
where $\Gamma_1$ and $\Gamma_2$ denote the two marginals.
If $\Gamma_1\neq \Gamma_2$ it is obvious that the 
right hand side cannot be bounded in $H$.
\qed

\begin{lemma}
\tlabel{zeroVar}
Suppose that~{\em (DV3\upp)} holds with an unbounded
function $W$.  
Suppose moreover that $M\in L_{\infty,2}^V$, and that the 
asymptotic variance of the partial sums
$\sum_{k=0}^{n-1} M(\Phi(k),\Phi(k+1)),\; n\ge 1,$
is equal to zero. Then the function $M$ is degenerate.
\end{lemma} 

\proof
Applying \cite[Proposition~2.4]{kontoyiannis-meyn:I}
to the bivariate chain $\bfPsi$ with transition kernel $P_2$,
we can find $\haM$ such that 
\[
\haM(\Phi(k),\Phi(k+1)) - \haM(\Phi(k-1),\Phi(k)) 
= - M(\Phi(k-1),\Phi(k)) + \pi_2(M)\qquad a.s.\ [\Prob_\pi],\ k\ge 1,
\]
where $\pi_2=\pi\odot P$ is the invariant probability measure for $P_2$.
Since $\Phi(k+1)$ is conditionally independent of $\Phi(k-1)$ given $\Phi(k)$,
it follows that $\haM$ does not depend on its first variable.  Thus 
we can find $\haF\in\LV$ satisfying
\[
\haF(\Phi(k+1)) - \haF(\Phi(k)) = - M(\Phi(k-1),\Phi(k)) + \pi_2(M)
\qquad a.s.\ [\Prob_{\pi}],\ k\ge 1\,,
\]
therefore, $M$ is similar to the constant function $\pi_2(M)$:
\[
M(x,y) = \pi_2(M) +G(x)-G(y),\qquad a.e.\ \pi\odot P\, ,
\]
with $G(x) = P\haF\, (x)$.
\qed

\newpage

\begin{theorem}
{\em (Identification of Dual Pairs)}
\tlabel{dualPairs}
Suppose that~{\em (DV3\upp)} holds with an unbounded
function $W$. 
\begin{MyEnumerate}
\item 
Assume that $M\in\LWotwo$ and $\Gamma\in \MWotwo$ are given,
such that $\Lambda^*(\Gamma)<\infty$ and
$(M,\Gamma)$ is a dual pair, i.e.,
$\langle \Gamma,M\rangle = \Lambda(M)+ \Lambda^*(\Gamma).$
Define $M_0$ as  the 
Radon-Nikodym derivative,
\[
M_0(x,y) = \log\Bigl( \frac{d \Gamma}{d [\cpi\odot P]  } (x,y) \Bigr) 
 \,  \qquad x,y\in\state\, ,
\]
where $\cpi$ is a marginal of $\Gamma$ (see 
\Lemma{marginals}).
Then, the function $M_0$ is similar to $M-\Lambda(M)$,
\[
M_0(x,y) = M(x,y)  -\Lambda(M)-\cF(x) + \cF(y) \, ,
\]
where $\cF=\log(\cf)$, with $\cf$ equal to an eigenfunction for $P_m$,
with eigenvalue $\lambda(M)$.

\item  
Conversely, suppose that $\Gamma\in\MWotwo$ is given,
satisfying $\Gamma\prec [\pi\odot P]$, and suppose 
that its one-dimensional marginals agree.
Consider the decomposition,
$\Gamma(dx,dy) = [\cpi\odot \cP](dx,dy),$
where $\cpi\eqdef  \Gamma_{1} = \Gamma_{2}$ is the 
(common) first marginal of $\Gamma$ on 
$(\state,\clB)$, and $\cP$ is a transition kernel. 
Let
\[
M(x,y) = \log\Bigl( \frac{d \Gamma}{d [\cpi\odot P]  } (x,y) \Bigr) 
 \,  \qquad x,y\in\state\, .
\]
If $M\in L^{W_{0}}_{\infty,2}$, then $\Lambda^*(\Gamma)$ is finite and $(\Gamma,M)$ 
is a dual pair. 
\end{MyEnumerate}
\end{theorem}

\proof
Part (i) is a bivariate version of \Proposition{Iinv}:
We know that $\Gamma$ is an invariant measure for a bivariate
process, whose one-dimensional transition kernel is of the form,
\[
\cP_m(x,dy) = e^{M(x,y)  -\Lambda(M)-\cF(x) + \cF(y) } P(x,dy).
\]
Invariance may be expressed as follows:
\[ \Gamma(dy,dz)  
=
\int_{x\in\state} \Gamma(dx,dy)\cP_m(y,dz) 
\, ,\qquad y,z\in\state.
\]
Since $\Gamma$ has equal marginals, denoted $\cpi$,  this 
identity may be expressed,
\[
\cpi(dy ) \cP(y,dz)=\cpi(dy ) \cP_m(y,dz) 
\, ,\qquad y,z\in\state,
\]
which is the desired identity in (i).

To prove (ii), let $\cLambda(\varble)$ denote the  functional defining the
log-generalized principal eigenvalue  for the transition kernel 
$\cP=P_{m}$.  \Proposition{twisted} gives,  
$\cLambda(N) = \Lambda(N+M) - \Lambda(M),$
for any $N\in L^{W_{0}}_{\infty,2}$.
We can then write,
\[
\begin{array}{rcl}
\Lambda^*(\Gamma) &=&
        \sup_{N\in \L} \Bigl( \langle \Gamma, N\rangle - \Lambda(N) \Bigr)
\breakMed
&=&
        \sup_{N\in \L} \Bigl( \langle \Gamma, N+M\rangle - \Lambda(N+M) \Bigr)
\breakMed
&=&
        \sup_{N\in \L} \Bigl( \langle \Gamma, N\rangle 
                +\langle \Gamma,M\rangle 
                - \cLambda(N) - \Lambda(M) \Bigr)
\breakMed
&=& \cLambda^*(\Gamma) + \langle \Gamma,M\rangle  - \Lambda(M)\,.
\end{array}
\]
We have $\cLambda^*(\Gamma)=0$ by \Proposition{Iinv}, and consequently 
$\langle \Gamma,M\rangle = \Lambda(M) +\Lambda^*(\Gamma)$.  
This shows that $(M,\Gamma)$ is a dual pair.
\qed

\section{Large Deviations Asymptotics}
\slabel{ldp}

In this section we use the multiplicative mean
ergodic theorems of~\Section{met}
and the structural results of~\Section{convex}
to study the large deviations properties
of the empirical measures $\{L_n\}$ induced
by the Markov chain $\bfPhi$ on $(\state,\clB)$;
recall the definition of $\{L_n\}$ in \eq OccEmp/.

As in the previous section, we also assume throughout
this section that the Markov chain $\bfPhi$ 
satisfies~(DV3+) with an unbounded function
$W$, and we choose and fix a function 
$W_0:\state\to[1,\infty)$ in $\LW$
as in \eq littleOh/. Our first result,
the large deviations principle (LDP)
for the sequence of measures $\{L_n\}$,
will be established in a topology finer
(and hence stronger) than either the topology
of weak convergence, or the $\tau$-topology.
As described
in the Introduction, we consider the
$\tau^{W_0}$-topology on the space 
$\clM_1$ of probability measures
on $(\state,\clB)$, defined by the
system of neighborhoods (\ref{eq:nbhds}).

Since the map 
$(x_1,\ldots,x_n)\mapsto\smalloneovern\sum_{i=1}^n\delta_{x_i}$
from $\state^n$ to $\clM_1$ may not
be measurable with respect to the natural
Borel $\sigma$-field induced by the 
$\tau^{W_0}$-topology on $\clM_1$,
we will instead consider the 
(smaller) $\sigma$-field
$\clF$, defined as the smallest 
$\sigma$-field that makes all the maps 
below measurable:
\be
\nu\mapsto\int F\,d\nu,\qquad\mbox{for real-valued}\; F\in \LWo .
\label{eq:clF}
\ee

\begin{theorem}
\tlabel{ldp}
{\em (LDP for Empirical Measures) }
Suppose that $\bfPhi$ satisfies {\em (DV3\upp)} 
with an unbounded function $W$. Then, for any
initial condition $\Phi(0)=x$, the sequence of
empirical measures $\{L_n\}$ satisfies the LDP
in the space $(\clM_1,\clF)$
equipped with the $\tau^{W_0}$-topology,
with the good, convex rate function
\be
I(\nu):=\inf{\cP} H(\nu\odot \cP\|\nu\odot P)
\label{eq:1dI}
\ee
where the infimum is over all transition
kernels $\cP$ for which $\nu$ is an invariant
measure, and $\nu\odot \cP$ denotes
the bivariate measure 
$[\nu\odot \cP](dx,dy)\eqdef \nu(dx)\cP(x,dy)$
on $(\state\times\state, \clB\times\clB)$:
Writing $\mu_{n,x}$ for the law of the empirical
measure $L_n$ under the initial condition $\Phi(0)=x$,
then for any $E\in\clF$,
\ben
-\inf_{\nu\in E^o}I(\nu)
&\leq&
\liminf_{n\to\infty}\smalloneovern\log\mu_{n,x}(E)\\
&\leq&
\limsup_{n\to\infty}\smalloneovern\log\mu_{n,x}(E)
\;\leq\;-\inf_{\nu\in\bar{E}}I(\nu)\,,
\een
where $E^o$ and $\bar{E}$ denote the interior
and the closure of $E$ in the $\tau^{W_0}$ topology,
respectively.
\end{theorem}

The proof is based on an application of the 
Dawson-G\"{a}rtner projective limit theorem 
along the same lines as the proof of Theorem~6.2.10
in \cite{dembo-zeitouni:book}. The main two
technical ingredients are provided by,
first, the multiplicative mean
ergodic theorem \Theorem{mainSpectral}~(iii)
which, as noted in (\ref{eq:mmetL}),
shows that the log-moment generating functions
converge to $\LA$. And second, by the regularity
properties of $\LA$ and the identification
of $\LA^*$ in terms of relative entropy,
established in \Section{convex} and \Section{app:Lambda}
of the Appendix.

As in \Section{convex}, in order to identify the
rate function for the LDP we find it easier to consider 
the bivariate chain $\bfPsi$. 
Recall the bivariate extensions of our earlier
definitions from equations (\ref{eq:biL}), 
(\ref{eq:biM}), \eq Psi/ and \eq I2/.

\vspace{0.1in}

\noindent{\sc Proof of \Theorem{ldp}.}
We begin by establishing an LDP for $\bfPhi$ 
with rate function given by $\Lambda^*$.
Recall that  \Proposition{discreteTimeGa} gives
\be
\Lambda_n(F)
\eqdef
\smalloneovern\log
\Expect_x 
  \Bigl [ \exp
    \Bigl ( n \langle L_n, F \rangle 
    \Bigr ) 
  \Bigr ]
	\to \LA(F)\,,
\;\;\;\;n\to\infty.
\label{eq:weakMMET}
\ee
In order to apply the projective limit theorem
we need to extend the domain of the convex dual 
functional $\LA^*$ as follows. For probability
measures $\nu\in\MWo$, $\LA^*(\nu)$ is defined
in \eq I/, and the same definition applies when
$\nu$ is a probability measure not necessarily
in $\MWo$. More generally, let $L'$ denote 
the algebraic dual of the space
$L=\LWo$, consisting of all linear functionals 
$\Theta:L\to\RL$, and equipped with the
weakest topology that makes the functional
$$ \Theta \mapsto \Theta(F)=(\Theta,F): L'\to\RL$$
continuous, for each in $F\in\LWo$.
Note that each probability measure $\nu$ on 
$(\state,\clB)$ induces a linear functional 
$\Theta_\nu:L\to\RL$ via 
$$(\Theta_\nu,F)=\langle\nu,F\rangle=\int Fd\nu.$$
Therefore, we can identify the space of probability 
measures $\clM_1$ with the corresponding 
subset of $L'$, and observe that the induced topology 
on $\clM_1$ is simply the $\tau^{W_0}$-topology.

Next, extend the definition of 
$\LA^*$ to all $\Theta\in L'$ via
\begin{equation}
\LA^*(\Theta)=\sup\{(\Theta,F)-\LA(F)\,:\,F\in\LWo\}\,,
\elabel{LambdaStar}
\end{equation}
and observe that \cite[Assumption~4.6.8]{dembo-zeitouni:book}
is satisfied by construction (with $\clW=L=\LWo$,
$\clX=L'$ and $\clB=\clF$), and that by \Proposition{Lcts}
the function $\LA(F_0+\alpha F)$ is Gateaux differentiable.
Therefore, we can apply the Dawson-G\"{a}rtner projective
limit theorem \cite[Corollary~4.6.11~(a)]{dembo-zeitouni:book}
to obtain that the sequence of empirical measures $\{L_n\}$
satisfy the LDP in the space $L'$
with respect to the convex, good rate function 
$\LA^*$. Moreover, since by \Proposition{goodconvexdual}
we know that 
$\LA^*(\Theta)=\infty$ for $\Theta\not\in\clM_1$, 
we obtain the same LDP in the space $(\clM_1,\clF)$, 
with respect to the induced topology, namely,
the $\tau^{W_0}$-topology; see, e.g., 
\cite[Lemma~4.1.5]{dembo-zeitouni:book}.

Next note that, in view of \Proposition{PsiAndPhi},
the bivariate chain $\bfPsi$ also satisfies the same LDP. 
But in this case, we claim that can express $\Lambda^*(\Gamma)$
for any bivariate probability measure $\Gamma$ as follows:
\ben
\Lambda^*(\Gamma)=
\begin{cases}
		H(\Gamma\|\Gamma_1\odot P)\,,
                \quad & \hbox{if the two marginals 
			$\Gamma_1$ and $\Gamma_2$ of
			$\Gamma$ agree};
        \breakMed
		\infty\,,
                \quad & \hbox{otherwise}.
\end{cases}
\een
To see this, first consider the case when $\Gamma_1\neq\Gamma_2$;
then \Theorem{bivariateLambda}~(ii) and \Proposition{EntropyBdd1} 
imply that $\Lambda^*(\Gamma)=\infty.$ 
Suppose now that $\Gamma_1=\Gamma_2$.
Then \Proposition{EntropyBdd1} shows that 
$\Lambda^*(\Gamma)$ must equal $H(\Gamma\|\Gamma_1\odot P)$
whenever $\Lambda^*(\Gamma)=\infty.$ And if the marginals
agree and $\Lambda^*(\Gamma)$ is finite, then the identification
follows form \Theorem{bivariateLambda}~(iii).

Finally, an application of the contraction
principle \cite[Theorem~4.2.1]{dembo-zeitouni:book}
implies that the univariate convex dual $\Lambda^*(\nu)$
coincides with $I(\nu)$ in (\ref{eq:1dI}). Simply note
that the $\tau^{W_0}$-topology on the space of probability
measures is Hausdorff, and that the map $\Gamma\mapsto\Gamma_1$
is continuous in that topology.
\qed
 
\Theorem{ldp} strengthens
the ``local'' large deviations of 
\cite{kontoyiannis-meyn:I}
to a full LDP. The assumptions
under which this LDP is proved
are more restrictive that those in 
\cite{kontoyiannis-meyn:I}, but 
apparently they cannot 
be significantly relaxed. 
In particular, the density
assumption of (DV3+)~(ii) cannot be 
removed, as illustrated 
by the
counter-example given in
\cite{dupuis-zeitouni:96}.
This example is of an irreducible, 
aperiodic Markov chain 
with state space $\state=[0,1]$,
satisfying Doeblin's condition.
It can be easily seen that this
Markov chain satisfies
condition (DV3) with 
Lyapunov function 
$V(x)=-\half\log x,$
$x\in[0,1]$, and with 
$W$ given by
$$W(x):=
\begin{cases} 
		2-\log\Big\{\frac{4\sqrt{x}}{1-2x}
		\Big[
		\sqrt{\frac{3}{4}-\frac{x}{2}}
		-\sqrt{\frac{1}{4}+\frac{x}{2}}
		\Big]\Big\}
		\quad & \hbox{for $x\in[0,1/2)$};
        \breakMed
		2-\log(2\sqrt{x})
		\quad & \hbox{for $x\in[1/2,1]$}.
\end{cases}
$$
Taking $\delta=1$,  $C=[0,1]$ and
$b=2$ yields a solution to (DV3),
with the Lyapunov function $V$ and
the {\em unbounded} function $W$ as above.
But for this Markov chain 
the density assumption in (DV3+)~(ii)
is {\em not} satisfied, and as shown in 
\cite{dupuis-zeitouni:96},
it satisfies the LDP with 
a rate function different
from the one in \Theorem{ldp}.


The LDP of \Theorem{ldp} can easily be 
extended to the sequence of empirical 
measures of $k$-tuples $L_{n,k}$, 
defined for each $k\geq 2$ by
\be
L_{n,k} \eqdef 
\frac{1}{n} \sum_{t=0}^{n-1} \delta_{(\Phi(t),\Phi(t+1),\ldots,\Phi(t+k-1))},
\qquad  n\geq 1\, .
\label{eq:Lnk}
\ee
We write $\clM_{1,k}$ for the space of all probability
measures on $(\state^k,\clB^k)$, and we let $\clF_k$
denote the $\sigma$-field of subsets of 
$\clM_{1,k}$ defined analogously to 
$\clF$ in (\ref{eq:clF}), with $\state^k
$ in place of $\state$, and with 
real-valued functions 
$F$ in the space
$$
L_{\infty,k}^{W_0}:=
\Big\{F:\state^k\to\Co\,:\,
\|F\|_{W_0}:=
\sup_{(x_1,\ldots,x_k)\in\state^k}
\Big(
\frac{|F(x_1,\ldots,x_k)|}{W_{0}(x_1)+\cdots+ W_{0}(x_k)}\Big)<\infty
\Big\}$$
instead of $\LWo $. Similarly, the
$\tau_k^{W_0}$-topology on $\clM_{1,k}$ 
is defined by the system of neighborhoods
\ben
N^k_F(c,\delta)\eqdef 
\bigl\{\nu\in\clM_{1,k}:|\nu(F)-c|<\delta\bigr\}\,,
\quad\mbox{for real-valued}\;F\in 
L_{\infty,k}^{W_0},\,c\in\RL,\,\delta>0\,.
\een

A straightforward generalization of the
argument in the above proof yields the
following corollary. The proof 
is omitted.

\begin{corollary}
\tlabel{klpd}
Under the assumptions of \Theorem{ldp},
for any
initial condition $\Phi(0)=x$, the sequence of
empirical measures $\{L_{n,k}\}$ satisfies the LDP
in the space $(\clM_{1,k},\clF_k)$
equipped with the $\tau_k^{W_0}$-topology,
with the good, convex rate function
\ben
I_k(\nu_k)=  
        \left\{ \begin{array}{ll}
	H(\nu_k \|\nu_{k-1}\odot P ),&\;\;\;\mbox{if}\;
                                \nu\;\mbox{is shift-invariant}
\breakMed
                        \infty,                        &\;\;\;\mbox{otherwise.}
                \end{array}
        \right.
\een
where $\nu_{k-1}$ denotes the first $(k-1)$-dimensional
marginal of $\nu_k$.
\end{corollary}

Next we show that under the assumptions of \Theorem{ldp}
it is possible to obtain exact large deviations results
for the partial sums $S_n$,
\be
S_n:=\sum_{t=0}^{n-1} F(\Phi(t))
	=\langle L_n,F\rangle
\,,\qquad n\geq 1,
\label{eq:Sn}
\ee
of a real-valued functional $F\in\LWo$.
In the next two theorems we prove analogs
of the corresponding expansions of Bahadur and Ranga Rao 
for the partial sums of independent random
variables \cite{bahadur-rao:60}. Our results
generalize those obtained by
Miller \cite{miller:61} for finite state 
Markov chains, and those in 
\cite{kontoyiannis-meyn:I} proved 
for geometrically ergodic Markov
processes but only in a neighborhood
of the mean; see 
\cite{kontoyiannis-meyn:I} for further
bibliographical references.


First we note that, since for any $F\in\LWo$
the map $\nu\mapsto\langle\nu,F\rangle$
from $\clM_1$ to $\RL$, is continuous under the $\tau^{W_0}$
topology, we can apply the contraction principle
to obtain an LDP for the partial sums $\{S_n\}$
in (\ref{eq:Sn}): Their laws satisfy the LDP
on $\RL$ with respect to the good, convex rate
function $J(c)$ as in (\ref{eq:rf}),
\ben
J(c)
&=&
	\inf \bigl\{
 	I(\nu) : \hbox{$\nu$ is a probability measure on $(\state,\clB)$
  	satisfying $\nu(F)\geq c$}
                        \bigr\}\\
&=&
	\inf \Big\{
        H(\Gamma\|\Gamma_1\odot P) :
	\hbox{$\Gamma\in\clM_{1,2}$ with marginals
	$\Gamma_1=\Gamma_2$ such that
        $\Gamma_1(F)\geq c$}
			\Big\}.
\een

Alternatively, based on (the weak version of)
the multiplicative mean ergodic theorem in (\ref{eq:weakMMET}),
we can apply the 
G\"{a}rtner-Ellis theorem
\cite[Theorem~2.3.6]{dembo-zeitouni:book} to 
conclude that the laws of the partials sums
$\{S_n\}$ satisfy the LDP
on $\RL$ with respect to the good rate
function $J^*(c),$ 
\be
J^*(c)\eqdef\sup_{a\in\RL}\;[ac-\Lambda(aF)]\,,\qquad c\in\RL\,,
\label{eq:Jstar}
\ee
so that, in particular, $J(c)=J^*(c)$ for all $c$.

Now suppose for simplicity that the function $F$
has zero mean $\pi(F)=0$ and nontrivial central
limit theorem variance $\sigma^2(F)>0$; 
recall the definition of $\sigma^2(F)$ from
\Section{met-summary}.
To evaluate the supremum in (\ref{eq:Jstar}), we
recall from \Lemma{LambdaEqualsXi} 
that $\Lambda(aF)$ is convex in $a\in\RL$,
and since by \Theorem{mainSpectral}
it is also analytic, it is strictly convex.
Therefore, if we define
$$\Fmax
\eqdef\lim_{a\to\infty}
\frac{d}{da}\Lambda(aF)
=\sup_{a\in\RL}
\frac{d}{da}\Lambda(aF),$$
then $J^*(c)=\infty$
for values of $c$ larger
than $\Fmax$, and the probabilities
of the large deviations events $\{S_n\geq nc\}$
decay to zero {\em super}-exponentially fast.

Therefore, from now on we concentrate on the
interesting range of values $0<c<\Fmax$.
Note that, although in the case of independent
and identically distributed random variables
it is easy to identify $\Fmax$ as the right endpoint
of the support of $F$, for Markov chains this
need not be the case, as illustrated by the
following example.
 

\medskip

\noindent{\bf Example.}
Let $\bfPhi=\{\Phi(n):n\geq 0\}$ be a discrete-time 
version of the Ornstein-Uhlenbeck
process in $\RL^2$, with $\Phi(0)=x\in\RL^2$ and,
$$\Phi(n+1)=\left(\begin{array}{c}
		\Phi_1(n+1)\\
		\Phi_2(n+1)
		\end{array}
		\right)=
	\begin{bmatrix}0 & 1\\ -a_{2} &  -a_{1}\end{bmatrix}
	\Phi(n)
		+
		\left(\begin{array}{c}
                0\\
                N(n+1)
                \end{array}
                \right)\,,$$
where $\{N(k)\}$ is a sequence of independent and identically
distributed $N(0,1)$ random variables.   
Let $A$ denote the above 2-by-2 matrix, and 
assume that the roots of
the quadratic equation $z^{2}+a_{1}z + a_{2}=0$ lie   within
the open unit disk in $\Co$.

Note that  there exists $\gamma<1$ and a positive definite  
matrix $P$ satisfying, $A^{\transpose} P A \le \gamma I$.
One may take 
$P=\sum_{0}^{\infty} \gamma^{-k}(A^{k})^{\transpose}A^{k}$, 
where $\gamma<1$ is chosen so that the sum is convergent.

Then $\bfPhi$ satisfies (DV3\upp)~(i) with Lyapunov function $V(x) = 1+\epsilon x^{\transpose}P x$,
and $W=V$, for suitably small $\epsilon>0$ (hence, the drift condition (DV4) also holds).
Condition   (DV3\upp)~(ii)   holds with $T_{0}=2$ since $P^{2}(x,\varble)$
has a Gaussian distribution with full-rank covariance.

Consider the functions
\[
F_{+}(x) = \ind_{\{|x_1|<1\}},\quad F_{0}(x)=x_2-x_1,\quad
F(x)=F_{1}(x)+F_{0}(x),\qquad x=(x_1,x_2)^{\transpose}\in\Re^{2}.
\]
The asymptotic variance of $F_{0}$ is zero, 
and     for any initial condition we have
\[
\sum_{t=0}^{n-1}F(\Phi(t))
=\sum_{t=0}^{n-1}F_{+}(\Phi(t))+ 
[\Phi_2(n-1)-x_1]\,.
\]
We conclude that $\Fmax = (F_{+})_{\hbox{\small max}}=1$, although   $\pi\{F >c\} >0$ for each $c\ge 0$  under the
invariant distribution $\pi$.

\medskip

Recall 
form \Section{met-summary}
the definitions of lattice and
non-lattice functionals.


\begin{theorem}
\tlabel{Bahadur-RaoNL}
{\em (Exact Large Deviations for Non-Lattice Functionals)}
Suppose that $\bfPhi$ satisfies {\em (DV3\upp)} with an unbounded
function $W$, and that $F\in\LWo$ is a real-valued, 
strongly-non-lattice functional, with $\pi(F)=0$ and
$\sigma^2(F)\neq 0$.
Then, for any $0<c<\Fmax$ and all $x\in\state$,
\ben
\Probsub_x\{S_n\geq nc\}
       \;\sim\;
        \frac{\cf_a(x)}{a\sqrt{2\pi n\sigma_a^2}}
        e^{-nJ(c)},
        \quad
        n\to\infty,
\een
where $a>0$ is the unique solution of the
equation $\frac{d}{da}\LA(aF)=c$,
$\sigma^2_a:=\frac{d^2}{da^2}\LA(aF)>0$,
$\cf_a(x)$ is the eigenfunction
constructed in \Theorem{mainSpectral}, 
and $J(c)$ is defined in (\ref{eq:rf}).
A corresponding
result holds for the lower tail.
\end{theorem}

The proof of \Theorem{Bahadur-RaoNL} is identical
to that of the corresponding result
in \cite{kontoyiannis-meyn:I}, 
based on the following simple properties of a
Markov chain satisfying (DV3+). 
We omit properties P5 and P6 since 
they are not needed here.

\paragraph{Properties.} Suppose $\bfPhi$ satisfies
(DV3+) with an unbounded function $W$, and choose
and fix an arbitrary $x\in\state$ and a function
$F\in\LWo$ with zero asymptotic mean $\pi(F)=0$
and nontrivial asymptotic variance $\sigma^2=\sigma^2(F)\neq 0$.
Let $S_n$ denote the partial sums in (\ref{eq:Sn})
and write $m_n(\alpha)$ for the moment generating 
functions
\be
m_n(\alpha):=\Expect_x[\exp(\alpha S_n)]
	=\Expect_x[\exp(\alpha\langle L_n,F\rangle)]
	\,,\qquad n\ge 1,\
        \alpha\in\Co.
\label{eq:mgf}
\ee
The proofs of the following properties are
exactly as those of the corresponding
results in \cite{kontoyiannis-meyn:I},
and are based primarily on the 
multiplicative mean ergodic theorem \Theorem{mainSpectral},
and the Taylor expansion of $\Lambda(F)$ given in 
\Proposition{Lcts}.
Observe that by \Theorem{v5} we have that
the Lyapunov function $V$ in (DV3+) satisfies
$\pi(V^2)<\infty$.
\begin{itemize}
\item[P1.]
        For any $m>0$ there is $\bara>m$, $\baromega>0$ and
	a sequence $\{\epsilon_n\}$ such that
        $$m_n(\alpha)=\exp(n\LA(\alpha F))
        [\cf_\alpha (x)+ |\alpha| \epsilon_n]\, ,
        \quad n\geq 1\, ,$$
        and $|\epsilon_n|\to 0$
        exponentially fast
        as $n\to\infty$,
        uniformly over all $\alpha\in\Omega(\bara,\baromega)$,
        with $\Omega(\bara,\baromega)$ as in \Theorem{mainSpectral}.
\item[P2.]
        If $F$ is strongly non-lattice, then for any
	$m>0$ and any
        $0<\omega_0 <\omega_1< \infty$, there
	is $\bara>m$ and
	a sequence $\{\epsilon_n'\}$
        such that
        $$m_n(\alpha)=\exp(n\LA(a F))\epsilon'_n\, ,
        \quad n\geq 1\,,$$
        and $|\epsilon'_n|\to 0$
        exponentially fast
        as $n\to\infty$,
        uniformly over all
        $\alpha=a+i\omega$ with $|a| \le \bara$
        and $\omega_0\leq|\omega|\leq\omega_1$.
\item[P3.]
        If $F$ is lattice (or almost lattice) with span $h>0$,
        then for any $\epsilon>0$, as $n\to\infty$,
        $$\sup_{
        \epsilon\leq|\omega|\leq 2\pi/h - \epsilon
        }|m_n(i\omega)| \to 0
        \qquad\mbox{exponentially fast.}$$
\item[P4.]
	For any $m>0$ there exist $\bara>m$ and $\baromega>0$ 
	such that the function $\LA(\alpha F)$ is analytic in 
	$\alpha\in\Omega(\baralpha,\baromega)$, and for $\alpha=a\in\RL$
	we have $\LA(a F)|_{a=0}=\frac{d}{da}\LA'(a f)|_{a=0}=0$,
                and $\frac{d^2}{da^2}\LA''(a F)|_{a=0}=\sigma^2>0.$
                Moreover, $\sigma^2_a:=\frac{d^2}{da^2}\LA(aF)$
                is strictly positive for 
                real $a\in[-\bara,\bara]$.
        
\item[P7.]
                For each $m>0$ there exist
		$\bara>m$ and $\baromega>0$
		such that
		the eigenfunction $\cf_\alpha$
                is analytic
                in $\alpha\in\Omega(\bara,\baromega)$,
                it satisfies
                $\cf_\alpha\big|_{\alpha=0} \equiv 1$,
                and it is strictly positive
                for real $\alpha$. Moreover,
                there is some $\baromega_0\in(0,\baromega)$
                such that
                $$\delta(i\omega):=|\log \cf_{i\omega}(x)
                        -i\omega\haF(x)|\leq (\mbox{Const})\omega^2,$$
                for all $|\omega|\leq\baromega_0$,
                where $\haF$ is as in \Theorem{Fergodic}.
\end{itemize}

An analogous asymptotic expansion 
for lattice functionals is given in the next
theorem; again, its proof is omitted as it 
is identical to that of the corresponding result 
in \cite{kontoyiannis-meyn:I}.

\begin{theorem}
\tlabel{Bahadur-Raoul}
{\em (Exact Large Deviations for Lattice Functionals)}
Suppose $\bfPhi$ satisfies {\em (DV3\upp)} with an unbounded 
function $W$, and that $F\in\LWo$ is a real-valued,
lattice functional
with span $h>0$, $\pi(F)=0$ and
$\sigma^2(F)\neq 0$.
Let $\{c_n\}$ be a sequence of real
numbers in $(\epsilon,\infty)$
for some $\epsilon>0$, and assume
(without loss of generality) that,
for each $n$, $c_n$ is in the support of $S_n$.
Then, for all $x\in\state,$
\be
\Probsub_x\{S_n\geq nc_n\}
       \;\sim\;
        \frac{h}
        {(1-e^{-ha_n})\sqrt{2\pi n\LA_n''(a_n)}}
        e^{-nJ_n(c_n)},
        \quad
        n\to\infty,
\label{eq:brL1}
\ee
where $\LA_n(a)$ is the 
log-moment generating function of $S_n$,
$$
\LA_n(a):=\log \Expect_x\Big[e^{a S_n} \Big]
\,,\qquad n\geq 1,\,a\in\RL\,,
$$
each $a_n>0$ is the unique solution of the
equation $\frac{d}{da}\LA_n(a)=c_n$,
and $J_n(c)$ is the convex dual of $\LA_n(a)$,
$$
J_n(c):=\LA_n^*(c):=\sup_{\la\in\RL}[\la c-\LA_n(\la)]
\,,\qquad n\geq 1,\,c\in\RL\,.$$
A corresponding
result holds for the lower tail.
\end{theorem}

Observe that the expansion (\ref{eq:brL1}) in the
lattice case is slightly more general than the one 
in \Theorem{Bahadur-RaoNL}. If the sequence
$\{c_n\}$ converges  to some $c>\epsilon$ as $n\to\infty$,
then, as in \cite{kontoyiannis-meyn:I},
the $a_n$ also converge to some $a>0$, and
\ben
\Probsub_x\{S_n\geq nc_n\}
       \;\sim\;
        \frac{h\cf_a(x)}
        {(1-e^{-ha})\sqrt{2\pi n\sigma_a^2}}
        e^{-nJ(c)},
        \quad
        n\to\infty,
\een
where $\sigma^2_a:=\frac{d^2}{da^2}\LA(aF)$.

 
\section*{Acknowledgments}

The authors would like to express their thanks to 
Amir Dembo and Jamal Najim for several interesting
pointers in the literature, and also to Tom Kurtz, 
Jin Feng and Luc Rey-Bellet for sharing their 
unpublished work.

\appendix

\section*{Appendix}
 \addcontentsline{toc}{part}{Appendix}
 

\section{Drift Conditions and Multiplicative Regularity}

\Lemma{multSpecial} allows us to bound the 
expansive term $b\ind_C(x)$ in condition~(DV3).  We say that a set 
$S\in\clB$ is \textit{multiplicatively-special} (m.-special) if for every  
$A\in\clB^+$ there exists $\eta>0$ such that
\[
\sup_{x\in\state}
\Expect_x\Bigr[\exp\Bigl(\eta \tau_A L_{\tau_A}(S)  \Bigr) \Bigr]<\infty
\,.
\]
\begin{lemma}
\tlabel{multSpecial}
If $\bfPhi$ is $\psi$-irreducible, then
every small set is m.-special.
\end{lemma}


\proof
Let $S$ be a small set, and fix $A\in\clB^+$.  
For a given fixed $T>0$, define the stopping times $\{T_n: n\ge 0\}$
inductively via  $T_0=0$, and 
\[
T_{n+1} =\inf\Big\{t\ge T_n + T : \Phi(t) \in S \Big\},\quad n\ge 0\, .
\]
We consider the sequence of functions,
\[
g_n(x) = 
\Expect_x\Bigl[\exp\Bigl(\eta\int_{[0,  \tau_A\wedge T_n  )} \ind_S
(\Phi(t)) \,  dt \Bigr)\Bigr]  ,\qquad n\ge 1\,,
\]
and we let $B_n=B_n(\eta) = \sup_{x\in\state} g_n(x)$, $n\ge 1$.
Since $S$ is small, there exists $\epsilon>0$,  $T>0$, such that 
$ \Prob_x\{\tau_A>T_1\} \le 1-\epsilon < 1$ for all  $x\in S$.
From the strong Markov property we then have, 
\[
\begin{array}{rcl}  
g_{n+1}(x) &=& 
\Expect_x\Bigl[\exp\Bigl( \eta\int_{[0,  \tau_A\wedge T_{n+1}  )}
				\ind_S (\Phi(t)) \, dt\Bigr)\Bigr]  
\breakMed 
&\le&
e^{\eta T}
\Expect_x\Bigl[ \ind(\tau_A \le T_1 ) \Bigr]
\breakMed 
&&
+ 
 \Expect_x\Bigl[ \exp\Bigl( \eta\int_{[0,  T_{1}  )}  \ind_S (\Phi(t)) \, dt \Bigr)  
\Expect_{\Phi(T_1)}\Bigl[ \exp\Bigl( \eta\int_{[0,  \tau_A\wedge T_n  )}
               \ind_S (\Phi(t)) \, dt \Bigr)  \Bigr] 
                  \ind(\tau_A>T_1 )  \Bigr]
\breakBig
&\le &
e^{\eta T} + (e^{\eta T}\sqrt{B_1(2\eta) \Prob\{\tau_A>T_1\} }) B_n 
\le 
e^{\eta T}  +   (e^{\eta T}\sqrt{B_1(2\eta) (1-\epsilon)  } )B_n\,,
\end{array}
\]
for all $x\in S$,
where the last bound uses Cauchy-Schwartz.

This gives an upper bound for $x\in S$,
and the same bound also 
holds for all $x$ since $g_n(x) \le \sup_{y\in S}g_n(y)$.
Choosing $\eta>0$ so small that 
$\rho \eqdef \Bigl( e^{\eta 2T}{B_1(2\eta) (1-\epsilon)  } \Bigr)^{1/2} <1,$
we see from
induction that $\{B_n\}$ is a bounded sequence, and
$\limsup_{n\to\infty}B_n \le (1-\rho)^{-1} e^{\eta T}$.
\qed

\subsection*{\sc Proof of \Theorem{multReg}.}

Recall that, under (DV3), the stochastic process $(m(t),\clF_t)$
given in \eq expMart/ is a  super-martingale.  That is, for  
any stopping time $\tau$,
\begin{equation}
\Expect_x  [m(\tau)]\le m(0) = v(x),\qquad x\in\state.
\elabel{V5mart1}
\end{equation}
Fix any set $A\in \clB^+$. 
An application of 
\Lemma{multSpecial} implies that there exist constants
$b_1,b_2<\infty$, and $\eta_1>0$ such that for any stopping time $\tau$,
\begin{equation}
\Expect_x\Bigl [
 \exp\Bigl(\sum_{s=0}^{\tau-1} 
             [\eta_1\ind_C(\Phi(s)) - b_1 \ind_A(\Phi(s))] \Bigr)\Bigr]
                        \le \exp(b_2) \, . 
\elabel{V5mart2}
\end{equation}

From \eq V5mart1/, Jensen's inequality, and H\"older's inequality, 
 for all sufficiently small $\eta >0 $,  and all finite $b_3>0$,
\begin{eqnarray*}
\lefteqn{
\Expect_x\Bigl [
 \exp\Bigl( \eta V(\Phi(\tau)) + \eta \sum_{s=0}^{\tau-1} [
         W(\Phi(s)) - b_3 \ind_A(\Phi(s))]  \Bigr)
\Bigr]}\\
        &=& \!\!
\Expect_x\Bigl [
 \exp\Bigl( \eta V(\Phi(\tau)) + \eta\sum_{s=0}^{\tau-1} [
         W(\Phi(s)) - \half b\ind_C(\Phi(s)) ] \Bigr) 
 \exp\Bigl(\eta \sum_{s=0}^{\tau-1} [
       \half b\ind_C(\Phi(s))  -    b_3 \ind_A(\Phi(s)) ] \Bigr)
\Bigr]\\
&\le&
\Expect_x\Bigl [
 \exp\Bigl( 2\eta V(\Phi(\tau))+2\eta\sum_{s=0}^{\tau-1} [
         W(\Phi(s)) -  b\ind_C(\Phi(s)) ] \Bigr) 
\Bigr]^{\half}\\
&&\times\,
\Expect_x\Bigl [
 \exp\Bigl(2\eta \sum_{s=0}^{\tau-1} [
       \half b\ind_C(\Phi(s))  -    b_3 \ind_A(\Phi(s)) ]
               \Bigr) \Bigr]^{\half}\\
&\le&
e^{\eta V(x)} \Expect_x\Bigl [
 \exp\Bigl(2\eta \sum_{s=0}^{\tau-1} [
        b\ind_C(\Phi(s))  -    b_3 \ind_A(\Phi(s)) ]
                \Bigr) \Bigr]^{\half}.
\end{eqnarray*}
Setting $b_3=  b_1\eta_1^{-1} $ we obtain from this and
\eq V5mart2/, for all $\eta< \eta_1(2 b)^{-1}$,
\begin{equation}
\Expect_x\Bigl [
 \exp\Bigl(  \eta V(\Phi(\tau)) + \eta \sum_{s=0}^{\tau-1} [
         W(\Phi(s)) - b_3 \ind_A(\Phi(s)) ] \Bigr)
\Bigr]
\le
v_\eta (x) \exp(2\eta bb_2\eta_1^{-1}),\qquad x\in\state.  
\elabel{bddV5}
\end{equation}
Setting $\tau=\tau_A\wedge m$ for $m\ge 1$, and then letting $m\to\infty$ completes the proof.
\qed

\subsection*{\sc Proof of \Theorem{v5}.}
The construction
of a Lyapunov function $V_*$ follows from the bounds given above, 
beginning with \eq bddV5/ (note however that $W\equiv 1$ under (DV2)).
Assume that the set $A\in\clB^+$ is fixed, with $V$ bounded on $A$. 
We assume moreover that $A$ is small -- this is without loss of generality by
 \cite[Proposition~5.2.4~(ii)]{meyn-tweedie:book}.  Fix $k\ge 0$, and define,
\[
 \sigma_{A}\eqdef   \min \{ i\ge 0 : \Phi(i)\in A \} \,, \qquad \tau\eqdef
 \sigma_A\wedge k.
 \]
Consideration of this stopping time in \eq bddV5/ gives the upper bound, for some $b_1<\infty$, 
\[
\Expect_x\Bigl [ \ind(\sigma_A\ge k)
 \exp (  \eta V(\Phi(k)) + \half \eta k  )
\Bigr]
\le
b_1v_\eta (x) e^{-\half \eta k} \, ,
\qquad x\in\state,\  k\ge 0\, ,
\]
and on summing both sides we obtain
the pair of bounds,
\[
v_\eta(x) \le 
V_*(x)\eqdef
\Expect_x\Bigl [ \sum_{k=0}^{\sigma_A}
 \exp(  \eta V(\Phi(k)) + \half \eta k  )
\Bigr]
\le
\Bigl(\frac{b_1 } {1- e^{-\half \eta} } \Bigr) v_\eta (x)\, ,
\qquad x\in\state\, .
\]
We now demonstrate that this function
satisfies the desired drift condition: We have,
\[
PV_*(x) =
\Expect_x\Bigl [ \sum_{k=1}^{\tau_A}
 \exp (  \eta V(\Phi(k)) + \half \eta k)
\Bigr]
\le
e^{-\half \eta}  V_*(x)  + b' \ind_A(x)\,,
\]
with 
$b' = \big(\frac{b_1} {1- e^{-\half \eta} }\big) \sup_{y\in A} v_\eta (y)$.
This is indeed a version of (V4).
\qed

\begin{proposition}
\tlabel{fellerWu}
Suppose that   $\state$ is $\sigma$-compact and locally compact;   
that $P$  has the   Feller property; 
and that   there exists a sequence of compact sets $\{K_{n }: n\ge 1\}$
satisfying \eq Wu/: For any compact set $K\subset\state$,
$$\sup_{x\in K}\Expect_x[e^{n\tau_{K_n}}]<\infty\,.$$
Then, there exists a solution to the inequality,
\[
 \clH(V)\leq -\half W+b\ind_C
\]
such that $V,W\colon \state\to[1,\infty)$ are continuous,
their sublevel sets are precompact,
$C\in\clB$ is compact, and $b<\infty$.
\end{proposition}

\proof
Let $\{O_{n}: n\ge 1\}$ denote a sequence of open, precompact sets
satisfying  $O_{n}\uparrow \state$, and 
$K_{n}\subset \hbox{closure of}\,  O_{n}\subset O_{n+1}$, $n\ge 1$.
For each $n\ge 1$ we consider a continuous function 
$s_{n}\colon\state\to [0,1]$
satisfying  $s_{n}(x)  =  1 $ for $  x\in O_{n}$, and 
$s_{n}(x)  =  0 $ for  $ x\in O_{n+1}^{c}$.  We then define 
a stopping time $\tau_{n}\ge 1$
through the conditional distributions,
\[
\Prob\{\tau_{n}>n\mid \clF_{n}\} =
\prod_{i=1}^{n}\bigr(1-s_{n}(\Phi(i)) \bigl) \,, \qquad n\ge 1.
\] 
 From the conditions imposed on $s_{n}$ we may conclude
that $\tau_{K_{n}}\ge \tau_{O_{n}}\ge \tau_{n}\ge \tau_{O_{n+1}}$ for each $n\ge 1$.

For $n\ge 1$, $m\ge 1$ we define $V_{n,m}\colon\state\to\Re_{+}$ by,
\[
V_{n,m}(x) \eqdef \log \Expect_x\Bigl[ 
\exp \Bigl(  \sum_{i=0}^{\tau_{n}-1 }  \bigl(n-1\bigr) 
\bigl(1-  s_{m}(\Phi(i))  \bigr)   \Bigr)\Bigr]\,, \qquad x\in\state.
\]
Continuity of this function is established as follows:  First, observe that
under the Feller property we can infer that $\Prob_{x}\{\tau_{n}=k\}$ is
a continuous function of $x\in\state$ for any $k\ge 1$.   
The bound $\tau_{n}\le\tau_{K_{n}}$, $n\ge 1$, combined with
\eq Wu/ then establishes a form of uniform integrability sufficient to
infer the desired continuity.

Moreover, by the dominated convergence theorem we have 
$V_{n,m}(x)\downarrow 0$, $m\to\infty$, for each $x\in\state$. 
Continuity implies that this convergence is uniform on compacta.  
We choose $\{m_n : n\ge 1\}$ so that $V_{n,m_n}(x)\le 1$ on $O_{n+1}$, 
and we define
$V_n = V_{n,m_n}$.  Letting $W_n =  \bigl(n-1\bigr) \bigl(1-  s_{m}  \bigr)$,
we obtain the bound $\clH(V_n)  \le - W_n + 1.$
Let $\{p_n\}\subset\Re_+$ satisfy $\sum_{n\ge 1}p_n =1$,
  $\sum p_n n=\infty$, and define,
\[
W \eqdef 1 + \sum_{n\ge 1} p_n W_n,
\quad
V\eqdef 1 + \sum_{n\ge 1} p_n V_n.
\]
Convexity of $\clH$ then gives, $\clH(V)\le V - W + 1$.
The functions $W$ and $V$ are evidently coercive and continuous.
Hence the desired inequality is obtained with 
$ C=\{x\in\state : W(x) \le \half\}$.
\qed

\section{\boldmath{$v$}-Separable Kernels}

The following result is immediate from the definition \eq spectralRad/.

\begin{lemma}
\tlabel{SpectralRadiusAnalytic}
Suppose that $\{\haP^n:n\in\nat\}$ is  
a positive semigroup, with finite spectral 
radius $\haxi>0$.   Then the inverse $ [Iz-\haP]^{-1}$ admits the   power series representation,
\[
 [Iz-\haP]^{-1} = \sum_{n=0}^{\infty} z^{-n-1} \haP^{n}\, ,\qquad |z|>\haxi\,,
\]
where the sum converges in norm.
\end{lemma}  

\Lemma{SpectralRadius}~(i) is a simple corollary:

\newpage
 
\begin{lemma}
\tlabel{SpectralRadius}
Consider 
a positive semigroup $\{\haP^n:n\in\nat\}$ that is $\psi$-irreducible.
Then:
\begin{MyEnumerate}\item  
The spectral radius $\haxi$ in $\Lv$ of $\{\haP^t\}$
satisfies $\haxi < b_{0}$ for a given $b_{0}<\infty$ if and only if 
there is a $b<b_0$, and a function $v_1\colon \state\to [1,\infty)$ such that   
$v_1$ equivalent to $v$, and 
$\haP v_1 \le bv_1.$

\item
The generalized principal eigenvalue $\halambda$ 
(see \Section{PF})  satisfies $\halambda \le b <\infty$ if and only if 
there is a  measurable function $v_1\colon \state\to (0,\infty)$ such that,
$\haP v_1 \le bv_1.$
\end{MyEnumerate}

\end{lemma}

\proof
Part (ii) is a consequence of \cite[Theorem~5.1]{nummelin:book}.

To see (i), suppose first that $b>\haxi$, and set
$v_1 =b [I b - \haP] ^{-1}v = \sum_{n=0}^\infty b^{-n} \haP^n v.$
Then $v_1\in\Lv$ by 
\Lemma{SpectralRadiusAnalytic}, and $v\le v_1$ by construction.  
Moreover, it is easy to see that $v_1$ satisfies the desired inequality.

Conversely,  if the inequality holds then for any $0<\eta<1$, $n\ge 1$,
\[
(\eta^{-1} b)^{-n-1}\haP^n v_1 \le b^{-1}\eta^{n+1}n v_1,
\]
which shows that 
$\lll[I\eta^{-1} b - \haP]^{-1}\lll_{v_1}\le (1-\eta)^{-1} \eta b^{-1}.$
It follows that $\haxi \le \eta^{-1}b$ since $v$ and $v_1$ are equivalent.
Since $\eta<1$ is arbitrary, this shows that 
$b\ge \haxi$, and completes the proof.
\qed

The following result   will be used below to construct
$v$-separable kernels.

\begin{lemma}
\tlabel{uniSep}
Suppose that $\haP$ is a positive kernel, and 
that there is a measure $\mu\in\Mv$ satisfying
\[
\haP(x,A)\le \mu(A)\,  , \qquad x\in \state,\ A\in\clB\, .
\]
Then $\haP^2$ is $v$-separable.
\end{lemma}

\proof
Consider the bivariate measure,
$\Gamma(dx,dy) = \mu(dx)\haP(x,dy)v(y),$ for $x,y\in\state.$
Under the assumptions of the proposition we have the upper bound,
$\Gamma(dx,dy) \le v(y)\mu(dx)\mu(dy),$
and hence there exists a density $r$ satisfying $r(x,y) \le v(y)$,
$x,y\in\state$, and $\Gamma = r [ \mu\times\mu]$.  It follows that for
any $g\in\Lv$ we have
\[
\haP g\, (x)  = \int r(x,y) g(y) v^{-1}(y)\mu(dy)\, ,\qquad a.e.\ x\in\state\  [\mu].
\]

For a given  $\epsilon>0$ the function  
$r$ can be approximated from below in $L_1(\mu\times\mu)$ by 
the simple functions,
\[
r_\epsilon(x,y) = \sum_{i=1}^N \alpha_i \ind_{A_i}(x) \ind_{B_i}(y) \le r(x,y),
\qquad x,y\in\state,
\]
and
\[
\int\int |r(x,y)-r_\epsilon(x,y)|\mu(dx)\mu(dy) \le \epsilon\, .
\]

We then define
\[
\haP_\epsilon (x,dy) = r_{\epsilon}(x,y)  v^{-1}(y)\mu(dy)\, ,
\quad  
x,y\in\state\, ,
 \]
 and $\haP_{\epsilon 2} \eqdef \haP \haP_\epsilon$.
 The latter kernel may be expressed  $\haP_{\epsilon 2}= \sum s_i \otimes\nu_i$,
 with
 \[
 s_{i}(x)\eqdef
 \alpha_i \haP(x,A_{i})\, ,
 \quad
 \nu_{i}(dy) = \ind_{B_i}(y)v^{-1}(y)\mu(dy)\, ,
\quad  
x,y\in\state\, .
 \]
 We have $s_{i}\in   \Lv$ and $\nu_{i}\in\Mv$ for each $i$.
  
For any $g\in \Lv$, $x\in\state$, we then have,
\ben
\hspace{0.6in}
| \haP_{\epsilon 2}g\, (x) - \haP^2 g\, (x)|
&=& | \haP[\haP_{\epsilon}g - \haP g]\, (x)|\\
&=& \Bigl|\int \haP(x,dy)  
 \Bigl\{ \int [ r_\epsilon(y,z)-r(y,z) ]g(z)v^{-1}(z)\mu(dz)  \Bigr\} \Bigr|\\
&\le& \int \mu(dy)
    \Bigl\{ \int  | r_\epsilon(y,z)-r(y,z) | |g(z)|v^{-1}(z)\mu(dz)  \Bigr\}\\
&\le& 
\|g\|_v 
 \int\int  | r_\epsilon(y,z)-r(y,z) | \mu(dy)\mu(dz) \le \epsilon \|g\|_v\,.
\hspace{0.6in}
\Box
\een

\begin{lemma}
\tlabel{LambdaTightA}
Suppose that  {\em (DV3)}  holds with $W$ unbounded.
Fix  $0 < \eta\le 1$, and consider any measurable function $F$ satisfying 
\begin{equation}\begin{array}{rcl} F^{+}\eqdef \max(F,0) &  \in & \LW;   
\breakMed
\lim_{r\to\infty} \|F^{+}\ind_{C_W(r)^c} \|_{W}&  < &  \delta\eta.
   \end{array}\elabel{Fplus}
\end{equation}
We then have $\lll I_{C_W(r)^{c}} P_{f}\lll_{v_{\eta}}\to 0$, 
exponentially fast, as $r\to\infty$.
\end{lemma}

\proof 
For simplicity we consider only $\eta=1$.
Choosing $r_{0}\ge1$ so that  $\|F^{+}\ind_{C_W(r_0)^c} \|_{W}  = \delta_0<\delta$,
we have,
\[
P_fe^V \le e^{V-(\delta-\delta_0) W +b}\le e^{V-(\delta-\delta_0) r +b}
 \qquad \hbox{on $C_W(r_0)^c$,}
\]
and hence $\lll \ind_{C_W(r)^c} P_f\lll_v   \le   e^{-(\delta-\delta_0) r +b}$ for all $r\ge 1$.
\qed

\begin{lemma} 
\tlabel{eV3sep} 
Suppose that  {\em (DV3+)}  holds with $W$ unbounded.
Fix  $0 < \eta\le 1$, and consider any measurable function $F$ satisfying \eq Fplus/.Then $(P_f)^{2T_0+2}$ is $v_\eta$-separable.
\end{lemma} 
 
\proof 
For simplicity we present the proof only for $\eta=1$.
We define the truncation,
\[
\haP_r \eqdef (\ind_{C_W(r)} P_f)^{T_0+1}\, .
\]
For each $r\ge 1$ we have
\[
\haP_{r}(x,A) \le \beta_{r}'(A)  \eqdef \int_{C_W(r)} \beta_r(dx) P_{f}(x,A)
\, \quad x\in\state,\ A\in\clB.
\] 
It then follows from  \Lemma{uniSep} that 
 the kernel $\haP_r^{2}$ is $v$-separable.

Finally, applying \Lemma{LambdaTightA} we may 
conclude that $\lll (P_f)^{2T_0+2} - \haP_r^{2}\lll_{v}\to 0$, 
$r\to\infty$, which implies that $(P_f)^{2T_0+2} $ is also $v$-separable.
\qed


\subsection*{\sc Proof of \Theorem{gapB}.}
 
(a) $\Rightarrow$ (b).\quad  When (DV3) holds we can 
conclude from \Lemma{LambdaTightA}  that $\lll P - I_{C_W(r)}P\lll_{v_0} \to 0$ as $r\to\infty$.  
It follows that $\lll P^T - I_{C_W(r)}P^T\lll_{v_0}\to 0$ as 
$r\to\infty$ for any $T\ge 1$.  In particular, this holds for $T=T_0$.
Under the separability assumption on $\{ I_{C_W(r)}P^{T_0}: r\ge 1\}$ 
it then follows that $P^{T_0}$ is $v$-separable. 

\medskip

\noindent
(b) $\Rightarrow$ (a).\quad
We first show that each of the sets $\{ C_{v_0}(r): r\ge 1\}$ is small.
Under the assumptions of (b) we may find, for each 
$\epsilon>0$, an integer   $N\ge 1$, functions
$\{ s_i:1\le i\le N\}\subset L_\infty^{v_0}$,
and probability measures
$\{ \nu_i:1\le i\le N\}\subset \clM_1^{v_0}$
such that, with $K= \sum s_i\otimes\nu_i$, 
\begin{equation}
\lll P^{T_0}- K\lll_{v_0} <\epsilon\, .
\elabel{appPTK}
\end{equation}
This gives for any $r\ge 1$,
\[
|1-\sum s_i(x)| 
= |P^{T_0}\One \,(x)- K\One\,(x) |
\le \epsilon v_0(x) \le \epsilon r,\qquad x\in C_{v_0}(r).
\]

Let $A\in\clB$ be a small set with $\nu_i(A^c)< \epsilon$ for each $i$.
From the bound above and  using similar arguments,
\[
\begin{array}{rcl}
P^{T_0} (x,A^c) 
&\le &
	K(x,A^c) + \epsilon v_0(x)
	\breakMed
&\le &
	  \sum_i s_i(x)\nu_i(A^c) + \epsilon v_0(x)
	\breakMed
&\le &
	(1+\epsilon r) \epsilon + \epsilon r,\qquad x\in C_{v_0}(r)\,.
\end{array}
\]
It follows that for any $r\ge 1$, we may find a small set $A(r)$ such that
$P^{T_0}(x,A(r)) \ge \half$, for $x\in C_{v_0}(r).$
It then follows from \cite[Proposition~5.2.4]{meyn-tweedie:book}
that $C_{v_0}(r)$ is small.
 
We now construct a solution to the drift inequality in (DV3).
Using finite approximations as in  \eq appPTK/, we 
may construct, for each $n\ge 1$, an integer $r_n\ge n$ such that
\[
\lll (P I_{C_{r_n}^c})^{T_0} \lll_{v_0}  
\le
\lll P^{T_0} I_{C_{r_n}^c} \lll_{v_0}  \le e^{-2nT_0}.
\]
Since the norm is submultiplicative, this then gives the bound,
\[
\lll (P I_{C_{r_n}^c})^k \lll_{v_0}  
\le b_0 e^{-2nk},\qquad	 k\ge 0\, ,
\]
where $b_0\eqdef (\lll P\lll_{v_0})^{T_0}$.

We then define for each $n\ge 1$,
\[
v_n = I_{C_{r_n}^c} 
\sum_{k=0}^\infty e^{kn} (P I_{C_{r_n}^c})^k v_0 \, .
\]
From the previous bound on
$\lll (P I_{C_{r_n}^c})^k \lll_{v_0}  $ we   have
the pair of bounds,
\begin{equation}
\| v_n \|_{v_0} \le b_0\frac{1}{1-e^{-n}},
\quad \hbox{and}\qquad 
\| P v_n \|_{v_0} \le b_0 \frac{e^{-2n}}{1-e^{-n}}\, .
\elabel{vnBdds}
\end{equation}
Finally, we set 
\[
\begin{array}{rcl}
V &\eqdef& \log\Bigl( 1+\sum_{n=1}^\infty v_n\Bigr)
\breakMed
W & \eqdef&   b\ind_C   -\clH(V),
\end{array}
\]
where  $C=C_v(r)$ for some $r$, and the constants $b$ and $r$ are 
chosen so that 
$W(x)\ge 1$ for all $x\in\state$.   The bounds \eq vnBdds/
together with the lower bound $v_n \ge v_0 e^{n} \ind_{C_{r_n}^c}$
imply that
\[
\lim_{r\to\infty}\inf_{x\in C_v(r)^c} \exp(-\clH(V)) 
=
\lim_{r\to\infty}
\inf_{x\in C_v(r)^c} \frac{e^{V(x)}}{(Pe^V)\, (x)} 
= \infty,
\]
which implies the existence of $r$ and $b$ satisfying these requirements.
\qed

In much of the remainder of the appendix 
we replace (DV3+) with the following more general 
condition:
\begin{equation}\elabel{DV3relax} 
\mbox{\parbox{.85\hsize}{\raggedright\begin{description}
\item[(i)]
The Markov process $\bfPhi$ is $\psi$-irreducible,
        aperiodic, and it satisfies condition~(DV3) with some
	Lyapunov function $V\colon \state \to [1,\infty)$, 
	and an unbounded function
	$W\colon\state\to [1,\infty)$.

\item[(ii)] 
There exists $T_0>0$ such that
 $I_{C_W(r)} P^{T_0}$ is $v$-separable for  
  for each $r< \infty$.
        \end{description}
}} \end{equation}

\Theorem{gapB}  states that this is roughly equivalent to (DV3+) with
an unbounded function $W$. In fact, we do have an analogous upper
bound for $P^{T_0}$:
\begin{lemma}
\tlabel{DV3relax}
Suppose that the conditions of \eq DV3relax/ hold.
Then, for each $r\ge 1$, $\epsilon>0$,         there is a positive  measure $\beta_{r,\epsilon} \in\Mv$ such that
\[
P^{T_0} h\, (x) \le \beta_{r,\epsilon}(h) + \epsilon\| h\|_v\, ,\qquad x\in C_W(r),\ h\in\Lv. 
\] 
\end{lemma}

\proof
We apply the approximation \eq appPTK/ used in 
the proof of \Theorem{gapB}, where 
$\{ s_i:1\le i\le N\}\subset L_\infty^{v_0}$ are non-negative valued,
and $\{ \nu_i:1\le i\le N\}\subset \clM_1^{v_0}$ 
are probability measures. 
We may assume that  the $\{s_i\}$ satisfy the
bound $1=P^{T_0}(x,\state) \ge \sum s_i(x) - 1$, $x\in C_W(r)$, and 
it follows   that we may take
$\beta_{r,\epsilon} =  2\sum_{i=1}^N\nu_i.$
\qed
The following result is proven exactly as 
\Lemma{eV3sep}, using 
\Lemma{DV3relax}.

\begin{lemma} 
\tlabel{eV3sepb} 
Suppose that the conditions of \eq DV3relax/ hold.
Fix  $0 < \eta\le 1$, and consider any $F\in \LW$   
satisfying \eq Fplus/.   Then $(P_f)^{2T_0}$ is $v_\eta$-separable.
\end{lemma}

\section{Properties of \boldmath{$\Lambda$}
and \boldmath{$\Lambda^*$}}
\slabel{app:Lambda}

In this section we obtain  additional properties
of $\Lambda$ and $\Lambda^*$. One of the main goals 
is to establish approximations of  $\Lambda(G)$ through
bounded functions when $G$ is possibly unbounded. 
Similar issues are treated in \cite[Chapter 5]{deuschel-stroock:book}
where a tightness condition is used to provide related approximations.   

\begin{lemma}
\tlabel{SpectralRadiusConvex}
For a $\psi$-irreducible Markov chain:
\begin{MyEnumerate}
\item  The log-generalized principal eigenvalue $\Lambda$ 
is convex on the space of measurable functions 
$F\colon\state\to (-\infty,\infty]$.
\item  The log-spectral radius $\Xi$ is convex 
on the space of measurable functions $F\colon\state\to (-\infty,\infty]$.
\end{MyEnumerate} 
\end{lemma} 

\proof
The proofs of (i) and (ii) are similar, and both proofs are based on 
\Lemma{SpectralRadius}.  We provide a proof of (ii) only.

Fix $F_1,F_2\in\LWo$, $\eta,\theta\in (0,1)$, and
let $b_i = \eta^{-1} \xi(F_i)$, $i=1,2$. 
\Lemma{SpectralRadius} implies that there exists functions $\{v_1,v_2\}$
equivalent to $v$, and satisfying
\[
 \Expect_x\bigl[\exp\bigl( F_i(\Phi(0) )  + V_i(\Phi(1))  
 \bigr)\bigr]
 \eqdef
P_{f_i} v_i\, (x) \le b_i v_i(x),\qquad i=1,2,\ x\in\state.
\]
We then define
\[
 F_\theta = \theta F_1 + (1-\theta)F_2,
 \quad
V_\theta =  \theta V_1 + (1-\theta)V_2,
\]
 so that by H\"older's inequality,
\[
\begin{array}{rcl}
P_{f_\theta} v_\theta\, (x)
& =&
 \Expect_x\bigl[\exp\bigl( \theta [F_1(\Phi(0) )  + V_1(\Phi(1))  ]
 + (1-\theta)[F_2(\Phi(0) )  + V_2(\Phi(1))  ]\bigr)\bigr]
\breakMed
&\le & 
 \Expect_x\bigl[\exp\bigl( F_1(\Phi(0) )  + V_1(\Phi(1))  
 \bigr)\bigr]^\theta
  \Expect_x\bigl[\exp\bigl(F_2(\Phi(0) )  + V_2(\Phi(1))  \bigr)\bigr]^{ 1-\theta}
\breakMed
&\le & 
 b_1^\theta b_2^{1-\theta}v_{\theta}(x)\, ,\qquad x\in\state.
 \end{array}
\]

The function $v_\theta$ is equivalent to $v$.  Consequently, we may apply
\Lemma{SpectralRadius} once more to obtain that
$\xi(F_\theta) 
\le 
 b_1^\theta b_2^{1-\theta}.$
Taking logarithms then gives,
\[
\Xi(F_\theta)\le     \theta \log(b_1) + (1-\theta)\log(b_2)
= \theta \Xi(F_1) + (1-\theta)\Xi(F_2) -\log(\eta).
\]
This completes the proof since $0<\eta<1$ is arbitrary.
\qed


The following result establishes a form of upper semi-continuity for the functional $\Lambda$.

\begin{lemma}\tlabel{usc}
Suppose that $\bfPhi$ is $\psi$-irreducible, and consider a sequence $\{F_n\}$ of measurable, real-valued functions
on $\state$.  Suppose there exists a measurable function
 $F\colon\state\to \Re$  such that  $F_n\uparrow F$,  as $n\uparrow\infty$.   
 Then the corresponding generalized principal eigenvalues converge:
$\Lambda(F_n) \to \Lambda(F),$ as $n\uparrow\infty.$
\end{lemma}

\proof 
It is obvious that $\limsup_{n\to\infty} \Lambda(F_n) \le \Lambda(F)$.  
To complete the proof we establish a bound
on the limit infimum.

Under the assumptions of the proposition we have 
$P_{f_n}^T \ge  P_{f_1}^T$,
for any $T\ge 1$, $n\ge 1$.
It follows that we can find an integer $T_0\ge 1$,
a function $s\colon\state\to[0,1]$, and a probability
$\nu$ on $\clB$ satisfying $\psi(s)>0$ and
\[
P_{f_n}^{T_{0}}\ge s\otimes\nu\, ,\qquad 1\le n\le \infty.
\]

Let $(\cf_{n},\lambda_{n})$ denote the Perron-Frobenius 
eigenfunction and generalized principal eigenvalue for $P_{f_n}$, normalized 
so that $\nu(h_{n})=1$ for each $n$.
For each $n\ge 1$ we have the upper bound,
$P_{f_n} \cf_{n}\le \lambda_{n} \cf_{n}.$
This gives a lower bound on the $\{\cf_{n}\}$:
\[
\cf_{n}\ge \lambda_{n}^{-T_{0}} P_{f_n}^{T_{0}}\cf_{n}
\ge  \lambda_{n}^{-T_{0}} \nu(\cf_{n}) s = \lambda_{n}^{-T_{0}}   s.
\]

Let $h=\liminf_{n\to\infty} \cf_{n}$, 
$\lambda=\liminf_{n\to\infty} \lambda_{n}$.
Then, by Fatou's Lemma,
$P_f h\le \lambda h.$
We also have $\nu(h)\le 1$ by Fatou's Lemma, and the   lower bound
$h\ge \lambda^{-T_{0}}   s$.  It follows from \Lemma{SpectralRadius} that $\Lambda(F) \le \log(\lambda)$.    
\qed  

In applying \Lemma{usc} we typically assume that suitable regularity conditions
hold so that $\Xi(F)= \Lambda(F)$.   Under a finiteness assumption alone    we   obtain
a complementary continuity result for certain classes of decreasing sequences of functions.
One such result is given here:
\begin{lemma}\tlabel{lscLambda} 
Suppose that $\lll P\lll_{v}<\infty$, and that   $F\colon\state\to\Re$ is 
measurable, with $\Xi(F_{+}) < \infty$.  Then,  with $F_{n}\eqdef
\max(F,-n)$ we have, 
$ \Xi(F_n) \downarrow  \Xi(F)$, as $n\uparrow\infty.$
\end{lemma}
 
\proof
This follows immediately from the approximation,
$\lll P_{f_{n}} - P_{f}\lll_{v}\le e^{-n}\lll P\lll_{v},$
$n\ge 1.$
\qed 

To establish a tight approximation for $\Lambda(M)$,
where $M=\log m$ is as in the proof \Theorem{bivariateLambda},
we will approximate $M$ by bounded functions.    

\begin{proposition}
\tlabel{EntropyBdd3} 
Suppose that $\lll P\lll_{v}<\infty$, and that   $F\colon\state\to\Re$ is 
measurable, with $\Xi(F) < \infty$, and $\Lambda(F)=\Xi(F)$.
Then,  there exists a sequence
$\{n_{k}:k\ge 1\}$ such that with 
$F_{k}\eqdef F\ind \{ -n_{k} \le F\le k\}$ we have:
$$
\Lambda(F_{k}) \to \Lambda(F)
\;\;\;\;\mbox{and}\;\;\;\;
\Xi(F_{k}) \to \Xi(F)
\;\;\;\;\mbox{as}\; k\to\infty\,. 
$$
\end{proposition}

\proof
Let $F^{0}_{k}\eqdef F\ind \{  F\le k\}$.  From \Lemma{usc} we 
have $\Lambda(F_{k}^{0})\uparrow \Lambda(F)$,  $ k\to\infty$. 
It follows that we also have  $\Xi(F_{k}^{0})\uparrow \Xi(F)$,  $ k\to\infty$,
since $\Xi$ dominates $\Lambda$.

We now apply \Lemma{lscLambda}:  
For each $k\ge 1$ we may find $n_{k}\ge 1$ such 
that with $F_{k}\eqdef F\ind \{ -n_{k} \le F\le k\}$,
\ben
	\hspace{1.55in}
 \Lambda(F_{k}^{0})   &\le&  \Lambda(F_{k}) \le \Lambda(F_{k}^{0})  + k^{-1},\\
 \Xi(F_{k}^{0})   &\le&  \Xi(F_{k}) \le \Xi(F_{k}^{0})  + k^{-1},\qquad k\ge 1
	\,. 
	\hspace{1.55in}
	\Box
\een

%

The following proposition implies that $\Lambda$ is tight in a strong sense
under (DV3+):
\begin{proposition}
\tlabel{LambdaTight}
Suppose that the conditions of \eq DV3relax/ hold.   Then,
for any increasing sequence of measurable sets   $K_n\uparrow \state$,
  and any $G\in\LWo$,
\[\begin{array}{rrcl}
\hbox{\rm (i)}&\qquad
\displaystyle
\lim_{n\to\infty} \Lambda (G\ind_{K_n^c})  &  = &   0
\breakBig
\hbox{\rm (ii)}&\qquad
\displaystyle \lim_{n\to\infty} \Lambda (G\ind_{K_n}) &  = &     \Lambda(G) \end{array}\]
\end{proposition}

\medskip

The proof is postponed until after the following   lemma.

\begin{lemma}\tlabel{LambdaTightB}
Suppose that the conditions of \eq DV3relax/ hold, 
and consider  any increasing sequence of measurable sets   
$K_n\uparrow \state$, and any $G\in\LWo$.   
Then, on letting $g_{n}= \exp(I_{K_{n}^{c}} G )$, $n\ge 1$,
we have
 \[
 \lll P^{T_{0}} P_{g_n}  - P^{T_{0}+1}\lll_{v}\to 0\, ,\qquad n\to\infty.
 \]
\end{lemma}

\proof 
We may assume without loss of generality that $G\ge 0$.  As usual, we set
$g=e^{G}$.

Under   \eq DV3relax/  we have $\lll P_{g_n}\lll_{v}\le \lll P_{g}\lll_{v}<\infty$, $n\ge 1$. 
Consequently, 
given \Lemma{LambdaTightA}, it is enough to show that
for any $r\ge 1$,
 \[
 \lll I_{C_W(r)} [   P^{T_{0}} P_{g_n}  - P^{T_{0}+1} ]\lll_{v}\to 0\, ,\qquad n\to\infty.
 \]
 To see this, observe that for any $h\in\Lv$, $x\in\state$,
 \[\begin{array}{rcl}
\bigl | I_{C_W(r)} [ P^{T_{0}} P_{g_n} - P^{T_{0}+1} ] h\, (x) \bigr|
&=&
\bigl|I_{C_W(r)} [  P^{T_{0}}I_{K_{n^{c}}}  [P_g - P] ]h\, (x) \bigr|
\breakMed
&\le &
I_{C_W(r)}   P^{T_{0}} I_{K_{n^{c}}} [P_g - P] |h|\, (x)
\breakMed
 &  \le &   
 \|h\|_{v} \lll P_g \lll_{v}
(  I_{C_W(r)} P^{T_{0}}  I_{K_{n}^{c}} ) v \, (x)
  \breakMed  &  \le &   
 \|h\|_{v} \lll P_g \lll_{v} [
 \beta_{r,\epsilon}(\ind_{K_{n}^{c}}  v )   +\epsilon v]\,,
  \end{array}\]
where the measure $ \beta_{r,\epsilon}\in\Mv$ is given
in \Lemma{DV3relax}.
Consequently,
\[
\limsup_{n\to\infty} \lll  I_{C_W(r)} [ P^{T_{0}} P_{g_n} - P^{T_{0}+1} ]\lll_v
\le \epsilon  \lll P_g \lll_{v}.
\]
This proves the result since $\epsilon >0$ is  arbitrary.
\qed

\noindent{\sc Proof of \Proposition{LambdaTight}.}
To see (i), consider any $G\in\LWo$, and any sequence 
of measurable sets  $K_n\uparrow \state$.
We assume without loss of generality that $G\ge 0$.

Fix any $b>1$, and define for $n\ge 1$,
$G_{n} = (T_{0}+1)b \ind_{K_{n}^{c}}  G.$
In view of  \Lemma{LambdaTightB},
given any $\Lambda>0$, we may find $n\ge 1$ such that
the spectral radius of the semigroup generated 
by the kernel $\haP_{n}\eqdef P^{T_{0}} P_{g_n}$
satisfies $\xi_{n} < e^{\Lambda}$.  With $n,\Lambda$ fixed, we then have
for some $b_{n}<\infty$,
$\haP_{n}^{k}v \le b_{n} e^{k\Lambda}v$ for $k\ge 1.$
This has the sample path representation,
\[\Expect_{x}\Bigl[ \exp\Bigl(  \sum_{i=1}^{k}  
G_n(\Phi(   (T_{0 }+1)i  -1 ))  \Bigr)
v\bigl((T_{0 }+1)k)\bigr)  \Bigr]  
\le  b_{n} e^{k\Lambda} v(x), \qquad x\in\state, k\ge 1.
\]
Denote by $h_{0,k}(x)$ the expectation on the left hand side.  We then have, 
for each $j\ge 1$,
\[
h_{j,k}(x)
\eqdef
P^{j}h_{0,k}(x)
\le  b_{n}e^{k\Lambda}  (\lll P\lll_{v})^{j} v(x)\, ,\qquad x\in\state.
\]
Moreover, each of these functions has a sample path representation,
\[
h_{j,k}(x) =\Expect_{x}\Bigl[ \exp\Bigl(  \sum_{i=1}^{k}  G_n(\Phi( j -1+  (T_{0 }+1)i )) \Bigr)
v\bigl(  \Phi (j+(T_{0 }+1)k) \bigr)
\Bigr)  \Bigr]  \, ,\quad x\in\state,\ j\ge 1,\ k\ge 1.
\]
We then obtain  the following bound 
using H\"older's inequality,
\[\begin{array}{rcl}
\lefteqn{\Expect_{x}\Bigl[ \exp\Bigl(  \sum_{i=T_0}^{ (T_{0 }+1)(k+1)-1}  b \ind_{K_{n}^{c}}   (\Phi(i)) G (\Phi(i))  \Bigr)
v\bigl( \Phi( (T_{0 }+1)(k+1) ) \bigr)
\Bigr]}
&& 
\breakMed
&  \le  &
\Bigl(\prod_{j=0}^{T_{0}}   \Expect_{x}\Bigl[ \exp\Bigl(  \sum_{i=1}^{k}  G_n(\Phi( j-1+  (T_{0 }+1)i )) \Bigr)
v\bigl( \Phi( (T_{0 }+1)(k+1)) \bigr)
\Bigr]  
\Bigr)^{(T_{0} +1)^{-1}}
\breakMed
&  \le  &
\lll P\lll_{v}
\Bigl(\prod_{j=0}^{T_{0}}   \Expect_{x}\Bigl[ \exp\Bigl(  \sum_{i=1}^{k}  G_n(\Phi( j-1+  (T_{0 }+1)i )) \Bigr)
v\bigl(  \Phi (j+(T_{0 }+1)k) \bigr)\Bigr]  
\Bigr)^{(T_{0} +1)^{-1}}
\breakMed  &  = &   
\lll P\lll_{v}
\Bigl(\prod_{j=0}^{T_{0}}  h_{j,k}(x)
\Bigr)^{(T_{0} +1)^{-1}}
\breakMed  &  \le  &   
 b_n (\lll P\lll_{v})^{T_{0}+1} e^{k\Lambda} v(x), \qquad x\in\state, k\ge 1.\end{array}\]
We conclude that $\Lambda(\ind_{K_{n}^{c}} G) \le \Lambda/(T_{0}+1)$.
Since $\Lambda>0$ is arbitrary, it follows that $\Lambda(\ind_{K_{n}^{c}} G) \to 0$,
$n\to\infty$.

To see (ii), fix $\theta\in(0,1)$, and obtain the following bound using convexity,
\[\begin{array}{rcl} \Lambda( \theta G)  &  = &   \Lambda\bigl( \theta \ind_{K_{n}} G  +  (1-\theta) 
 \theta (1-\theta)^{-1}\ind_{K_{n}^{c}} G \bigr)
\breakMed  &  \le &   
\theta \Lambda\bigl(  \ind_{K_{n}} G  \bigr) +  (1-\theta) 
\Lambda\bigl(
 \theta (1-\theta)^{-1}\ind_{K_{n}^{c}} G \bigr) \,.
 \end{array}\] 
From (i) we conclude that
 \[ \Lambda( \theta G) \le \theta  \liminf_{n\to\infty} \Lambda\bigl(  \ind_{K_{n}} G  \bigr)  ,\quad 0<\theta<1,
  \] 
which gives $ \Lambda(  G) \le   \liminf_{n\to\infty} \Lambda\bigl(  \ind_{K_{n}} G  \bigr)  $.
    To obtain the reverse inequality we argue similarly:
 \[  \Lambda\bigl(  \ind_{K_{n}} G  \bigr)  
 \le
 \theta \Lambda\bigl(\theta^{-1} G\bigr) + (1-\theta)\Lambda \bigl( - (1-\theta)^{-1}\ind_{K_{n}^{c}} G \bigr) ,
  \]
 which shows that
 \[
 \limsup_{n\to\infty } \Lambda\bigl(  \ind_{K_{n}} G  \bigr)  \le 
 \theta \Lambda\bigl(\theta^{-1} G\bigr) ,\quad 0<\theta<1.
\]
This shows that $ \Lambda\bigl(  \ind_{K_{n}} G  \bigr) \to  \Lambda\bigl(   G  \bigr) $ as claimed.
\qed

 \Proposition{LambdaTight} allows us to   broaden the class of functions for which 
   $\Xi$ is finite-valued.
\begin{proposition}
\tlabel{coercive}
Suppose that the conditions of \eq DV3relax/ hold.
Then, there exists   $W_1\colon\state\to[1,\infty)$
satisfying the following:
\begin{MyEnumerate}
\item  $W_0\in L_\infty^{W_{1}}$, and $W_1\in\LV$; 

\item
$\displaystyle \sup  \{ V(x) :  x\in C_{W_{1}} (r)  \} < \infty$ for each $r\ge 1$;

\item
$\Xi(W_1) < \infty$.\end{MyEnumerate}
If the state space  $\state$
is $\sigma$-compact, then we may assume that $W_1$ is also coercive.

\end{proposition}

\proof
Fix a sequence of  measurable sets
satisfying $K_{n}\uparrow \state$, with $\sup_{x\in K_n} V(x)<\infty$ for each $n$.
 \Proposition{LambdaTight} implies that we may find,  
for each $k\ge 1$, an integer $n_{k}\ge 1$ such that  
$\Xi( 2^{k+1} \ind_{K_{n_{k}}^{c}}  W_0)   \le 1.$
We then define
 \[
W_1 = \Bigl( 1+  \sum_{k=1}^{\infty}\ind_{K_{n_{k}}^{c}} \Bigr)W_0 .
 \]
   
The functional $\Xi$ is convex by 
\Lemma{SpectralRadiusConvex}, which gives the bound,  
 \[
 \Xi(W_1 )
 \le 
 \half \Xi(2W_0)
 +
 \sum_{k=1}^{\infty} 2^{-k-1} \Xi (2^{k+1 }\ind_{K_{n_{k}}^{c}}   W_0) 
 \le
 \half(1+ \Xi(2W_0)) < \infty.
 \]
 To see that $W_{1}\in\LV$ we apply \Lemma{WinV}.
 
 Finally, if $\state$ is $\sigma$-compact, then the $\{K_{n}\}$ may be taken
 to be compact sets, which then implies the coercive property for $W_{1}$.
\qed 

We have the  following useful corollary.   
The proof is routine, given 
\Proposition{coercive} and 
\Proposition{SpectralRadius}~(i);
see also \cite[Theorem~2.4]{balaji-meyn}. 

\begin{lemma}\tlabel{VerySpecial} 
Suppose that the conditions of \eq DV3relax/ hold.
Then, for any $N<\infty$, there exists $r_{0}\ge 1$, $b_{0}<\infty$,
such that with $\tau=\tau_{C_V(r_{0})}$,
\[
 \Expect_x\bigl [\exp\bigl  (N \tau \bigr)    \bigr] < b_{0} e^{V(x)},
 \qquad x\in\state. 
\]
\end{lemma}
 
We now turn to properties of the dual functional 
$\Lambda^*$ defined in \eq LambdaStar/.
The continuity results stated in  \Proposition{LambdaTight} 
lead to the following representation. 

\begin{proposition}
\tlabel{goodconvexdual}
Suppose that the conditions of \eq DV3relax/ hold. 
Let $\Theta$ be a  linear functional on $\LWotwo$ satisfying 
$\Lambda^*(\Theta)<\infty$.
Then $\Theta$ may be represented as,
\[
\langle \Theta, G\rangle = \nu(G)\, ,\qquad G\in \LWo,
\]
where $\nu\in\MWo$ is a probability measure.
\end{proposition} 

\proof
We proceed in several steps, making repeated use of the bound,
\begin{equation}
\langle \Theta, G\rangle \le 
\Lambda^*(\Theta)  +  \Lambda(G)  < \infty \, ,\qquad G\in 
\LWo.\elabel{DualBdd}
\end{equation}

First note that on considering constant functions in \eq DualBdd/ we obtain,
\[
\Lambda^*(\Theta)
\ge  \sup_{c} [\langle \Theta, c\rangle  - \Lambda(c) ]
=  \sup_{c\in\Re}  [\langle \Theta, 1 \rangle  - 1 ]c.
\]
It is clear that finiteness of $\Lambda^*$
implies that   $\langle \Theta, 1\rangle =1$.
Next, consider any $G\colon \state\to \Re_+$ with $G\in\LWo$.
Then, since $\Lambda(cG)\le 0$ for $c\le 0$,
\[
\Lambda^*(\Theta) \ge \sup_{c} [\langle \Theta, cG\rangle  - \Lambda(cG) ]
\ge \sup_{c<0}  \langle \Theta, G\rangle  c.
\]
We conclude that  $\langle \Theta, G\rangle \ge 0$ for $G\ge 0$.

Consider now a set $A\in\clB$ of $\psi$-measure zero.  
Then $\Lambda(c\ind_{A})=0$ for any $c\ge 0$, and we can argue as above using \eq DualBdd/ 
that $\infty>\Lambda^*(\Theta) \ge \sup_{c>0}  \langle \Theta, 
\ind_{A }\rangle  c$, 
which shows that  $\langle \Theta, \ind_{A }\rangle =0$.  

Finally, we demonstrate
that $\Theta$ defines a countably additive set function on $\clB$.
Let  $\{A_i\}\subset \clB$ denote disjoint sets, 
and let $G_n = \sum_{i=n+1}^\infty \ind_{A_i}$.
Then $0\le G_n\le 1$ everywhere, and $G_n\downarrow 0$. 
\Proposition{LambdaTight} implies that   $\Lambda(b G_{n})\to 0$, $n\to\infty$,
for any $b\in \Re$.  Consequently,
\[
\begin{array}{rcl}
 \Lambda^*(\Theta) &  \ge &   \limsup_{n\to\infty}[ 
		\Theta( b G_n) - \Lambda(bG_n) ]
\breakMed
  &  = &
b  \limsup_{n\to\infty} \Theta( G_n)\,.
\end{array}
\]
It follows that $ \limsup_{n\to\infty} \Theta(G_n)=0$, 
which implies that $\Theta$ defines
a countably additive set function, so that $\Theta$ is in fact a probability 
measure.
\qed


More generally, we define $\LA^*$ for 
bivariate probability measures $\Gamma$ not
necessarily in $\MWotwo$ using the 
same definition as in \eq I2/.
Recall from 
\Lemma{marginals} that  the two marginals of $\Gamma$ agree
whenever 
$\Lambda^*(\Gamma)   <\infty$.  \Proposition{EntropyBdd1}  provides further structure.

\begin{proposition}
\tlabel{EntropyBdd1} 
For any probability measure $\Gamma$ 
on $(\state\times\state,\clB\times\clB)$ with 
first and second marginal equal to some $\cpi$,
\be
\Lambda^*(\Gamma)   \le H(\Gamma \mmid \cpi\odot P)\,,
\label{eq:ident1}
\ee
and, moreover,
\be
\Lambda^*(\Gamma) = H(\Gamma \mmid \cpi\odot P) = \infty
\;\;\;\;\mbox{for}\;\; \Gamma\not\in\clM_{1,2}^{W}.
\label{eq:ident2}
\ee
\end{proposition}

\proof 
If we view $W$ as a function on
$\state\times \state$ with $W(x,y)\equiv W(x)$, $x,y\in\state$, then we have the bound,
for all $\epsilon>0$, $n\ge 1$,
\[
\epsilon \langle \Gamma , W\wedge n\rangle
\le
 \Lambda(\epsilon W\wedge n) + \Lambda^{*}(\Gamma)
\le
 \Lambda(\epsilon W ) + \Lambda^{*}(\Gamma).
 \]
\Lemma{eV3sep}  shows that  $\Lambda(\epsilon W)
 <\infty$ for $\epsilon>0$ sufficiently small, and this   gives
(\ref{eq:ident2}).

Define $\cP$ through the decomposition $\Gamma = \cpi\odot\cP$, and
let $\cExpect$ denote the expectation  for the Markov chain
with transition kernel $\cP$.  
We assume that $\cP$ is of the form
\[
\cP(x,dy) = m(x,y) P(x,dy),\qquad x,y\in\state,
\]
and set $M=\log(m)$, since otherwise the relative
entropy is infinite and there is nothing to prove.
We then have, for any $G\in \LWotwo$,
\[
\begin{array}{rcll}
\Lambda(G) 
&=& \lim_{T\to\infty} \frac{1}{T} \log \Bigl( 
\Expect_x 
    \Bigl[ \exp (   T  \langle L_T, G \rangle   ) \Bigr]   \Bigr)
    &
\breakMed
&=& \lim_{T\to\infty} \frac{1}{T} \log \Bigl( 
\cExpect_x 
    \Bigl[ \exp (   T  \langle L_T, G-M \rangle  ) \Bigr]   \Bigr)
    &
\breakMed
&\ge& 
\limsup_{T\to\infty} \frac{1}{T}   
\cExpect_x 
    \Bigl[     T  \langle L_T, G-M \rangle   \Bigr] 
    \quad  &\hbox{(Jensen's inequality)} 
  \breakMed
  &  =&
  \langle \Gamma, G-M\rangle\qquad a.e.\ x\in\state\ [\cpi]\,,
    \quad &\hbox{(mean ergodic theorem for $\cP$)} 
\end{array}
\] 
where the application of the mean ergodic theorem is
justified by   the $f$-norm ergodic 
theorem \cite[Theorem 14.0.1]{meyn-tweedie:book}. 
The integrability conditions required in this result are obtained as follows.
First,  recall  that $\Gamma(|G|)<\infty$ when $\Lambda^{*}(\Gamma)$ is finite
and $G\in\LWo$.
Also, as in the proof that $ H(\Gamma \mmid \cpi\odot P) \ge 0$, 
one can show that $\Gamma( M_{-}) < \infty $,  where 
$M_{-}\eqdef |M\wedge 0|$.   Consequently, 
$(M-G)_{-}$ is $\Gamma$-integrable, which is what 
is required in the mean ergodic theorem.

\archival{Let $S$ denote the set on which $M<0$, and let $beta=\Gamma(S)$.  
Then, by Jensen's,   $\beta^{-1} \Gamma(|M\wedge 0|)  =$
\[
\beta^{-1} \int_{S }\log(m^{-1}) d\Gamma
\le \log 
\beta^{-1} \int_{S } m^{-1} d\Gamma =\log ( \beta^{-1 } \cpi\odot P  \{S\} )
\le \log ( \beta^{-1 } )
\]
}

The above bound may be interpreted as,
\[
H(\Gamma \mmid \cpi\odot P) 
=
\langle \Gamma, M\rangle
\ge 
\langle \Gamma, G\rangle
-
\Lambda(G) .
\]
Taking the supremum over all $G\in\LWotwo$ gives (\ref{eq:ident1}).
\qed

\addcontentsline{toc}{part}{Bibliography}

\bibliographystyle{plain}

 
\end{document}